\documentclass{article}
\usepackage[utf8]{inputenc}

\usepackage{tabu}

\usepackage{fourier} 
\usepackage{array}
\usepackage{makecell}

\usepackage{multirow}

\usepackage{algorithm2e}
\usepackage{xcolor}

\SetAlFnt{\footnotesize}

\SetAlCapSty{xAlCapSty}

\SetCommentSty{xCommentSty}

\SetNlSty{mynlfont}{}{} 

\LinesNumbered

\SetSideCommentRight

\DontPrintSemicolon

\RestyleAlgo{algoruled}

\usepackage[driverfallback=dvipdfm,
bookmarks=true,    
bookmarksopen=true,
bookmarksopenlevel=\maxdimen,
bookmarksdepth=10
]{hyperref}
\hypersetup{
pdfcenterwindow=true, 
pdfpagelabels=true,
pagebackref=true,            
hypertexnames=true,
plainpages=false,
unicode=false,     
pdftoolbar=true,                
    pdfmenubar=true,          
    pdffitwindow=false,         
    pdfstartview={FitH},        
    pdftitle={Accuracy- Source Term Quadrature for 
 Third-Order Edge-Based Discretization},
    pdfauthor={Yi Liu and Hiroaki Nishikawa},    
    pdfsubject={Accuracy-Preserving Source Term Quadrature for 
 Third-Order Edge-Based Discretization},       
    pdfcreator={Hiroaki Nishikawa},   
    pdfproducer={Hiroaki Nishikawa}, 
    baseurl={http://www.hiroakinishikawa.com},
    pdfkeywords={}, 
    pdfnewwindow=true,     
    colorlinks=true,       
    linkcolor=black,       
    citecolor=black,       
    filecolor=black,       
    urlcolor=black,        
  breaklinks=true,
  hyperfigures=true,
  backref=true,
breaklinks=false, 
pdfpagelabels,      
pagebackref,        
hypertexnames=true, 
plainpages=false,   
naturalnames        
}
\usepackage[all]{hypcap}

\usepackage[thinlines]{easytable}

\usepackage{color}      
\usepackage{placeins}
\usepackage{indentfirst}
\usepackage{multirow}

\usepackage{amssymb}
\usepackage{amsmath,latexsym}
\usepackage{url}
\usepackage{enumerate}
\usepackage{graphicx}

\usepackage{booktabs}
 \newcommand{\ra}[1]{\renewcommand{\arraystretch}{#1}}
 
\usepackage{subcaption}
\usepackage{cleveref}
\usepackage{mathtools}

\def\dyad{\! \otimes \!}

\usepackage{blindtext,graphicx}
\usepackage[absolute]{textpos}
\begin{document}

\begin{center}
\begin{textblock}{20}(0.6,1)
{  {\em Preprint accepted in International Journal for Numerical Methods in Fluids (originally submitted in December, 2020). } \\
Download the journal version at \href{ https://doi.org/10.1002/fld.5028}{ \color{blue}  https://doi.org/10.1002/fld.5028} 
 }

\end{textblock}
\end{center}
 

\title{\bf Economically High-Order Unstructured-Grid Methods: Clarification and Efficient FSR Schemes}

\author{
{Hiroaki Nishikawa}\thanks{ {hiro@nianet.org},
 100 Exploration Way, Hampton, VA 23666 USA, Associate Fellow AIAA}\\
  {\normalsize\itshape National Institute of Aerospace,
 Hampton, VA 23666, USA} 
}

\date{\today}
\maketitle

\begin{abstract} 
In this paper, we clarify reconstruction-based discretization schemes for unstructured grids and discuss their economically 
high-order versions, which can achieve high-order accuracy under certain conditions at little extra cost.  
The clarification leads to one of the most economical approaches: the flux-and-solution-reconstruction (FSR) approach, where highly economical schemes can be constructed based on an extended $\kappa$-scheme combined with economical flux reconstruction formulas, achieving up to fifth-order 
accuracy (sixth-order with zero dissipation) when a grid is regular. Various economical FSR schemes are presented and their formal orders of 
accuracy are verified by numerical experiments. 
\end{abstract}

\section{Introduction}
\label{introduction}
 
 
In this paper, we follow the previous papers \cite{Nishikawa_3rdMUSCL:2020IJNMF,Nishikawa_3rdQUICK:2020,Nishikawa_FakeAccuracy:2020} and clarify economically high-order reconstruction-based unstructured-grid methods for practical computational fluid dynamics (CFD) solvers and identity 
some of the most economical schemes that can be easily implemented in an existing unstructured-grid solver to achieve up to fifth-order 
accuracy when a grid is regular. Economically high-order unstructured-grid methods are defined as those designed for unstructured grids 
with at least second-order accuracy but capable of delivering high-order accuracy under certain conditions (e.g., a regular grid). 
A popular approach is a finite-volume-type discretization with high-order solution reconstruction schemes such as 
UMUSCL \cite{burg_umuscl:AIAA2005-4999,yang_harris:AIAAJ2016,Barakos:IJNMF2018,ZhongSheng:CF2020}. These methods are of 
great interest to practical unstructured-grid CFD solvers because they are relatively simple to implement and bring significant improvements, when 
a grid is relatively regular, with much less memory and computing time than genuinely high-order unstructured-grid 
methods (e.g. discontinuous Galerkin methods). However, there exist some confusions over their orders of accuracy; some methods are not 
high-order even on a regular grid and even in one dimension as we have revealed in the previous papers \cite{Nishikawa_3rdMUSCL:2020IJNMF,Nishikawa_3rdQUICK:2020,Nishikawa_FakeAccuracy:2020}. To best take advantage of these economical methods, it is necessary to resolve the confusions and 
correctly understand their underlying principles. The clarification is based upon the previous 
work \cite{Nishikawa_3rdMUSCL:2020IJNMF,Nishikawa_3rdQUICK:2020,Nishikawa_FakeAccuracy:2020}, which will be significantly expanded here for unstructured grids.

 
Our focus is on methods that can be cast in the flux-balance form:  
\begin{eqnarray} 
\frac{1}{V_j} \sum_{ k \in \{ k_j \} } {\Phi}_{jk} |{\bf n}_{jk}| = 0,
\label{intro:flux_balance_form}
\end{eqnarray}
where $j$ is a node/cell in an unstructured grid, $\{ k_j \} $ is a set of neighbors of $j$, $V_j$ is the volume of a control volume around $j$, ${\bf n}_{jk}$ denotes a scaled vector outward normal to the face between $j$ and $k$, ${\Phi}_{jk}$ is a numerical flux projected along the normal direction computed
 with flux and/or solution reconstructions performed with Burg's unstructured-grid extension of Van Leer's $\kappa$-reconstruction scheme \cite{burg_umuscl:AIAA2005-4999} or its higher-order extension with an extra parameter $\kappa_3$ \cite{yang_harris:AIAAJ2016}; see Equations (\ref{umuscl_L2_2d_euler}) and (\ref{umuscl_R2_2d_euler}). As we will discuss in detail, the flux-balance form arises from various 
 conservative discretization approaches, e.g., finite-volume, grid-less, and edge-based methods. Moreover, it is used as a discretization of not only 
 the surface flux integral of a target conservation law 
 but also of a flux divergence or the cell-average of the flux divergence. Methods of the flux-balance form are very efficient, allowing the computation 
 of residuals in a single loop over faces/edges; it includes many practical algorithms such as cell-centered finite-volume methods \cite{WhiteNishikawaBaurle_scitech2020,scFLOW:Aviation2020} and the edge-based discretization \cite{fun3d_website,nishikawa_liu_source_quadrature:jcp2017,katz_sankaran:JSC_DOI,Boris_Jim_NIA2007-08}. These methods are of 
 great interest because of their ability to achieve higher-order accuracy very economically under certain conditions, e.g., on regular grids. 
 However, there exist confusions over their algorithmic details (e.g., the type of numerical solution), which are important and must be 
 clearly understood if one wishes to develop an unstructured-grid scheme that can achieve third- or higher-order accuracy.

The flux-balance discretization is typically considered as a discrete approximation to a flux integral over a control volume, but it can be 
misleading because sometimes it is more accurate as an approximation to the flux divergence at a 
point \cite{nishikawa_liu_source_quadrature:jcp2017}. Also, a confusion arises about the numerical solution: point value or cell average, leading 
to a confusion over the reconstruction parameter $\kappa$: third-order with $\kappa=1/2$ or $\kappa=1/3$. As we have clarified in the previous 
papers \cite{Nishikawa_3rdMUSCL:2020IJNMF,Nishikawa_3rdQUICK:2020,Nishikawa_FakeAccuracy:2020}, the above flux balance scheme 
becomes third-order accurate for a steady problem in one dimension with $\kappa=1/2$ for point-valued numerical solutions or $\kappa=1/3$ for 
cell-averaged solutions. Another confusion is over third-order accuracy on regular quadrilateral/hexahedral grids. If the flux-balance form is 
a finite-volume scheme, then third- or higher-order accuracy cannot be achieved unless a high-order flux quadrature is employed over 
each face, which is costly with multiple flux evaluations per face and rarely implemented in practical solvers. If a flux-balance scheme shows high-order accuracy on such grids without high-order quadrature, then it must be finite-difference, not finite-volume, as we will discuss later. 

Clarification of various flux-balance forms will lead us to the conclusion that a flux-balance scheme must be constructed so as to reduce to a finite-difference scheme on a regular grid if it is desired to achieve higher-order accuracy in multi-dimensions. 
 However, many current schemes in the flux-balance form are based on the MUSCL approach, where the numerical flux is evaluated with high-order reconstructed solutions. As shown in the previous papers \cite{Nishikawa_3rdQUICK:2020,Nishikawa_FakeAccuracy:2020} and also pointed out in Refs.\cite{ZhangZhangShu2011,NLV6_INRIA_report:2008}, such schemes can be second-order at best for nonlinear equations although still bring improvements to complex flow simulations as demonstrated in Refs.\cite{yang_harris:AIAAJ2016,yang_harris:CCP2018,DementRuffin:aiaa2018-1305}. In this paper, we will focus on schemes that can be genuinely high-order on regular grids, which are strongly desired for scale-resolving turbulent-flow simulations requiring highly refined grids (where high-order schemes are more efficient than second-order schemes). As we discussed in the previous paper \cite{Nishikawa_FakeAccuracy:2020} and will also discuss later, genuine high-order accuracy requires flux reconstruction since a scheme must be finite-difference. But it is generally very expensive to perform on unstructured grids especially in three dimensions and thus desired to be avoided for a scheme to be practical. These considerations lead us to efficient flux-balance schemes with a three-parameter family of flux and solution reconstruction (FSR) schemes without direct flux reconstruction, which can achieve up to fifth-order accuracy when a grid is regular. Our focus is on the development of efficient FSR schemes with flux derivatives expressed in terms of solution derivatives, for example, by the chain rule. This is the key development to the work presented here. As we will show, various such economical schemes can be developed for third- and fourth-order accuracy. Fifth-order accuracy is difficult to achieve in general but can be achieved for systems having a set of solution variables in which fluxes are quadratic. 
 
 The objective of the paper is to provide the clarification leading to these efficient FSR schemes and verify their formal orders of accuracy for one- and two-dimensional inviscid test cases. The focus on the FSR schemes is not because they are the most efficient and practical schemes but because these schemes have not been available before. Note also that we do not consider modern high-order methods such as discontinuous Galerkin methods and spectral-difference methods because our focus is again on discrertization methods in the flux-balance form (\ref{intro:flux_balance_form}) that can be easily implemented in existing practical codes. Detailed comparative studies with similar schemes (or modern high-order methods) and applications to practical problems in three dimensions will be reported in a subsequent paper. {\color{black} Note again that we are interested in schemes that can be 
easily implemented in practical unstructured-grid codes, where only a single layer of neighbor information is typically available unlike 
structured-grid codes having easy access to many neighbors along each grid line. Hence, high-order finite-difference-type schemes cannot be directly implemented;  
the FSR schemes are proposed here as high-order finite-difference 
 schemes that can be easily implemented in practical unstructured-grid codes. 
 }

The paper is organized as follows.
In Section 2, we describe a target conservation law. 
In Section 3, we clarify various unstructured-grid methods of the flux-balance type and forms of target equations used for discretization, and conclude that those approximating the differential form of a conservation law can be more efficient than those approximating the integral form. 
In Section 4, we present a new economically high-order FSR schemes. 
In Section 5, we present truncation errors of the FSR schemes.
In Section 6, we present accuracy verification results for the FSR schemes applied to the Burgers equation in one dimension and the Euler equations in one and two dimensions.
In Section 7, we conclude the paper with remarks. 




\section{Target conservation law: the Euler equations}
\label{general:target_cl_euler}

Consider a general conservation law over a control volume $V$: 
\begin{eqnarray}
\int_V \frac{ \partial {\bf u} }{\partial t} \, dV  + \oint_{\partial V}  {\bf f}_n \, d \! s  =  \int_V {\bf s} \, dV ,
\label{general_cl}
\end{eqnarray}
where $t$ denotes the time, ${\bf u}$ is a vector of conservative variables, ${\bf s}$ denotes a source/forcing term, $\partial V$ denotes the control volume 
boundary, ${\bf f}_n = {\cal F} \cdot {\bf n}$ is a projection of a flux tensor ${\cal F}$ along the outward normal $ {\bf n}$ of the boundary, and $d \! s$ is the infinitesimal boundary 
area, i.e., the length in two dimensions. 

Various physical equations can be written as a conservation law. For the purpose of this paper, it suffices to consider the Euler equations:
 \begin{eqnarray}
  {\bf u} =  \left[  \begin{array}{c} 
               \rho       \\ [1ex]
               \rho {\bf v}     \\ [1ex]
               \rho E     
              \end{array} \right], \quad
  {\cal F} =  
 \left[  \begin{array}{c} 
               \rho {\bf v}      \\  [1ex]
               \rho  {\bf v} \dyad  {\bf v}  + p  {\bf I }  \\ [1ex]
               \rho  {\bf v} H 
              \end{array} \right] ,
              \label{Euler_Conservative}
\end{eqnarray}
where $\rho$ is the density, $p$ is the pressure, ${\bf v}$ is the velocity vector,
$ \dyad$ denotes the dyadic product, ${\bf I}$ is the identity matrix, and $\rho E = \rho H-p   = p / (\gamma-1) + \rho {\bf v}^2 /2$ is the specific total energy ($\gamma=1.4$).
In this work, the forcing term ${\bf s}$ will be relevant only for accuracy verification tests with the method of manufactured solutions. 
Also, we consider only two dimensions, but will keep the general form (\ref{general_cl}). 
Therefore, the velocity has only two Cartesian components ${\bf v} = (u,v)^t$, where the superscript indicates transpose (the velocity vector is taken as a column vector). 
The discretization of viscous terms is beyond the scope of this paper; it will require additional considerations and will be discussed elsewhere. 


\section{Forms of Conservation Law and Discretizations}
\label{general:integral}

There exist various approaches to discretizing the conservation law (\ref{general_cl}). Each approach relies on a specific form, integral or
other alternative forms, of the conservation law (\ref{general_cl}). However, a discretization constructed for a certain form 
may turn out to be more accurate for another form. To avoid confusion and clarify the approaches, we will provide brief discussions for 
some popular forms and their discretizations on unstructured grids. As mentioned earlier, we only discuss methods that can be written in 
the flux-balance form (\ref{intro:flux_balance_form}). 

\subsection{Cell-averaged form with cell-averaged solutions: MUSCL}

\begin{figure}[th!]
    \centering
      \hfill  
                \begin{subfigure}[t]{0.32\textwidth}
        \includegraphics[width=\textwidth]{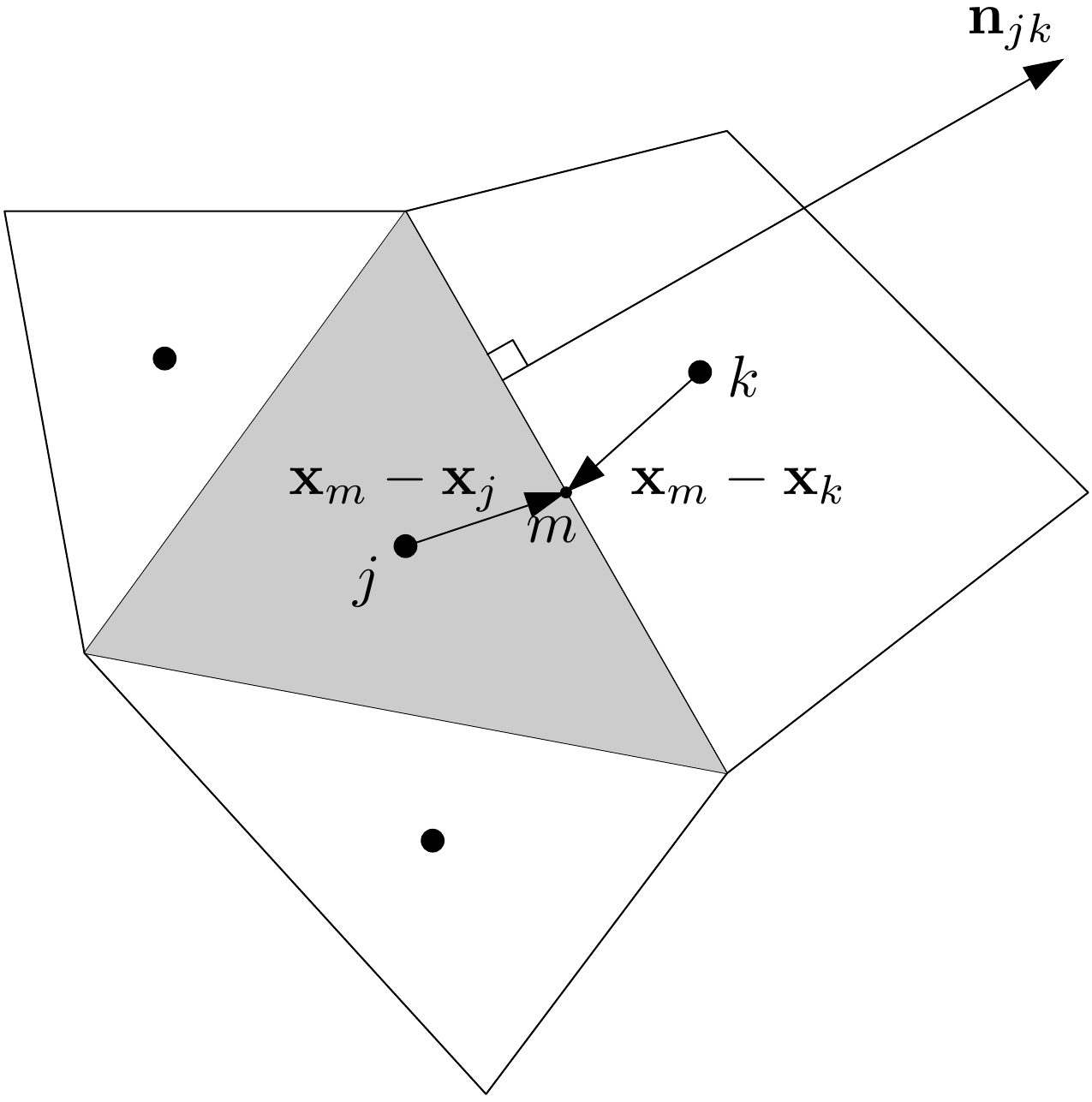}
          \caption{Midpoint rule.}
       \label{fig:stencil_MUSCL}
      \end{subfigure}
      \hfill
          \begin{subfigure}[t]{0.32\textwidth}
        \includegraphics[width=\textwidth]{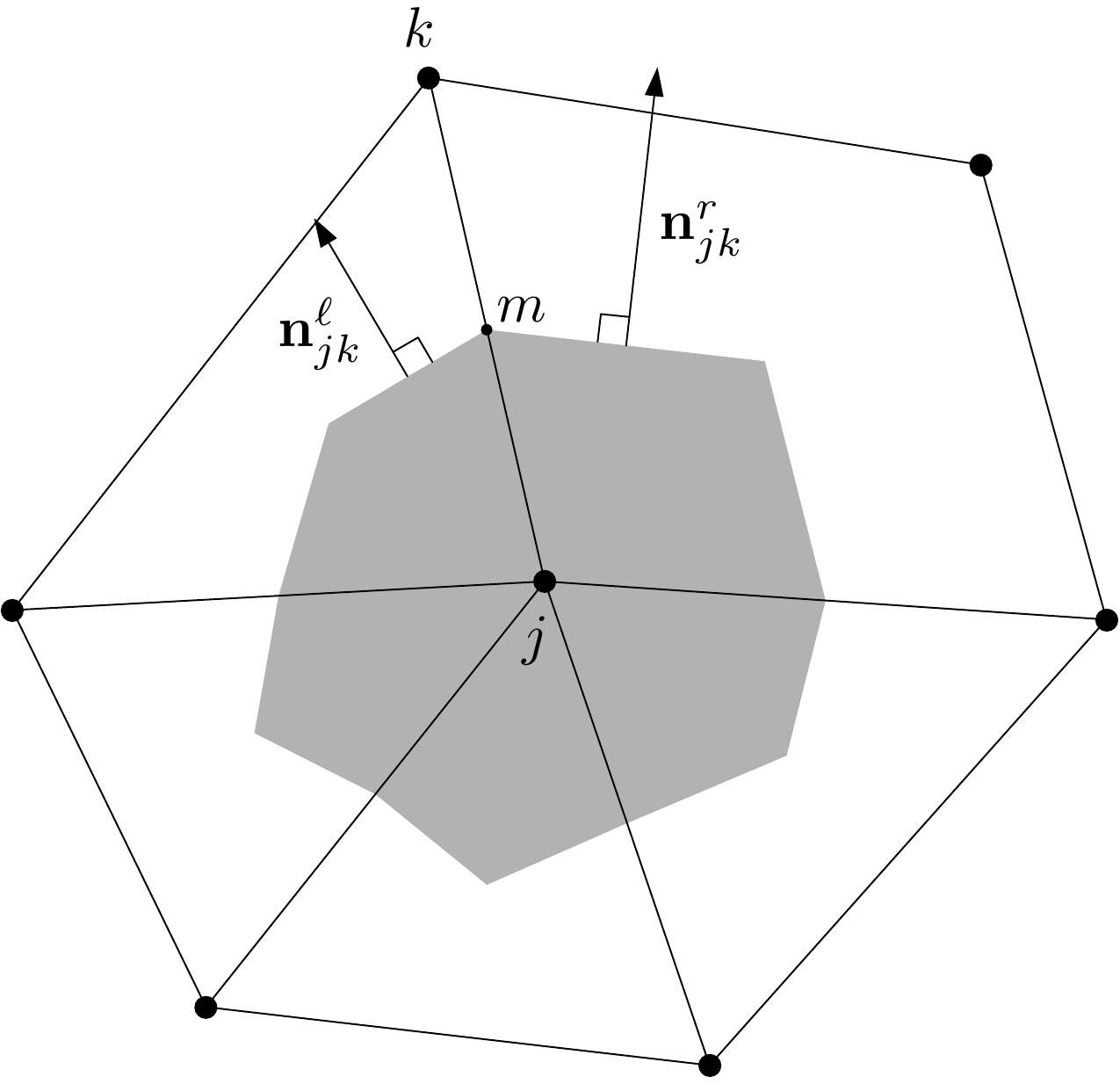}
          \caption{Edge-based quadrature.}
       \label{fig:stencil_EB}
      \end{subfigure}
      \hfill
          \begin{subfigure}[t]{0.28\textwidth}
        \includegraphics[width=\textwidth]{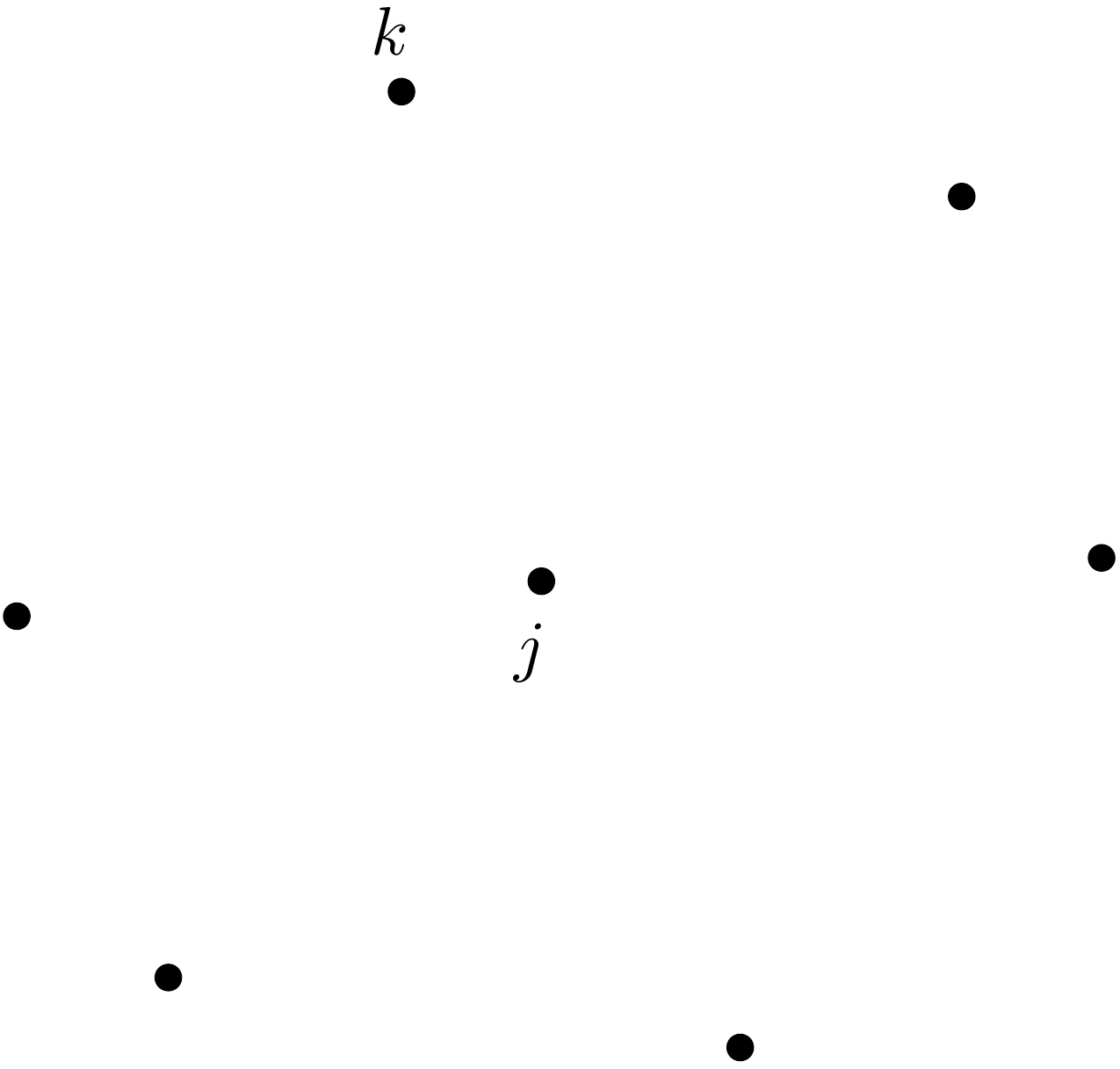}
          \caption{Grid-less.}
       \label{fig:stencil_gridless}
      \end{subfigure}

            \caption{
\label{fig:stenciles_2d}%
Basic stencils for the flux-balance schemes: (a) the midpoint rule for MUSCL-type methods, (b) the edge-based 
flux quadrature, where a flux is computed at the edge-midpoint only, 
for the edge-based schemes (shown here in a triangular-quadrilateral mixed grid), ${\bf n}_{jk} ={\bf n}_{jk}^\ell + {\bf n}_{jk}^r  $; it can be applied to other methods also, (c) grid-less methods. 
} 
\end{figure}

Dividing the conservation law (\ref{general_cl}) by the volume $V$, we obtain the cell-averaged form: 
\begin{eqnarray}
  \frac{ d  \overline{\bf u}}{dt }  + \frac{1}{V}  \int_{\partial V} {\bf f}_n \, ds =   \overline{\bf s},  
\label{cal_form}
\end{eqnarray}
where $ \overline{\bf u}$ and $\overline{\bf s}$ are cell-averaged solution and source term, expressed by the cell-average operator $ {\cal C}$,
\begin{eqnarray}
  \overline{\bf u}  =  {\cal C} {\bf u}    =   \frac{1}{V} \int_{V} {\bf u}  \, dV , \quad 
\overline{\bf s} =  {\cal C} {\bf s}    =   \frac{1}{V} \int_{V} {\bf s}  \, dV  .
\label{ca_operator}
\end{eqnarray}
Equation (\ref{cal_form}) is often called the integral form but here it is referred to as the cell-averaged form in order to distinguish from the original integral conservation law (\ref{general_cl}).
The cell-averaged form is exact and can be directly applied to a finite control volume in a computational grid, as in Figure \ref{fig:stencil_MUSCL}, with the cell-averaged solution 
stored as a numerical solution:
\begin{eqnarray} 
  \frac{ d  \overline{\bf u}_j}{dt } +  \frac{1}{V_j} \sum_{ k \in \{ k_j \} } {\Phi}_{jk} |{\bf n}_{jk}| =  \overline{\bf s}_j,
\label{sec2_muscl}
\end{eqnarray}
where $V_j$ is the volume of the triangle $j$. 
The resulting method is often called a finite-volume method or more specifically the MUSCL finite-volume method \cite{Nishikawa_3rdMUSCL:2020IJNMF}, where the numerical flux is evaluated at a control volume boundary with reconstructed solutions: ${\Phi}_{jk} = {\Phi}_{jk}( {\bf w}_L,{\bf w}_R)$, e.g., via the primitive variables ${\bf w} = (\rho,u,v,p)$,
\begin{eqnarray}
{\bf w}_L \!\! \!\! &=&\!\!\!\!   \kappa  \frac{ \overline{\bf w}_{j} +  \overline{\bf w}_{k}}{2} +   (1-\kappa) \left[  \overline{\bf w}_j   + \widehat\partial_{j} {\bf w}_j    \right]  
 ,     \quad  \widehat\partial_{j}  =  ( {\bf x}_m - {\bf x}_j ) \cdot  \nabla
\label{sec2_umuscl_L2_2d_euler}\\ [1ex]
  {\bf w}_R \!\! \!\!  &=& \!\! \!\!  \kappa  \frac{ \overline{\bf w}_{k} + \overline{\bf w}_{j}}{2} +   (1-\kappa) \left[  \overline{\bf w}_{k} + \widehat\partial_{k} {\bf w}_k   \right] 
, \quad     \widehat\partial_{k}  =  ( {\bf x}_m - {\bf x}_k ) \cdot  \nabla  
, 
 \label{sec2_umuscl_R2_2d_euler}
\end{eqnarray}
where $\nabla {\bf w}_j $ and $\nabla {\bf w}_k$ are gradients computed by a least-squares (LSQ) method, for example, and  ${\bf x}_m$  denotes the position vector of the face centroid (edge midpoint in 2D), and ${\bf x}_j$ and ${\bf x}_k$ denote the nodal position vectors of $j$ and its neighbor $k$, respectively.
This method can be made arbitrarily high-order for general unstructured grids by high-order solution polynomials and high-order discretization of the surface integral (see, e.g., Refs.\cite{jalali_gooch:CF2017,Tsoutsanis-HOFV:JCP2018}). In one dimension, third-order accuracy is achieved with $\kappa=1/3$ as shown in Ref.\cite{Nishikawa_3rdMUSCL:2020IJNMF}.

The high-order flux quadrature requirement significantly increases the complexity of the discretization, especially in three dimensions.
To reduce the cost, Refs.\cite{FVWENO:JSC2014,FVWENO:AMC2016,TamakiImamura:CF2017,Tamaki_PhDThesis2018} proposed economical fourth-order 
discretizations for Cartesian grids, where the flux is evaluated only at the centroid of a face and the solution is reconstructed in a one-dimensional fashion with a flux correction term added
to eliminate low-order quadrature errors. 
It is not clear how these techniques can be extended to unstructured grids (e.g., irregular tetrahedral grids); but economically high-order schemes presented in this paper could provide some insight on their extensions to unstructured grids. 


\subsection{Cell-averaged form with point-valued solutions: UMUSCL-Burg, deconvolution}
\label{cell_averaged_form_with_cell_averaged_sol}

In the discretization of the cell-averaged form, it is indeed natural to choose the cell-average as a numerical solution because the cell-averaged form is an evolution equation of the cell average. 
However, one can choose to store 
a point-valued solution instead as a numerical solution \cite{Nishikawa_3rdQUICK:2020}. Then, the time derivative needs to be expressed in terms the point-valued solution:
\begin{eqnarray}
 {\cal C}  \frac{ d   {\bf u} }{dt }  + \frac{1}{V}  \int_{\partial V} {\bf f}_n \, ds =   \overline{\bf s}, 
\label{ca_form_deconvoluted}
\end{eqnarray}
where the control volume has been assumed to be independent of time $d  \overline{\bf u} /dt  = d ( {\cal C}   {\bf u}) /dt  = {\cal C} d   {\bf u} /dt$, which can be discretized as
\begin{eqnarray} 
 \left[  {\cal C}_{j,j}  \frac{ d   {\bf u}_j}{dt }+  \sum_{ k \in \{ k_j \} } {\cal C}_{j,k}  \frac{ d   {\bf u}_k}{dt} \right] +  \frac{1}{V_j} \sum_{ k \in \{ k_j \} } {\Phi}_{jk} |{\bf n}_{jk}| =  \overline{\bf s}_j,
\label{sec2_quick}
\end{eqnarray}
where ${\Phi}_{jk} = {\Phi}_{jk}( {\bf w}_L,{\bf w}_R)$ and  the operator $ {\cal C} $ has been assumed to be expressed locally with neighbor contributions: 
$ {\cal C}  {\bf u}_j  \approx  {\cal C}_{j,j}  {\bf u}_j +  \sum_{ k \in \{ k_j \} } {\cal C}_{j,k}  {\bf u}_k$ 
 (see Figure \ref{fig:stencil_MUSCL}).
As before, the discretization can be made arbitrarily high-order for general unstructured grids by high-order solution polynomials and high-order discretization of the surface integral. Note that the solution polynomial is here an interpolating polynomial over point-valued solutions stored at cell centers 
 in an unstructured grid. As discussed in Ref.\cite{Nishikawa_3rdQUICK:2020}, this is the reason that $\kappa=1/2$ gives third-order accuracy in one dimension.

To update the point-valued solution in time, we need to invert the cell-average operator ${\cal C}$, but it is not invertible in general: e.g., there 
exist many functions whose cell-average vanishes. One way to overcome the difficulty is to construct an approximate operator $\tilde{\cal C}$ such as 
written in Equation (\ref{sec2_quick}), which can be inverted formally, up to a desired order of accuracy \cite{Denaro:IJNMF1996,Denaro:IJNMF2002}, as
\begin{eqnarray}
\frac{ d   {\bf u} }{dt }  =  \tilde{\cal C}^{-1}  \left[  - \frac{1}{V}  \int_{\partial V} {\bf f}_n \, ds +  {\bf s}  \right]  ,
\label{ca_form_deconvoluted2}
\end{eqnarray}
where the approximate cell-average operator $\tilde{\cal C}$ forms a globally coupled matrix (e.g., as in Equation (\ref{sec2_quick})), which is often called the mass matrix. It needs to be inverted at every time step (or every stage of a multi-stage time-stepping scheme), which is typically 
implemented via a linear solver applied to $\tilde{\cal C} \left(\frac{ d   {\bf u} }{dt } \right)= \left[   - \frac{1}{V}  \int_{\partial V} {\bf f}_n \, ds +  {\bf s}   \right] $. The operator $\tilde{\cal C}^{-1}$ is called a deconvolution operator and 
the methods based on the above form are often called the deconvolution finite-volume methods \cite{Denaro:IJNMF1996,Denaro:IJNMF2002}. See Ref.\cite{Denaro:IJNMF2002} for an example applied to unstructured triangular grids. The approach is similar to continuous Galerkin methods; the similarity has been known for a long time since Ref.\cite{FeliceDenaroMoela:NHT1993}.

The so-called UMUSCL scheme of Burg \cite{burg_umuscl:AIAA2005-4999}, which is called here UMUSCL-Burg, belongs to this category as the numerical solution is stored as point values at nodes and source terms are integrated over a dual control volume. However, the flux integral is approximated by the edge-based quadrature \cite{nishikawa_liu_source_quadrature:jcp2017} (i.e., one-point quadrature; see Figure \ref{fig:stencil_EB}) and therefore it cannot be high-order even on regular grids. Note also that the edge-based flux quadrature loses linear exactness for non-simplex-element grids unless they are regular; therefore even second-order accuracy will be lost on 
such grids \cite{nishikawa:AIAA2010,DiskinThomas:ANM2010,diskin_thomas:AIAA2012-0609}. In one dimension, it is equivalent to the QUICK scheme and therefore achieves third-order accuracy as demonstrated by Burg for a steady problem \cite{burg_umuscl:AIAA2005-4999}. 

For the same reason as in the MUSCL scheme, third-order accuracy cannot be achieved with a single numerical flux per face in multi-dimensions even with the consistent time-derivative treatment as in the deconvolution method. For Cartesian grids, one may be able to develop economical high-order 
methods by deriving high-order solution and flux reconstruction formulas, similar to those in Refs.\cite{FVWENO:JSC2014,FVWENO:AMC2016,TamakiImamura:CF2017,Tamaki_PhDThesis2018}, designed for point-valued solutions instead of cell-averaged solutions.

\subsection{Point-valued differential form: grid-less methods}
\label{cell_averaged_form_with_point_sol}

One may directly discretize the differential form: 
\begin{eqnarray}
\frac{ \partial {\bf u} }{\partial t} + \mbox{div}{\cal F} = {\bf s},
\label{diff_form}
\end{eqnarray}
which is derived from the conservation law (\ref{general_cl}) as a limit of zero control volume. 
In contrast to the cell-averaged form defined over a finite control volume, the differential form is defined at a point. 
For unstructured grids, the so-called grid-less methods \cite{Hong_Meshless:JCP2007,aiaa2009-596,Chiu_etal_meshless_SIAM_2014} and the generalized 
finite-difference methods \cite{HighOrder_GFD:JSC2020} are typical examples of directly discretizing  
the differential form with flux reconstruction at a point with nearby points around it as shown in Figure \ref{fig:stencil_gridless}.
The latter may be considered as a generalized compact scheme, where flux derivatives are computed by solving a globally-coupled
linear system. In these methods, conservation is not automatically satisfied in the discrete sense unless carefully designed as 
such \cite{Hong_Meshless:JCP2007,Chiu_etal_meshless_SIAM_2014}; a conservative
grid-less discretization \cite{Hong_Meshless:JCP2007} may be written in the form:
 \begin{eqnarray} 
   \frac{ d   {\bf u}_j}{dt } +  \frac{1}{V_j} \sum_{ k \in \{ k_j \} } {\Phi}_{jk} |{\bf n}_{jk}| =  {\bf s}_j,
\label{sec2_gridless}
\end{eqnarray}
where the flux balance term is understood as an approximation to the flux divergence at a point $j$ and ${\bf s}_j$ is simply ${\bf s}$ evaluated at a point $j$.  
In contrast to finite-volume methods in the previous sections, these methods may be more efficient especially since no flux quadrature is required. But they 
could actually be more expensive since it will require flux reconstructions for high-order accuracy \cite{HighOrder_GFD:JSC2020}. 
Therefore, these methods can potentially be economical options to unstructured-grid solvers if high-order accurate flux and solution derivatives can be efficiently evaluated, and thus 
should be given further attention. In this paper, we do not discuss these methods as our interest is in methods that can be implemented in existing codes without
introducing new grid data and algorithms (e.g., implicit gradient methods \cite{HighOrder_GFD:JSC2020}).

\subsection{Point-valued differential form with an approximate flux integral}

For ensuring conservation in the discretization of the differential form (\ref{diff_form}), it is convenient to express the flux divergence by a surface integral:
\begin{eqnarray}
\frac{d  {\bf u}}{dt} +  \frac{1}{V} \oint_{\partial V} {\bf f}_n \, ds = {\bf s},
\label{diff_form_approx}
\end{eqnarray}
which is, however, exact only for linearly-varying fluxes \cite{Nishikawa_aiaa2020-1786,nishikawa_centroid:JCP2020} and thus cannot be used to construct high-order discretizations. One would obtain the discretization of the flux-balance form (\ref{sec2_gridless}) if the flux integral term is discretized with a single flux evaluation and high-order solution reconstruction. In general, this type of discretization is second-order at best no matter how accurate the solution reconstruction is. 
Note also that third-order accuracy can be achieved with $\kappa=1/2$ for a steady problem with ${\bf s}=0$ because then it is equivalent to the deconvoluted finite-volume scheme. Ref.\cite{nishikawa_centroid:JCP2020} mentions the possibility of achieving high-order accuracy 
via reconstruction of the function whose cell average is the point-valued flux, but it is in fact true only for uniform grids as we will discuss later in Section \ref{clarification_cfd}.



Some existing unstructured-grid methods based on high-order solution reconstruction 
schemes \cite{yang_harris:AIAAJ2016,yang_harris:CCP2018,DementRuffin:aiaa2018-1305} are 
described as discretizations of Equation (\ref{diff_form_approx}). This is confusing because they present high-order accurate results, which 
contradicts the statement above. As we will explain later, their methods should be understood as finite-difference schemes, which can be high-order accurate. 
However, as we have shown in the previous 
paper \cite{Nishikawa_FakeAccuracy:2020}, their methods cannot be high-order for nonlinear equations and high-order 
results are due to unexpected linearization of the Euler equations. 





\subsection{Integral form with a point-valued solution: third-order edge-based method}

The third-order edge-based (EB3) discretization method \cite{katz_sankaran:JSC_DOI} is a unique method constructed as a discretization of the original integral conservation law: 
\begin{eqnarray}
\int_{V} \frac{ \partial {\bf u} }{\partial t} \, dV  + \oint_{\partial V}  {\bf f}_n \, ds  =  \int_{V} {\bf s} \, dV, 
\label{general_cl2}
\end{eqnarray}
with point-valued solutions stored at nodes and the flux integral discretized by the edge-based flux quadrature (see Figure \ref{fig:stencil_EB}). 
The flux is extrapolated to the edge midpoint but only linearly (second derivatives are not needed), and other terms are discretized by a family of accuracy-preserving quadrature formulas over a dual control volume around a node \cite{nishikawa_liu_source_quadrature:jcp2017,liu_nishikawa_aiaa2017-0738,nishikawa_liu_aiaa2018-4166}, thus resulting in
 \begin{eqnarray} 
  \left[  {\cal M}_{j,j}  \frac{ d   {\bf u}_j}{dt } + {\cal M}'_{j,j}  \nabla \left( \frac{ d   {\bf u}_j}{dt } \right)_{\! j} +  \sum_{ k \in \{ k_j \} } {\cal M}_{j,k}  \frac{ d   {\bf u}_k}{dt} \right] 
 +  \frac{1}{V_j} \sum_{ k \in \{ k_j \} } {\Phi}_{jk} |{\bf n}_{jk}| =
   \left[  {\cal M}_{j,j}   {\bf s}_j + {\cal M}'_{j,j}  \nabla {\bf s}_j  +  \sum_{ k \in \{ k_j \} } {\cal M}_{j,k} {\bf s}_k \right] ,
\label{sec2_edge_based}
\end{eqnarray} 
where ${\bf n}_{jk} ={\bf n}_{jk}^\ell + {\bf n}_{jk}^r $, and ${\cal M}_{j,j}$, ${\cal M}'_{j,j}$, and ${\cal M}_{j,k}$ are the spatial quadrature weights derived in Ref.\cite{nishikawa_liu_source_quadrature:jcp2017} (
for second-order accuracy, one can set ${\cal M}_{j,j}=1$, ${\cal M}'_{j,j}=0$, and ${\cal M}_{j,k}=0$). Time integration can be performed by any time stepping scheme; the mass matrix needs to be inverted at each stage if an explicit multi-stage time-stepping scheme is used \cite{nishikawa_liu_aiaa2018-4166}.
The flux extrapolation can be performed efficiently by computing the flux gradient in terms of the solution gradient using the chain rule \cite{liu_nishikawa_aiaa2016-3969}; we will describe and extend this technique later. It is important to note that it is third-order accurate only with $\kappa=0$ if $\kappa$-reconstruction scheme is employed for the flux reconstruction. 
{\color{black}
Note that the EB3 method is similar to UMUSCL-Burg in that both employ the edge-based flux quadrature but different in that 
        the EB3 method uses direct flux extrapolation while 
        UMUSCL-Burg computes the flux with solutions reconstructed with $\kappa=1/2$, which is the reason that UMUSCL-Burg cannot be high-order accurate
        for nonlinear equations in two and three dimensions \cite{Nishikawa_FakeAccuracy:2020}.
        Because of the edge-based flux quadrature and the flux reconstruction, the EB3 method is
        third-order accurate on arbitrary simplex-element (triangular/tetrahedral) grids, but reduces to first-order accurate on irregular 
        (or second-order accurate on regular) quadrilateral/hexahedral/prismatic/pyramidal grids \cite{nishikawa:AIAA2010,DiskinThomas:ANM2010,diskin_thomas:AIAA2012-0609}.}
 
One would expect that the EB3 method approximates the integral form, but 
quite interestingly, it actually approximates the differential form (\ref{diff_form}) with third-order accuracy: as proved in 
Refs.\cite{nishikawa_liu_source_quadrature:jcp2017,nishikawa_boundary_formula:JCP2014}, the truncation error on a regular triangular/tetrahedral grid is 
given by
\begin{eqnarray}
\frac{ \partial {\bf u} }{\partial t} + \mbox{div}{\cal F} = {\bf s}  + O(h^3),
\end{eqnarray}
showing that it is a third-order approximation to the differential form (\ref{diff_form}). 
{\color{black} Therefore, although it may look similar, once
implemented in a code, to the 
integral-equation-based methods in Sections \ref{cell_averaged_form_with_cell_averaged_sol} and \ref{cell_averaged_form_with_point_sol},
but is fundamentally different in the mechanism by which third-order accuracy is achieved. }

As mentioned earlier, the EB3 method is
 third-order accurate only on arbitrary simplex-element grids, but it thus serves as a very 
economical third-order method for fully adaptive grids with triangular and tetrahedral grids \cite{NishikawaPadway:Aviation2020}. 
Furthermore, the method has been shown to preserve third-order accuracy with linear grids over a curved boundary \cite{nishikawa_boundary_formula:JCP2014}, 
which is another evidence of the point value approximation, not an integral-based approximation such as the MUSCL finite-volume method. 
The method has been demonstrated for inviscid and viscous simulations on three-dimensional unstructured tetrahedral grids \cite{liu_nishikawa_aiaa2016-3969,liu_nishikawa_aiaa2017-0738,nishikawa_liu_aiaa2018-4166,LiuNishikawa_2017-3443}. While third-order accurate on arbitrary 
tetrahedral grids without computing and storing second derivatives is far more efficient than modern high-order methods, it cannot be extended to fourth- and higher-order due to the limited accuracy of the edge-based quadrature formula. In this paper, we seek methods that can achieve fourth- and fifth-order accuracy if a grid happens to be regular.

\subsection{Point-valued differential form with an exact integral}
\label{clarification_cfd}

Finally, we consider another flux-integral form of the differential form (\ref{diff_form}), which is rarely employed in designing unstructured-grid methods but 
is a key to understanding and developing economically high-order unstructured-grid methods. 
Suppose there exists a function ${\cal G}$ whose cell average is the point-valued flux ${\cal F}$:
\begin{eqnarray}
 {\cal F}  = {\cal C} {\cal G}  . 
      \label{cons_diff_00}
\end{eqnarray}
Then, take the divergence,
\begin{eqnarray}
  \mbox{div}{\cal F} = \mbox{div}  \left(  {\cal C}  {\cal G} \right) =  {\cal C} \mbox{div}    {\cal G} - {\cal E}_{comm},
  \label{commute_divf_g}
\end{eqnarray}
where $ {\cal E}_{comm}$ is a commutation error, and thus
\begin{eqnarray}
  \mbox{div}{\cal F}  
  =  \frac{1}{V} \int_{V} \mbox{div}{\cal G}  \, dV  - {\cal E}_{comm}
 =   \frac{1}{V} \oint_{\partial V} {\bf g}_n \, ds - {\cal E}_{comm}, 
\end{eqnarray}
where $   {\bf g}_n  = {\cal G} \cdot {\bf n}  $. 
Substituting this into the differential form (\ref{diff_form}), we obtain
 \begin{eqnarray}
\frac{ \partial {\bf u} }{\partial t}+  \frac{1}{V} \oint_{\partial V} {\bf g}_n \, ds - {\cal E}_{comm} = {\bf s},
\label{diff_form_conservatrive_exact_eps}
\end{eqnarray}
which becomes, if $ {\cal C}$ commutes with the divergence operator: ${\cal E}_{comm} = 0$, 
 \begin{eqnarray}
\frac{ \partial {\bf u} }{\partial t}+  \frac{1}{V} \oint_{\partial V} {\bf g}_n \, ds = {\bf s}.
\label{diff_form_conservatrive_exact}
\end{eqnarray}
This is exact at a point and equivalent to the differential form  (\ref{diff_form}) for a finite control volume (as long as ${\cal E}_{comm} = 0$), and therefore 
a high-order method can be constructed by discretizing the flux integral with high-order quadrature and high-order reconstruction of the 
function ${\cal G}$. 


It is well known that the commutation error may vanish only on regular grids (except some specially stretched grids \cite{Merryman:JSC2003}) and therefore the above form is not exact on a general unstructured grid. 
However, it is a very useful form specifically for our purpose because we are only interested to achieve high-order accuracy when a grid is regular. 
Observe that the relation (\ref{commute_divf_g}) gives, when there is no commutation error,
\begin{eqnarray}
\frac{ \partial   {\bf f}_x }{\partial x}=  \frac{ \partial   ( {\cal C}  {\bf g}_x) }{\partial x}  =   {\cal C} \frac{ \partial   {\bf g}_x }{\partial x} =     \frac{1}{h_x} \int_{x-h_x/2}^{x+h_x/2}  \frac{ \partial   {\bf g}_x }{\partial x}  \,  dx
=   \frac{1 }{h_x} \left[ {\bf g}_x \!\! \left( x+  \frac{h_x}{2} \right)  -{\bf g}_x \!\! \left( x-  \frac{h_x}{2} \right)   \right],
\end{eqnarray}
where ${\bf f}_x$ is the $x$-component of ${\cal F}$ and ${\bf g}_x$ is the $x$-component of ${\cal G}$, and similarly in other coordinate directions. Therefore, Equation (\ref{diff_form_conservatrive_exact}) remains exact when the surface integral is discretized with {\it a one-point flux quadrature formula} applied at each face. In two dimensions, it is written over a rectangular control volume $ [ x- \frac{h_x}{2} , x +  \frac{h_x}{2}] \times[ y- \frac{h_y}{2} , y +  \frac{h_y}{2}] $ around a point $(x,y)$, where $h_x$ and $h_y$ are constants, as 
 \begin{eqnarray}
\frac{ \partial {\bf u} }{\partial t}
  +  \frac{1 }{h_x} \left[ {\bf g}_x \!\! \left( x+  \frac{h_x}{2} ,y \right)  -{\bf g}_x \!\! \left( x-  \frac{h_x}{2} ,y \right)   \right]
+    \frac{1 }{h_y} \left[ {\bf g}_y \!\! \left( x,y+  \frac{h_y}{2} \right)  -{\bf g}_y \!\! \left( x,y-  \frac{h_y}{2} \right)   \right] = {\bf s},
\label{diff_form_conservatrive_exact_fd_form_2d}
\end{eqnarray}
which is exact and equivalent to the differential form (\ref{diff_form}). Therefore, the flux-balance discretization: 
 \begin{eqnarray} 
   \frac{ d   {\bf u}_j}{dt } +  \frac{1}{V_j} \sum_{ k \in \{ k_j \} } {\Phi}_{jk} |{\bf n}_{jk}| =  {\bf s}_j,
\label{sec2_cfd}
\end{eqnarray}
where ${\Phi}_{jk} $ should be an approximation to ${\cal G}$, can achieve high-order accuracy on regular grids with high-order flux reconstruction, which can be performed by a one-dimensional algorithm along the direction from a node (or a cell center) to its neighbor, e.g., the $\kappa$-reconstruction 
scheme applied to the flux. It is noted that the reconstruction problem here is equivalent to that in the MUSCL scheme in the sense that the reconstructed quantity needed at a face is a function whose cell average is the solution/flux stored at a cell/node. Therefore, the value of $\kappa$ to achieve third-order accuracy is $1/3$. Such a finite-difference scheme has been known since 1977 as derived by Van Leer \cite{VLeer_Ultimate_III:JCP1977}.

Several remarks are in order. First, this exact form has been known for a long time as the basis of high-order conservative finite-difference schemes designed for smooth grids \cite{Shu_Osher_Efficient_ENO_II_JCP1989}, where the grid spacing is nearly uniform in each coordinate direction. Here, instead of extending a one-dimensional high-order finite-difference scheme to unstructured grids as is done in the EBR5 scheme \cite{AbalakinBakhvalovKozubskaya:IJNMF2015}, we propose to construct economical high-order methods by directly discretizing the generalized form (\ref{diff_form_conservatrive_exact}) on unstructured grids. 
Second, any unstructured-grid method written in the flux-balance form, no matter how it is constructed, must be a discretization of the generalized form (\ref{diff_form_conservatrive_exact}) if it can achieve high-order accuracy on regular grids. An example of such an scheme is the NLV6 scheme \cite{NLV6_INRIA_report:2008}; the EBR scheme \cite{AbalakinBakhvalovKozubskaya:IJNMF2015} could also be considered as an example. 
Moreover, the cell-centered finite-volume scheme in Ref.\cite{Nishikawa_FANG_AQ:Aviation2020} is also such an example although it is third-order only for linear equations for the reason discussed in Ref.\cite{Nishikawa_FakeAccuracy:2020}. Third, therefore, a flux-balance discretization of Equation (\ref{diff_form_conservatrive_exact}) is the only conservative method, among those considered so far, that can achieve third- and higher-order accuracy when a grid is regular (and quadrilateral or hexahedral). Fourth, the form (\ref{diff_form_conservatrive_exact_fd_form_2d}) is exact only for regular grids but $h_x$ and $h_y$ can be different. Finally, it is noted that methods in Refs.\cite{yang_harris:AIAAJ2016,yang_harris:CCP2018,DementRuffin:aiaa2018-1305} may be considered as examples of the discretization of Equation (\ref{diff_form_conservatrive_exact}) but only for linear equations because they do not perform flux reconstruction (see Ref.\cite{Nishikawa_FakeAccuracy:2020}).

\subsection{Remarks and classification}

\begin{table}[t]
\ra{1.5}
\begin{center}
{
\begin{tabular}{llcccl}\hline\hline 
\multicolumn{1}{l}{Type}                                               &
\multicolumn{1}{l}{Examples}                                               &
\multicolumn{1}{l}{Target form}                                             &
\multicolumn{1}{l}{Solution}                                 &
\multicolumn{1}{l}{Reconst.}    &
\multicolumn{1}{l}{$O(h^3)$: Grid type}                               

  \\ \hline  
CC-SR   &  MUSCL                   &    cell-average  : Eq.(\ref{cal_form}) &    $\overline{u}_j$  & {\bf w}& $\kappa=\frac{1}{3}$: 1D grids  \\
CC-SR   &  MUSCL \cite{FVWENO:JSC2014,TamakiImamura:CF2017}                  &    cell-average  : Eq.(\ref{cal_form}) &    $\overline{u}_j$  & {\bf w} & $\kappa=\frac{1}{3}$: Cartesian grids  \\
CP-SR   &   \makecell[l]{  UMUSCL-Burg  \cite{burg_umuscl:AIAA2005-4999}, \\ DecFV \cite{Denaro:IJNMF2002} }  &    cell-average : Eq.(\ref{ca_form_deconvoluted})    &    ${u}_j$  & {\bf w} & $\kappa=\frac{1}{2}$: 1D grids \\
PP-SR   &  UMUSCL-YH  \cite{yang_harris:AIAAJ2016}           &      point-value : Eq.(\ref{diff_form})  &     ${u}_j$  & {\bf w} & $\kappa=-\frac{1}{6}$ (Linear): 1D grids \\  \hline
PP-FR   &  NLV6 \cite{NLV6_INRIA_report:2008}          &       point-value : Eq.(\ref{diff_form})  &     ${u}_j$  &  $ {\cal F} $ &$\kappa=\frac{1}{3}$: Regular grids \\
PP-FSR &   EBR \cite{AbalakinBakhvalovKozubskaya:IJNMF2015}         &       point-value : Eq.(\ref{diff_form})  &     ${u}_j$  & {\bf w},  ${\cal F} $  & $\kappa=\frac{1}{3}$: Regular grids  \\
PP-FSR &    FSR           &       point-value : Eq.(\ref{diff_form})  &     ${u}_j$  &  $ {\bf w}$, ${\cal F} $  & $\kappa=\frac{1}{3}$: Regular grids  \\
PP-FSR &  EB3  \cite{nishikawa_liu_source_quadrature:jcp2017}         &       point-value : Eq.(\ref{diff_form})  &     ${u}_j$  & $ {\bf w}$, ${\cal F} $   & $\kappa=0$: Simplex grids

  \\  \hline  \hline
\end{tabular}
}
\caption{A classification of flux-balance discretization methods with a single flux evaluation per face. The target form is the form that a method approximates consistently and most accurately, not necessarily the form based on which the discretization is constructed. Those approximating the point-value differential form (\ref{diff_form}) can be centered at a node or a point in a cell. The solution type is indicated by $\overline{u}_j$ for cell averages and ${u}_j$ for point values.  In CP-SR, the numerical solution is typically defined at the geometric centroid of a cell but better choices are available \cite{nishikawa_centroid:JCP2020}. Reconstruction is performed for the solution ${\bf w}$, or $\cal F$, or ${\bf w}$ and $\cal F$. The last column indicates on what type of grid and at what value of $\kappa$ a scheme can achieve at least third-order accuracy. 
}
\label{Tab.classificaiton}
\end{center}
\end{table}

Before we proceed, we classify the discretization methods based on the target equation form, the numerical solution type, the reconstruction type. 
See Table \ref{Tab.classificaiton}, where methods of the flux-balance form are classified. For example, CC-SR refers to methods approximating
the cell-averaged form with cell-averaged solutions and solution reconstruction. MUSCL is one (perhaps the only) example. It can be high-order accurate on one-dimensional grids but cannot be in multi-dimensions unless high-order flux quadrature is employed. Economical versions \cite{FVWENO:JSC2014,TamakiImamura:CF2017} can achieve third- and higher-order accuracy on Cartesian grids with a single flux per face by 
adding a high-order flux correction. CP-SR refers to methods approximating the cell-averaged form with point-valued solutions and solution reconstruction: e.g., UMUSCL-Burg and the deconvolution-based finite-volume methods. 
These methods achieve high-order accuracy only on one-dimensional grids but can be made high-order accurate with high-order flux quadrature and solution reconstruction in multi-dimensions. 

PP-SR refers to methods approximating the point-valued differential form (\ref{diff_form}) with point-valued solutions and solution reconstruction, 
for which UMUSCL-YH is an example. As discussed in Ref.\cite{Nishikawa_FakeAccuracy:2020}, this method cannot be high-order accurate for nonlinear equations even with high-order flux quadrature and solution reconstruction. 
PP-FR is the flux-reconstruction version, which can be high-order accurate; NLV6 \cite{NLV6_INRIA_report:2008} is an example. On the other hand, PP-FSR refers to methods approximating the point-valued differential form (\ref{diff_form}) with point-valued solutions and flux/solution reconstruction. NLV6 \cite{NLV6_INRIA_report:2008} is one example, EBR \cite{AbalakinBakhvalovKozubskaya:IJNMF2015} is another, EB3 \cite{nishikawa_liu_source_quadrature:jcp2017} is yet another, and finally FSR refers to new schemes we will present in the next section. Note that  NLV6, EBR, FSR can be higher-order accurate only on regular grids with symmetric stencils; the order of accuracy is third- or higher but may be lower than the design order, except EBR, which is designed to be equivalent to a one-dimensional scheme along each grid line and thus can preserve the design order of accuracy (away from boundaries). Note finally that EB3 is the only method that can achieve third-order accuracy on arbitrary grids as long as the elements are simplex (triangles/tetrahedra). It is also very economical since the flux reconstruction can be performed by the chain rule \cite{nishikawa_liu_source_quadrature:jcp2017}.

The clarification suggests that economical high-order methods can be constructed as a discretization of the generalized form (\ref{diff_form_conservatrive_exact}) with flux reconstruction. The NLV6 scheme \cite{NLV6_INRIA_report:2008} and the EBR scheme \cite{AbalakinBakhvalovKozubskaya:IJNMF2015} are considered as examples of such methods. However, it is not simple to implement these schemes in existing finite-volume-type unstructured-grid solvers because these schemes require the data of points located in a line extended from each edge to extend reconstruction stencils, which are not usually available in an existing code. Moreover, flux reconstruction is very expensive, especially in three dimensions since it will require to store three flux vectors and their derivatives. In this paper, we explore more efficient and simpler-to-implement schemes.

\section{FSR for unstructured grids}
\label{sec:FSR_unstructured_grids}

In this section, we present the new PP-FSR schemes, which are similar to NLV6 and EBR but more efficient (because direct flux reconstruction is avoided) and simpler to implement in existing unstructured-grid solvers (because no special grid data, e.g., a list of nodes along an edge direction, are needed). For brevity, we will denote the new schemes simply by FSR: e.g., third-, fourth-, and fifth-order methods are denoted by FSR3, FSR4, and FSR5, respectively. Later, in a subsequent section, we will show that these schemes will achieve high-order accuracy on uniform one-dimensional grids, which is sufficient to prove the design orders of accuracy on regular quadrilateral/hexahedral grids because they approximate the flux derivative in each coordinate direction.

\subsection{Discretization}
\label{FSR:discre}

We define FSR as a flux-balance discretization of Equation (\ref{diff_form_conservatrive_exact}) based on the edge-based flux quadrature applied to a dual control volume
around a node $j$ in an unstructured grid (see Figure \ref{fig:stencil_EB}): 
\begin{eqnarray}
\frac{ d {\bf u}_j }{d t}  + \frac{1}{V_j} \sum_{ k \in \{ k_j \} } {\Phi}_{jk} |{\bf n}_{jk}| = {\bf s}_j,
\label{FSR_discretization}
\end{eqnarray}
with the following numerical flux
\begin{eqnarray}
 {\Phi}_{jk} =  \frac{1}{2} \left[   {\bf f}_L   +  {\bf f}_R  \right]  - \frac{1}{2} \hat{\bf D}_n \left[  {\bf u}({\bf w}_R) -  {\bf u}({\bf w}_L) \right],
\label{cfd_numerical_flux_2d_euler}
\end{eqnarray}
 where $ {\bf f} = {\cal F} \cdot \hat{\bf n}_{jk}$ is the flux projected along the face normal ${\bf n}_{jk} = \hat{\bf n}_{jk} |{\bf n}_{jk}|   $, $\hat{\bf D}_n = | \partial {\bf f} / \partial {\bf u}  |$ is the dissipation term evaluated with the Roe averages \cite{Roe_JCP_1981} of ${\bf w}_j$ and ${\bf w}_k$ (instead of the reconstructed values) for robustness and smoothness of the residual \cite{Nishikawa_RobustFluxes:jcp2020}, the solution reconstruction is performed with the primitive variables 
${\bf w} = (\rho,u,v,p)$ unless otherwise stated, and the flux reconstruction is performed for $ {\bf f}_L   $ and $  {\bf f}_R  $ as we will discuss further below. 

The discretization (\ref{FSR_discretization}) with the numerical flux (\ref{cfd_numerical_flux_2d_euler}) is the basis for the EBR scheme \cite{AbalakinBakhvalovKozubskaya:IJNMF2015}. 
The NLV6 scheme \cite{NLV6_INRIA_report:2008} is also similar but
the dissipation term is defined with $ {\bf f}_R- {\bf f}_L$ in place of $ {\bf u}({\bf w}_R) -  {\bf u}({\bf w}_L) $, which is more efficient but less robust 
as pointed out in Ref.\cite{AbalakinBakhvalovKozubskaya:IJNMF2015}. For example, it is difficult to avoid the so-called expansion shock since 
$ {\bf f}_R- {\bf f}_L = 0 $ at a sonic expansion no matter what coefficient it is multiplied by. An entropy fix is simple to implement in the form (\ref{cfd_numerical_flux_2d_euler})
by modifying the eigenvalues in the matrix $ \hat{\bf D}_n $.

\subsection{Solution reconstruction}
\label{FSR:reconst_sol}

For the solution reconstruction, we employ the following extended $\kappa$-reconstruction scheme: 
\begin{eqnarray}
{\bf w}_L \!\! \!\! &=&\!\!\!\!   \kappa  \frac{ {\bf w}_{j} +  {\bf w}_{k}}{2} +   (1-\kappa) \left[  {\bf w}_j   + \widehat\partial_{\! j} {\bf w}_j    \right] 
+   \kappa_3  \,  \widetilde{\partial_{\! j}^3}   {\bf w}
 ,    
\label{umuscl_L2_2d_euler}\\ [1ex]
  {\bf w}_R \!\! \!\!  &=& \!\! \!\!  \kappa  \frac{ {\bf w}_{k} + {\bf w}_{j}}{2} +   (1-\kappa) \left[  {\bf w}_{k} + \widehat\partial_{k} {\bf w}_k   \right]
+   \kappa_3 \, \widetilde{\partial_k^3} {\bf w}
, 
 \label{umuscl_R2_2d_euler}
\end{eqnarray}
where 
\begin{eqnarray}
 \widetilde{\partial_j^3} {\bf w}   = \frac{1}{2}    \left\{   \widehat\partial_{\! j} {\bf w}_k - \widehat\partial_{\! j} {\bf w}_j  \right \}      -  \widehat\partial_{\! j}^2 {\bf w}_j   ,   \quad
\widetilde{\partial_k^3} {\bf w}  =   \frac{1}{2}    \left\{   \widehat\partial_{k} {\bf w}_j - \widehat\partial_{k} {\bf w}_k  \right \}      -  \widehat\partial_{k}^2 {\bf w}_k  , 
\end{eqnarray} 
\begin{eqnarray}
  \widehat\partial_{j}  =  ( {\bf x}_m - {\bf x}_j ) \cdot  \nabla,   \quad 
    \widehat\partial_{k}  =  ( {\bf x}_m - {\bf x}_k ) \cdot  \nabla, 
\end{eqnarray} 
\begin{eqnarray}
  \widehat\partial_{j}^2 =  [ \, ( {\bf x}_m - {\bf x}_j ) \cdot  \nabla \, ]^2   , \quad
     \widehat\partial_{k}^2  = [ \, ( {\bf x}_m - {\bf x}_k ) \cdot  \nabla \, ]^2, 
\end{eqnarray} 
 ${\bf x}_m$  denotes the position vector of the edge midpoint, and ${\bf x}_j$ and ${\bf x}_k$ denote the nodal position vectors of $j$ and its 
 neighbor $k$, respectively (see Figure \ref{fig:stencil_EB}).
These expressions can be slightly simplified for the node-centered method considered here but the notations will be kept general, so that the proposed schemes can be directly applied to cell-centered methods. For cell-centered methods, where the face centroid may not be located half way between 
two adjacent cell centroids, the modification proposed in Ref.\cite{nishikawa_LP_UMUSCL:JCP2020} or a similar technique described in Ref.\cite{Tamaki_PhDThesis2018} should be employed to preserve second-order accuracy on irregular grids.

The first two terms in the above formulas are the $\kappa$-reconstruction scheme generalized for unstructured grids
 by Burg  \cite{burg_umuscl:AIAA2005-4999}, which is a quadratic reconstruction leading to a third-order advection scheme at $\kappa=1/3$. 
It requires solution gradients to be available at nodes; we will compute them by a linear unweighted least-squares method \cite{nishikawa_stencil:JCP2019}. 
The last term with the coefficient $\kappa_3$ is a term approximating a cubic term, leading to a fifth-order advection scheme at $\kappa_3 = -2/3$. 
The reconstruction formulas are similar to those proposed by Yang and Harris \cite{yang_harris:AIAAJ2016}, but the cubic term is different, which they construct by replacing the solution in the first two terms by the solution gradient. It is a useful construction, but
it alters the baseline $\kappa$-reconstruction scheme: a quadratic reconstruction is achieved with $\kappa=-1/6$, not with $\kappa=1/3$ any more. In our reconstruction formulas, $\kappa_3=0$ recovers the baseline $\kappa$-reconstruction scheme leading to a third-order advection scheme at $\kappa= 1/3$.

\subsection{FSR3/FSR4/FSR5: Direct flux reconstruction}
\label{FSR:reconst_flux}

We begin with direct flux reconstruction schemes. These schemes are not very efficient but serve as the basis for developing more efficient schemes.
For third-order accuracy, it suffices to take $\kappa_3=0$ in the solution reconstruction and reconstruct the flux by applying the $\kappa$-reconstruction scheme: 
\begin{eqnarray}
{\bf f}_L &=& \theta \frac{ {\bf f}_j + {\bf f}_k   }{2} +   (1-\theta) \left[  {\bf f}_j   + \widehat\partial_{\! j} {\bf f}_j  \right] ,  
\label{FSR3_fL_2d_euler} \\ [1ex]
{\bf f}_R &=& \theta  \frac{{\bf f}_k  + {\bf f}_j  }{2} +   (1-\theta) \left[  {\bf f}_k +  \widehat\partial_{k}  {\bf f}_k  \right], 
 \label{FSR3_fR_2d_euler}
\end{eqnarray}
where ${\bf f} = {\cal F} \cdot \hat{\bf n}_{jk}$, ${\bf f}_j={\bf f}({\bf w}_j)$, ${\bf f}_k={\bf f}({\bf w}_k)$, and the parameter $\kappa$ has been 
replaced by $\theta$ for the flux. As we will show later, $\theta=1/3$ gives third-order accuracy. 
This third-order scheme (with $\kappa_3=0$) will be referred to as FSR3 in the rest of the paper. 
FSR3 requires computation and storage of the flux gradients at all nodes, which can be computed by a linear unweighted LSQ method. 
Note that the solution also needs to be reconstructed to evaluate the dissipation term in the numerical flux (\ref{cfd_numerical_flux_2d_euler}).

As we will show later, the leading third-order truncation error of the FSR3 scheme comes from the dissipation term.
Therefore, to achieve fourth-order accuracy, we only need to take $\kappa_3 = \kappa-1$ to make dissipation $O(h^4)$.
This requires the second derivatives of the solution variables. 
The resulting fourth-order scheme will be referred to as FSR4. 
 
If, in addition, we compute and store the second derivatives of the flux, we can achieve fifth-order accuracy by adding a cubic term to FSR3 in a similar manner as done in the solution reconstruction:
\begin{eqnarray}
{\bf f}_L  \!\!\!\! &=& \!\! \!\! \theta  \frac{{\bf f}_j  + {\bf f}_k  }{2} +   ( 1  \!   -   \! \theta) \left[  {\bf f}_j +  \widehat\partial_{\! j} {\bf f}_j  \  \right] 
+ \theta_3 \, \widetilde{\partial_{\! j}^3}  {\bf f}_j
,  
\label{FSR5_fL_2d_euler}\\ [1ex]
{\bf f}_R   \!\!\!\! &=&  \!\!\!\! \theta  \frac{ {\bf f}_k +{\bf f}_j  }{2} +   (1  \! -  \! \theta) \left[  {\bf f}_k +  \widehat\partial_{k}  {\bf f}_k   \right]
+  \theta_3 \, \widetilde{\partial_k^3}  {\bf f}_k ,
 \label{FSR5_fR_2d_euler}
\end{eqnarray}
where $\theta_3$ is another parameter, which we have to set $\theta_3 = - 8 /15 $ to achieve fifth-order accuracy as will be shown later.
This fifth-order scheme will be referred to as FSR5. It is an expensive scheme especially in three dimensions, where three derivatives 
and six second derivatives need to be computed and stored for the solution as well as three flux vectors. 
FSR5 is similar to the NLV6 scheme \cite{NLV6_INRIA_report:2008} but it is different in that we compute the dissipation term with reconstructed solutions 
and also that we compute the second derivatives by successive application of the linear LSQ method following Yang and Harris \cite{yang_harris:AIAAJ2016}. 
In fact, the successive application is a key to high-order accuracy as we will discuss later.

\subsection{CFSR3/CFSR4/CFSR5: Economical flux reconstruction with the chain rule}
\label{FSR:reconst_flux_quadratic_CR}

It is reasonable to ask if the flux derivatives can be computed by the chain rule without degrading the order of accuracy.
It is certainly possible for third-order accuracy as we already have shown in the previous paper \cite{Nishikawa_FakeAccuracy:2020}.
The chain-rule-based flux reconstruction has also been demonstrated for the third-order edge-based method \cite{nishikawa_liu_source_quadrature:jcp2017}.
For the FSR5 scheme, unfortunately, the chain rule does not allow us to keep fifth-order accuracy. However, a very efficient fourth-order scheme 
can be constructed. 

As shown in Ref.{\cite{Nishikawa_FakeAccuracy:2020}, the flux derivatives can be computed as
\begin{eqnarray}
 \widehat\partial_{\! j} {\bf f}_j  =   \left(  \frac{\partial {\bf f}}{\partial {\bf w}} \right)_{\!\! j}   \widehat\partial_{\! j} {\bf w}_j    
, \quad
 \widehat\partial_{k} {\bf f}_k =   \left(  \frac{\partial {\bf f}}{\partial {\bf w}} \right)_{\!\! k}  \widehat\partial_{k} {\bf w}_k  ,
 \label{chain_rule_grad}
\end{eqnarray}
where the flux Jacobian is given, 
with the unit face normal vector $\hat{\bf n}_{jk}$ and ${\bf v}$ defined as column vectors and the notation $u_n = {\bf v} \cdot  {\bf n}_{jk}$, as
\begin{eqnarray}
 \frac{\partial {\bf f}}{\partial {\bf w}}
 =
  \left[  \begin{array}{ccc} 
u_n          &  \rho  \hat{\bf n}_{jk}^t   & 0    \\ [1ex]
u_n {\bf v} & \rho ( u_n {\bf I} + {\bf v} \dyad   \hat{\bf n}_{jk} )  &  \hat{\bf n}_{jk}    \\ [1ex]
u_n {\bf v}^2 /2&  \rho(H  \hat{\bf n}_{jk}^t + u_n {\bf v}^t ) & \gamma u_n /(\gamma-1)
              \end{array} \right].
\end{eqnarray}
FSR3 with these chain-rule flux derivatives will be referred to as CFSR3 ($\kappa_3 = 0$). Similarly, 
FSR4 with the chain-rule flux derivatives will be referred to as CFSR4 ($\kappa_3 = \kappa-1$). As we will show later, CFSR3 and CFSR4 
are third- and fourth-order accurate.

To eliminate flux derivatives from FSR5, we also have to compute the second derivatives by the chain rule: 
\begin{eqnarray}
 \widehat\partial_{\! j}^2 {\bf f}_j  =
  \widehat\partial_{\! j} \left(  \frac{\partial {\bf f}}{\partial {\bf w} } \right)_{\!\! j} 
 +   \left(  \frac{\partial {\bf f}}{\partial {\bf w}} \right)_{\!\! j}   \widehat\partial_{\! j}^2 {\bf w}_j    
, \quad
 \widehat\partial_{k}^2 {\bf f}_k  = 
   \widehat\partial_k \left(  \frac{\partial {\bf f}}{\partial {\bf w} } \right)_{\!\! k} 
 +   \left(  \frac{\partial {\bf f}}{\partial {\bf w}} \right)_{\!\! k}   \widehat\partial_{k}^2 {\bf w}_k    ,
  \label{chain_rule_hessian}
\end{eqnarray}
where, with the subscripts dropped to avoid unnecessary clutter (e.g., $ \widehat\partial$ should be understood as $\widehat\partial_j$ or $\widehat\partial_k$),
\begin{eqnarray}
  \widehat\partial \left(  \frac{\partial {\bf f}}{\partial {\bf w} } \right)   
=
 \left(  \frac{\partial^2 {\bf f}}{\partial {\bf w}^2 } \right)  \widehat \partial {\bf w} 
 =
  \left[  \begin{array}{ccc} 
       \widehat\partial {u_n}       &      \hat{\bf n}_{jk}^t   \widehat\partial {\rho}   & 0    \\ [1ex]
 \widehat\partial  (u_n {\bf v} )  &  \widehat\partial  (\rho {u_n}) {\bf I} +  \widehat\partial  (\rho {\bf v} ) \dyad   \hat{\bf n}_{jk}   & {\bf 0}  \\ [1ex]
 \widehat\partial  (u_n {\bf v}^2 ) / 2&  \widehat\partial (   \rho H )  \hat{\bf n}_{jk}^t +   \widehat\partial  (\rho {u_n})  {\bf v}^t  
 & \gamma  \widehat\partial {u_n}   /(\gamma-1)
              \end{array} \right],
\end{eqnarray}
 \begin{eqnarray}
 \widehat\partial  (u_n {\bf v} )   = {\bf v}     \widehat\partial  {u_n} +   u_n  \widehat\partial {\bf v}  , \quad
 \widehat\partial  (\rho {u_n})   = \rho   \widehat\partial {\bf v} +  {\bf v}   \widehat\partial  {\rho}, \quad
  \widehat\partial  ( \rho {\bf v} )   = {\bf v}     \widehat\partial  \rho +  \rho  \widehat\partial {\bf v}, \\ [2ex]
   \widehat\partial  (u_n {\bf v}^2 ) /2 =  {\bf v}^2    \widehat\partial  {u_n} /2 +  u_n {\bf v}  \widehat\partial {\bf v}  , \quad
  \widehat\partial (   \rho H ) 
   =   \frac{\gamma}{\gamma-1} \widehat\partial p + \frac{1}{2}  {\bf v}^2    \widehat\partial  \rho   + \rho  {\bf v}    \widehat\partial {\bf v}    .
 \end{eqnarray}
FSR5 with these chain-rule flux derivatives will be referred to as CFSR5. This is a very efficient scheme as the fluxes and 
their first- and second-derivatives do not need to be directly computed and stored. 
However, as mentioned earlier and will be shown later, CFSR5 fails to preserve fifth-order accuracy of FSR5 and is fourth-order accurate at best.
After all, we have two fourth-order schemes, CFSR4 and CFSR5. CFSR4 is significantly less expensive than CFSR5 since CFSR4 does not 
require second derivatives of the fluxes. But numerical experiments show that CFSR5 is more accurate than CFSR4. 

{\color{black} Note that the failure of achieving fifth-order accuracy is due to limited accuracy of the chain-rule used to 
replace the direct flux reconstruction; fourth-order errors arise exclusively from the averaged flux term, not the dissipation term. 
Seeking a way to achieve fifth-order accuracy, we explore another way of avoiding direct flux reconstruction in the next section.}

\subsection{QFSR3/QFSR4/QFSR5: Economical quadratic-form flux reconstruction}
\label{FSR:reconst_flux_quadratic_special}

As an alternative to the chain-rule, we propose the following quadratic flux reconstruction:
\begin{eqnarray}
{\bf f}_L &=&  {\bf f}_j    +     \left(  \frac{\partial {\bf f}}{\partial {\bf w}} \right)_{\!\! j}  \Delta {\bf w}_L 
+ \frac{\theta_2}{2}   \left[   \left(  \frac{\partial^2 {\bf f}}{\partial {\bf w}^2 } \right)_{\!\! j} \!    \Delta {\bf w}_L  \right]  \!  \Delta {\bf w}_L ,
 \label{QFSR3_fL_2d_euler}    \\ [1ex]
{\bf f}_R &=&  {\bf f}_k    +     \left(  \frac{\partial {\bf f}}{\partial {\bf w}} \right)_{\!\! k}  \Delta {\bf w}_R
+ \frac{\theta_2}{2}   \left[   \left(  \frac{\partial^2 {\bf f}}{\partial {\bf w}^2 } \right)_{\!\! k} \!   \Delta {\bf w}_R   \right]   \!  \Delta {\bf w}_R,
 \label{QFSR3_fR_2d_euler} 
\end{eqnarray}
where $\theta_2$ is a parameter and 
\begin{eqnarray}
\Delta {\bf w}_L = {\bf w}_L - {\bf w}_j, \quad
\Delta {\bf w}_R = {\bf w}_R - {\bf w}_k,
 \label{w_diffs}
\end{eqnarray}
with ${\bf w}_L$ and ${\bf w}_R$ given by Equations (\ref{umuscl_L2_2d_euler}) and (\ref{umuscl_R2_2d_euler}), respectively.
 Each of the fluxes ${\bf f}_L$ and ${\bf f}_R$ is expressed as a quadratic polynomial in the solution difference, but it should be noted that 
the flux reconstruction must be performed for a function whose cell-average is the point-value flux. Therefore, $\theta_2$ cannot be $1$.
As we will show later, with $\theta_2= 2/3$, this reconstruction formula leads to a third-order scheme with $\kappa_3=0$ in the solution reconstruction 
and to a fourth-order accurate scheme with $\kappa_3=\kappa-1$, which will be referred to as QFSR3 and QFSR4, respectively. 

It seems difficult to extend this approach to fifth-order accuracy without introducing a third-order term which will involve third derivatives of the flux.
However, if the flux is a quadratic function of solution variables (e.g., Burgers' equation), then the quadratic reconstruction can be modified 
to achieve fifth-order accuracy. This approach is applicable to the Euler equations whose fluxes are all quadratic in the parameter vector 
\cite{Roe_JCP_1981,idolikeCFD_VOL1_v2p6_pdf}: 
\begin{eqnarray}
{\cal F} =  {\bf z}^t {\cal P} {\bf z},  \quad {\bf z} = \left[  \sqrt{\rho} ,  \sqrt{\rho} {\bf v} ,  \sqrt{\rho} H   \right],
\end{eqnarray}
where ${\cal P}$ is a constant third-rank tensor. Then, the following reconstruction leads to a fifth-order scheme:
\begin{eqnarray}
{\bf f}_L &=&  {\bf f}_j({\bf z}_j)    +     \left(  \frac{\partial {\bf f}}{\partial {\bf z}} \right)_{\!\! j} \left[   \Delta {\bf z}_L   + L_j \right] 
+ \frac{\theta_2}{2}  \left\{ 
 \left[   \left(  \frac{\partial^2 {\bf f}}{\partial {\bf z}^2 } \right)_{\!\! j}   \!  \Delta {\bf z}_L  \right] \!  \Delta {\bf z}_L  +  Q_j
 \right\},
 \label{QFSR5_fL_2d_euler}    \\ [1ex]
{\bf f}_R &=&  {\bf f}_k({\bf z}_k)      +     \left(  \frac{\partial {\bf f}}{\partial {\bf z}} \right)_{\!\! k}  \left( \Delta {\bf z}_R   + L_k \right)
+ \frac{\theta_2}{2}  \left\{ 
  \left[   \left(  \frac{\partial^2 {\bf f}}{\partial {\bf z}^2 } \right)_{\!\! k}  \!  \Delta {\bf z}_R   \right]   \! \Delta {\bf z}_R   +  Q_k \right\} ,
 \label{QFSR5_fR_2d_euler}
\end{eqnarray}
where $\theta_2= 2/3$ as before,
\begin{eqnarray}
 L_j  =a_{Q5} \,   \widetilde{\partial_{\!\! j}^3}  {\bf z}_j , \quad 
Q_j = b_{Q5}   \left[   \left(  \frac{\partial^2 {\bf f}}{\partial {\bf w}^2 } \right)_{\!\! j}  { \widehat\partial_{\! j}^2 } {\bf z}_j   \right]  { \widehat\partial_{\! j}^2 } {\bf z}_j  
       + c_{Q5}   \left[   \left(  \frac{\partial^2 {\bf f}}{\partial {\bf w}^2 } \right)_{\!\! j}  { \widehat\partial_{\! j} } {\bf z}_j   \right]   \widetilde{\partial_{\! j}^3}  {\bf z}_j   ,  \\ [2ex]
 a_{Q5}  =  \frac{2}{15} , \quad
 b_{Q5}  =  \frac{16}{45} , \quad
 c_{Q5}  = \frac{4}{5} ,
  \label{QFSR5_fR_2d_euler_extra_terms}
\end{eqnarray}
and similarly for $L_k$ and $Q_k$. 
{\color{black} These coefficients are derived from a truncation error analysis in one dimension, and the derivation 
will be given later in Section \ref{accuracy_QFSR3_QFSR5}.}
Here, the flux and its derivatives are all expressed in terms of the parameter vector: 
\begin{eqnarray}
 {\bf f}({\bf z})  
 =
  \left[  \begin{array}{ccc} 
\displaystyle  z_n z_1    \\ [2ex]
\displaystyle  z_n \tilde{\bf v} + \frac{\gamma-1}{\gamma}  \left(  z_1  z_4 - \frac{1}{2}  \tilde{\bf v}\cdot  \tilde{\bf v} \right)  \hat{\bf n}_{jk} \\ [2ex]
\displaystyle   z_n z_4
              \end{array} \right],
\end{eqnarray}
where $ \tilde{\bf v} = (z_2, z_3)$, $z_n = \tilde{\bf v} \cdot \hat{\bf n}_{jk} $, 
\begin{eqnarray}
 \frac{\partial {\bf f}}{\partial {\bf z}} 
 =
  \left[  \begin{array}{ccc} 
z_n        &  z_1  \hat{\bf n}_{jk}^t   & 0    \\ [1ex]
 (\gamma-1) z_4 /\gamma  \hat{\bf n}_{jk}  & z_n {\bf I} +\tilde{\bf v}   \dyad   \hat{\bf n}_{jk}   -   \frac{\gamma-1}{\gamma}  \hat{\bf n}_{jk}    \dyad \tilde{\bf v}   & (\gamma-1) z_1  /\gamma  \hat{\bf n}_{jk}  \\ [1ex]
 0 &  z_4  \hat{\bf n}_{jk}^t   & z_n
              \end{array} \right].
\end{eqnarray}
\begin{eqnarray}
 \left(  \frac{\partial^2 {\bf f}}{\partial {\bf z}^2 } \right) \! \Delta  {\bf z} 
=
  \left[  \begin{array}{ccc} 
                     \Delta z_n   &  \Delta z_1    \hat{\bf n}_{jk}^t     & 0    \\ [1ex]
 (\gamma-1)   \Delta  z_4 /\gamma   \hat{\bf n}_{jk}  &    \Delta z_n {\bf I} +    \Delta \tilde{\bf v}   \dyad   \hat{\bf n}_{jk}   -   \frac{\gamma-1}{\gamma}  \hat{\bf n}_{jk}    \dyad \Delta \tilde{\bf v}    & (\gamma-1)    \Delta z_1  /\gamma \hat{\bf n}_{jk}    \\ [1ex]
 0 &     \Delta z_4  \hat{\bf n}_{jk}^t   &    \Delta z_n
              \end{array} \right].
\end{eqnarray} 
In this case, the solution difference in the dissipation term may be evaluated based on the relationship, $ d {\bf u} = \frac{1}{2}  \! \left(   \partial {\bf u}  / \partial {\bf z}  
 \right)  d {\bf z} $ (see Ref.\cite{idolikeCFD_VOL1_v2p6_pdf}), as 
\begin{eqnarray}
  {\bf u}_R -  {\bf u}_L  = \frac{1}{2}  \overline{  \frac{  \partial {\bf u} }{ \partial {\bf z} }  }  \left(    {\bf z}_R -  {\bf z}_L \right) ,
  \quad
    \frac{  \partial {\bf u} }{ \partial {\bf z} }
=
  \left[  \begin{array}{ccc} 
                  2   {z}_1  &  0     & 0    \\ [1ex]
                 { \tilde{\bf z} } &  {z}_1   {\bf I}  & {\bf 0}       \\ [1ex]
               {z}_4 / \gamma  &  (\gamma-1)  { \tilde{\bf z} }^t   / \gamma       &    {z}_1/ \gamma
              \end{array} \right],
\end{eqnarray}
where the over bar indicates the arithmetic average over the values at $j$ and $k$. Or one may evaluate the dissipation term in the original form simply by converting the parameter vector variables to the primitive and conservative variables. In this study, we employed the former. 

The resulting fifth-order scheme will be referred to as QFSR5(Z). As we will show later, QFSR5(Z) is fifth-order accurate only when the flux is quadratic in the variables used in the reconstruction; it is fourth order accurate otherwise. 
 Therefore, QFSR5(Z) is suitable for target systems, where such variables exist as in the Euler equations.

\section{FSR on Regular Grids}
\label{sec:FSR_regular_grids}

In this section, we will discuss the FSR schemes on a regular grid in one dimension.
First, we will show what these schemes reduce to on a uniform one-dimensional grid and 
show the importance of the successive application of the linear LSQ method to compute second derivatives.
Then, we will present truncation errors of the various FSR schemes for a scalar nonlinear conservation law
on uniform grids and show that QFSR5 is fifth-order accurate if a flux is quadratic in solution variables.
Analyses in one dimension are sufficient to prove accuracy in two and three dimensions because 
the FSR schemes are finite-difference schemes approximating flux derivatives.

\subsection{FSR schemes on regular grids}
\label{FSR_schemes_reg}

The FSR schemes are all in the form (\ref{diff_form_conservatrive_exact_fd_form_2d}) and thus the accuracy is determined
solely by the discretization of the flux divergence. 
Consider a scalar nonlinear conservation law:
\begin{eqnarray} 
 \frac{\partial u }{ \partial t } +  \frac{\partial f }{ \partial x }  = s(x), 
\label{scalar_cl}
\end{eqnarray}
where $f$ is a nonlinear function of the solution variable $u$ and $s(x)$ is a forcing term. 
Without loss of generality, we consider the residual $Res^x_{j}$ approximating $\partial f /\partial x$ at a node $j$ on a uniform grid of spacing $h$:
\begin{eqnarray} 
Res^x_{j}  =  \frac{ \Phi_{j+1/2} - \Phi_{j-1/2}   }{h} , 
\label{fsr_res_x}
\end{eqnarray}
where $ \Phi_{j-1/2} $ and $ \Phi_{j+1/2} $ are numerical fluxes at the left and right faces, respectively. 
For example, the flux $\Phi_{j+1/2}$ is given by
\begin{eqnarray}
 {\Phi}_{j+1/2} =  \frac{1}{2} \left[  {f}_L   +  {f}_R  \right]  - \frac{1}{2} {D}_{j+1/2} \left[  u_R -  u_L \right],
\label{cfd_numerical_flux_2d_euler_1d_scalar}
\end{eqnarray}
where the fluxes ${f}_L$ and ${f}_R$ are computed differently by various FSR schemes, ${D}_{j+1/2}$ is a dissipation coefficient, and
the left and right states $u_L $ and $u_R$ are computed, for all the FSR schemes, as
\begin{eqnarray}
u_L \!\! \!\! &=&\!\!\!\!   \kappa  \frac{u_{j} +  u_{k}}{2} +   (1-\kappa) \left[  u_j   + \widehat\partial_{\! j} u_j    \right] 
+   \kappa_3  \,  \widetilde{\partial_{\! j}^3} u 
 ,    
\label{umuscl_L2_2d_euler_scalar}\\ [1ex]
u_R \!\! \!\!  &=& \!\! \!\!  \kappa  \frac{ u_{k} + u_{j}}{2} +   (1-\kappa) \left[  u_{k} + \widehat\partial_{k} u_k   \right]
+   \kappa_3 \, \widetilde{\partial_k^3} u
, 
 \label{umuscl_R2_2d_euler_scalar}
\end{eqnarray}
\begin{eqnarray}
 \widetilde{\partial_j^3} u   = \frac{1}{2}    \left\{   \widehat\partial_{\! j} u_k - \widehat\partial_{\! j} u_j  \right \}      -  \widehat\partial_{\! j}^2 u_j   ,   \quad
\widetilde{\partial_k^3} u  =   \frac{1}{2}    \left\{   \widehat\partial_{k} u_j - \widehat\partial_{k} u_k  \right \}      -  \widehat\partial_{k}^2 u_k  , 
\end{eqnarray} 
\begin{eqnarray}
 \widetilde{\partial_j^2} u_j   = \frac{h^2}{4}  \left(  \frac{ \partial^2 u}{\partial x^2} \right)_{\!\! j}  ,   \quad
\widetilde{\partial_k^2} u_k  =  \frac{h^2}{4}  \left(  \frac{ \partial^2 u}{\partial x^2} \right)_{\!\! k}    , 
\end{eqnarray} 
\begin{eqnarray}
 \widetilde{\partial_j} u_j   = \frac{h}{2}  \left(  \frac{ \partial u}{\partial x} \right)_{\!\! j},   \quad
\widetilde{\partial_k} u_k  = -  \frac{h}{2}  \left(  \frac{ \partial u}{\partial x} \right)_{\!\! k}  , \quad
 \widetilde{\partial_j} u_k   = \frac{h}{2}  \left(  \frac{ \partial u}{\partial x} \right)_{\!\! k},   \quad
\widetilde{\partial_k} u_j  = -  \frac{h}{2}  \left(  \frac{ \partial u}{\partial x} \right)_{\!\! j} .
\end{eqnarray} 
Below, we will expand the residual and derive a truncation error for each FSR scheme applied to a scalar nonlinear conservation
law (see Ref.\cite{Nishikawa_FakeAccuracy:2020} for the importance of considering nonlinear equations). 

In each scheme, the solution gradient and second derivatives are computed by successive application of a linear LSQ method, which produces, in 
a one-dimensional uniform grid,
\begin{eqnarray} 
\left(  \frac{ \partial u}{\partial x} \right)_{\!\! j} =  \frac{ u_{j+1} - u_{j-1}   }{2 h}  ,  \quad
\left(  \frac{ \partial^2 u}{\partial x^2} \right)_{\!\! j} 
= \frac{1}{2h}   \left[ \left(  \frac{ \partial u}{\partial x} \right)_{\!\! j+1}  -\left(  \frac{ \partial u}{\partial x} \right)_{\!\! j-1}  \right] 
 =    \frac{ u_{j+2}  -2 u_j + u_{j-2}   }{ (2 h)^2}.
 \label{fsr_res_ux_uxx}
\end{eqnarray}
It is important that the second-derivative approximation is the so-called $2h$-Laplacian, spanning over a wider stencil than a typical 
central difference formula for the second derivative. The wide-stencil approximation is not a good choice for solving the Laplace/Poisson 
equation as it leads to non $h$-elliptic discretizations \cite{Multigrid_book_2001}, but here the wider stencil is essential to achieving 
high-order accuracy. 
As one can easily see, at least five neighbors are required to generate a fourth-order accurate approximation to the flux derivative;
thus it is necessary to involve third-level neighbors $j+3$ and $j-3$, which are brought in by the $2h$-Laplacian applied at the neighbors $j-1$ 
and $j+1$. All similar fifth-order schemes have some mechanisms to bring the extra-level of neighbors in their residual: the NLV5 scheme uses 
the information about elements to which an edge is incident, whose nodal gradients involve the third-level neighbors (see Ref.\cite{NLV6_INRIA_report:2008} for details) and the EBR scheme
directly constructs such a one-dimensional stencil in the direction of each edge \cite{AbalakinBakhvalovKozubskaya:IJNMF2015}. 
The successive application of a linear LSQ method proposed by Yang and Harris \cite{yang_harris:AIAAJ2016} is a convenient way of 
accomplishing the task, and for this reason, we employ it in this work. 
If the standard three-point approximation $ \frac{ u_{j+2}  -2 u_j + u_{j-2}   }{ h^2}$ is employed instead of the $2h$-Laplacian, fourth- and
higher-order accuracy will be lost. 
{\color{black} It is pointed out again that we propose the FSR schemes for unstructured-grid codes, not structured-grid 
codes, and therefore the successive application of Yang and Harris is a very important technique to form an extended stencil in a code
where only a single layer of neighbor information is available.  }

To analyze the FSR schemes, we will focus on the truncation error defined by 
\begin{eqnarray} 
 {\cal E}   \equiv \frac{1}{h} Res^x_{j},
\end{eqnarray}
where the residual is evaluated with a smooth exact solution that can be expressed by a Taylor series. 
As one can easily expect, that $ {\cal E}  $ consists of two parts: one coming from the dissipation term and
the other from the averaged flux term in the numerical flux (\ref{cfd_numerical_flux_2d_euler}). We will begin with the contribution from
the dissipation term that is common to all the FSR schemes.

\subsection{Truncation error from the dissipation term}
\label{FSR_TE_diss}

The FSR schemes differ only by the definition of the left and right fluxes and thus they have a common dissipation term. 
In this section, we derive the truncation error generated by the common dissipation term.
 Note first that the solution jumps at the left and right faces are expanded around a node $j$ as
\begin{eqnarray} 
 \!\! \!\!   ( u_R - u_L )_{j \pm 1/2}  = 
    \frac{ h^3 }{4 }   \!  \!
    \left[
     \!  (\kappa \!  -\!   \kappa_3 \!  -\!  1)\! \!   \left\{ \frac{\partial^3 u}{ \partial x^3 }  \!  \pm \!  \frac{h}{2}  \!   \frac{\partial^4 u}{ \partial x^4}    \right\}
 \!  + \!    \frac{h^2}{4} \!   (\kappa \!  - \! 2  \kappa_3 \!  - \!  1 )\frac{\partial^5 u}{ \partial x^5}  
\!    \pm  \!  \frac{h^3}{12} \!   \left(  \kappa \!  - \!  \frac{5}{2}  \kappa_3\!  -  \! 1  \right) \!  \frac{\partial^6 u}{ \partial x^6}   
         \right]   \!   \!+\!    O(h^7), 
\end{eqnarray}
 which suggests that we set
\begin{eqnarray} 
\kappa_3 = \kappa-1,
\end{eqnarray} 
so that the leading third-order term vanishes, 
\begin{eqnarray} 
  ( u_R - u_L )_{j \pm 1/2}   = 
\frac{  1 - \kappa   }{16 } \left[     \frac{\partial^5 u}{ \partial x^5}  
    \pm   \frac{h}{2}      \frac{\partial^6 u}{ \partial x^6}   
 +    O(h^2) \right]  h^5.
\end{eqnarray} 
Note that the jump will vanish for $\kappa=1$ as clearly manifested in the truncation error. 
It may cause a solver to get unstable and therefore $\kappa=1$ is typically avoided (except with high-order Runge-Kutta schemes, which can be stable with $\kappa=1$), but 
otherwise the choice is arbitrary. In this work, we set $\kappa=1/2$ and thus $\kappa_3=-1/2$, which reduces the contributions of the gradient and Hessian terms 
compared with $\kappa=0$. The reduced contribution of derivatives has been known to alleviate iterative convergence difficulties encountered for 
unstructured grids \cite{nishikawa_hyperbolic_poisson:jcp2020}. It is pleasing that a free parameter is available, with which a solver can fight against instability coming from unfavorable grid quality.
This is the reason that we use different parameters $\kappa$ and $\theta$ for the solution and flux reconstruction, respectively. However, for QFSR schemes, the choice of 
$\kappa$ is not arbitrary; it must be $1/3$ as we will discuss later. 

Substituting the jumps into the residual with the left and right fluxes ignored for a moment, we obtain the truncation error arising from the dissipation terms as 
\begin{eqnarray} 
 {\cal E}_{f_L=f_R=0} =
  \frac{  \kappa -1  }{32} 
\left[
      (   D_{j +1/2} -  D_{j -1/2}   ) h^4  \frac{\partial^5 u}{ \partial x^5 } 
 +   \frac{ D_{j +1/2} +  D_{j -1/2} } {2}  h^5 \frac{\partial^6 u}{ \partial x^6 } 
\right]  + O(h^6) ,
\label{fsr_res_x_diss}
\end{eqnarray} 
which can be written, by substituting the expansions of the dissipation coefficients, 
\begin{eqnarray} 
D_{j \pm 1/2} &=& D(u_j) \pm  \frac{1}{2} \frac{\partial D}{ \partial x} h 
+  \frac{1}{8} \left[  \frac{\partial D}{ \partial u}  \frac{\partial^2  u }{ \partial x^2}  +  \frac{\partial^2 D}{ \partial u^2} \left(   \frac{\partial  u }{ \partial x} \right)^2  \right] h^2
+ O(h^3) ,  
\end{eqnarray}
as 
\begin{eqnarray} 
 {\cal E}_{f_L=f_R=0} =
 \frac{  \kappa -1  }{32} 
\left[
  \frac{\partial D}{ \partial x}   \frac{\partial^5 u}{ \partial x^5 } 
 +  D(u_j)  \frac{\partial^6 u}{ \partial x^6 } 
\right] h^5 + O(h^7),
\label{fsr_res_x_diss_Oh5}
\end{eqnarray}
which is small enough for developing fifth-order schemes. Also, one can easily show that if we took $\kappa_3 = 0$ instead of $\kappa_3=\kappa-1$, 
the leading term would be $O(h^3)$ for any $\kappa$: 
\begin{eqnarray} 
 {\cal E}_{f_L=f_R=0,\kappa_3=0} =
-  \frac{  \kappa -1  }{8} 
\left[
  \frac{\partial D}{ \partial x}   \frac{\partial^3 u}{ \partial x^3 } 
 +  D(u_j)  \frac{\partial^4 u}{ \partial x^4 } 
\right] h^3 + O(h^5),
\label{fsr_res_x_diss_Oh3}
\end{eqnarray}
which is sufficient for third-order accurate schemes. 
Either way, these high-order dissipation terms are not sufficient to develop high-order 
schemes; second- and fourth-order error terms are generated by the averaged flux term, which must also be eliminated to achieve high-order accuracy as we will discuss next.


\subsection{FSR3, FSR4 and FSR5}
\label{accuracy_FSR3_FSR5}

For the FSR schemes, the left and right fluxes, $f_L$ and $f_R$, are directly reconstructed. 
Let us begin with $\kappa_3 = 0$, so that the dissipation term yields a third-order truncation error. 
Then, the flux reconstruction with $\theta_3 = 0 $ gives 
\begin{eqnarray} 
 {\cal E} =  
  \frac{\partial f}{ \partial x}  
  + \frac{3 \theta - 1}{12}  \frac{\partial^3 f}{ \partial x^3}   h^2 
-  \frac{  \kappa -1  }{8} 
\left[
  \frac{\partial D}{ \partial x}   \frac{\partial^3 u}{ \partial x^3 } 
 +  D(u_j)  \frac{\partial^4 u}{ \partial x^4 } 
\right] h^3 
  + \frac{ 15  \theta  - 13 }{240}  \frac{\partial^5 f}{ \partial x^5}   h^4 
  + O(h^5), 
\label{fsr_te_FSR3}
\end{eqnarray} 
showing that we will have a third-order scheme if we set
\begin{eqnarray} 
  \theta   = \frac{1}{3}. 
\end{eqnarray}
This is the FSR3 scheme. Its leading third-order term comes from the dissipation term; 
it can be eliminated by setting $\kappa_3 = \kappa-1$ as discussed in the previous section:
\begin{eqnarray} 
 {\cal E}  =  
  \frac{\partial f}{ \partial x}  
  - \frac{ 1}{30}  \frac{\partial^5 f}{ \partial x^5}   h^4 
+   \frac{  \kappa -1  }{32} 
\left[
  \frac{\partial D}{ \partial x}   \frac{\partial^5 u}{ \partial x^5 } 
 +  D(u_j)  \frac{\partial^6 u}{ \partial x^6 } 
\right] h^5
  + O(h^6) .
\label{fsr_te_FSR4}
\end{eqnarray} 
This is the FSR4 scheme. 

To achieve fifth-order accuracy, we must  choose a nonzero value for $\theta_3$. 
To determine the value, we expand the flux reconstruction and obtain
\begin{eqnarray} 
 {\cal TE} = 
  \frac{\partial f}{ \partial x}  
  + \frac{3 \theta - 1}{12}  \frac{\partial^3 f}{ \partial x^3}   h^2 
  + \frac{ 15(  \theta  - \theta_3 )  - 13 }{240}  \frac{\partial^5 f}{ \partial x^5}   h^4 
+ 
 \frac{  \kappa -1  }{32} 
\left[
  \frac{\partial D}{ \partial x}   \frac{\partial^5 u}{ \partial x^5 } 
 +  D(u_j)  \frac{\partial^6 u}{ \partial x^6 } 
\right] h^5 + O(h^6), 
\label{fsr_te}
\end{eqnarray} 
and see that the second- and fourth-order errors can be eliminated by taking 
\begin{eqnarray} 
 3 \theta - 1 = 0 ,\quad 15(  \theta  - \theta_3 )  - 13  = 0,
\end{eqnarray}
 which gives
\begin{eqnarray}
  \theta = \frac{1}{3}, \quad \theta_3 = - \frac{8}{15}.
\end{eqnarray}   
This is the FSR5 scheme. 

\subsection{CFSR3, CFSR4, and CFSR5}
\label{accuracy_CFSR}

For the CFSR schemes, we first set $\kappa_3=0$ and $\theta_3 = 0$, and obtain the following truncation error:
\begin{eqnarray} 
 {\cal E} = 
  \frac{\partial f}{ \partial x}  
  + \frac{3 \theta - 1}{12}  \frac{\partial^3 f}{ \partial x^3}   h^2 
  -  \frac{  \kappa -1  }{8} 
\left[
  \frac{\partial D}{ \partial x}   \frac{\partial^3 u}{ \partial x^3 } 
 +  D(u_j)  \frac{\partial^4 u}{ \partial x^4 } 
\right] h^3 
+ O(h^4),
\label{fsr_te_CFSR3}
\end{eqnarray} 
which shows that third-order accuracy can be achieved with $\theta=1/3$. This is the CFSR3 scheme.
Then, also as before, its leading third-order term can be eliminated with $\kappa_3 = \kappa-1$; this is
the CFSR4 scheme. 

To analyze the CFSR5 scheme, we set $\kappa_3 = \kappa-1$ but keep $\theta_3$ arbitrary to show that no values can eliminate the fourth-order truncation error.
The truncation error is given by
\begin{eqnarray} 
 {\cal E} = 
  \frac{\partial f}{ \partial x}  
  + \frac{3 \theta - 1}{12}  \frac{\partial^3 f}{ \partial x^3}   h^2 
  + C_4 h^4 
+ 
 \frac{  \kappa -1  }{32} 
\left[
  \frac{\partial D}{ \partial x}   \frac{\partial^5 u}{ \partial x^5 } 
 +  D(u_j)  \frac{\partial^6 u}{ \partial x^6 } 
\right] h^5 + O(h^6),
\label{fsr_te_CFSR5}
\end{eqnarray} 
where $C_4$ is a complicated coefficient, which, for a quadratic flux (e.g., $f=u^2/2$), simplifies to 
\begin{eqnarray} 
C_4 = 
- \frac{1}{12}  \left[    
    \left(   \frac{9}{4} +  \theta_3  \right)  \frac{\partial u}{ \partial x}  \frac{\partial^4 u}{ \partial x^4} 
+   \left(    3+ \theta_3  \right)  \  \frac{\partial^2 u}{ \partial x^2} \frac{\partial^3 u}{ \partial x^3} 
  \right]  \frac{\partial^2 f}{ \partial u^2} 
- \frac{1}{16} \left(   \theta_3 + \frac{8}{15} \right)  \frac{\partial^5 u}{ \partial x^5}   \frac{\partial f}{ \partial u}.
\end{eqnarray} 
Thus, even for a quadratic flux, the fourth-order error cannot be eliminated by any value of $\theta_3$. 
The best we can get is fourth-order accuracy. 
Then, $\theta_3$ is a free parameter; in this work we set $\theta_3 = -8/15$. 
Note that if we set $\theta_3=0$, then the scheme reduces to CFSR4. Numerical experiments will show that CFSR5 is fourth-order 
accurate but more accurate than CFSR4.

\subsection{QFSR3, QFSR4, and QFSR5}
\label{accuracy_QFSR3_QFSR5}

The QFSR schemes perform flux reconstruction based on reconstructed solutions. Hence, the parameters $\kappa$ and 
$\kappa_3$ have impact not only on odd-order error terms but also on even-order terms generated from the averaged flux term. 
Let us begin with the QFSR scheme with $\kappa_3 = 0$, which gives the following truncation error: 
\begin{eqnarray} 
 {\cal E}  &=& 
  \frac{\partial f}{ \partial x}  
  +  \frac{  1  }{4} 
\left[
\left(  \kappa - \frac{1}{3}  \right)   \frac{\partial f}{ \partial u}    \frac{\partial^3 u}{ \partial x^3} 
+ 
( \kappa + \theta_2 - 1 ) \frac{\partial^2 f}{ \partial u^2}  \frac{\partial u}{ \partial x }    \frac{\partial^2 u}{ \partial x^2 } 
+
\frac{1}{2} \left( \theta_2 - \frac{2}{3} \right) \frac{\partial^3 f}{ \partial u^3}    \left(  \frac{\partial u}{ \partial x }\right)^3 
\right] h^2 
\nonumber \\ [2ex]
&-& \frac{  \kappa -1  }{8} 
\left[
  \frac{\partial D}{ \partial x}   \frac{\partial^3 u}{ \partial x^3 } 
 +  D(u_j)  \frac{\partial^4 u}{ \partial x^4 } 
\right] h^3
+ O(h^4).
\label{fsr_te_QFSR3}
\end{eqnarray} 
The second-order error can be eliminated, fortunately, by
\begin{eqnarray}
  \kappa= \frac{1}{3}, \quad \theta_2 = \frac{2}{3},
\end{eqnarray}   
which define the QFSR3 scheme. Note that the leading third-order error here comes solely from the dissipation term; therefore 
we can achieve fourth-order accuracy by setting $\kappa_3=  \kappa-1$, which defines the QFSR4 scheme. Its truncation error 
is complicated but is given for a quadratic flux as
 \begin{eqnarray} 
 {\cal E}   =
  \frac{\partial f}{ \partial x}  
 +
 \left[
 \frac{1}{120}  \frac{\partial f}{ \partial u}    \frac{\partial^5 u}{ \partial x^5} 
 -
 \frac{5}{216}  \frac{\partial^2 f}{ \partial u^2}    \frac{\partial^2 u}{ \partial x^2}   \frac{\partial^3 u}{ \partial x^3} 
 \right] h^4 
+ O(h^5).
\label{fsr_te_QFSR4}
\end{eqnarray} 
Thus, the fourth-order error remains even for a quadratic flux.

In an attempt to achieve fifth-order accuracy, we added extra terms as in Equation (\ref{QFSR5_fR_2d_euler}) and derived 
the truncation error as 
\begin{eqnarray} 
 {\cal E}   = 
  \frac{\partial f}{ \partial x}  
 + 
\left[ 
    K_1\frac{\partial f}{ \partial u} 
 + K_2\frac{\partial^2 f}{ \partial u^2} 
 + K_3\frac{\partial^3 f}{ \partial u^3} 
 + K_4\frac{\partial^4 f}{ \partial u^4} 
 + K_5\frac{\partial^5 f}{ \partial u^5} 
\right] h^4
+ O(h^5),
\label{fsr_te_QFSR4}
\end{eqnarray} 
where we have set $\kappa=1/3$, $\kappa_3 = \kappa-1$, and $\theta_2 = 2/3$, and
\begin{eqnarray} 
 K_1 =  - \frac{1}{16}  \left( a_{Q5} - \frac{2}{15}    \right)  \frac{\partial^5 u}{ \partial x^5}, \quad
 K_2 =  - \frac{1}{16}  \left[
     \left( 2  a_{Q5} - \frac{1}{3}   c_{Q5}   \right)   \frac{\partial u}{ \partial x} \frac{\partial^4 u}{ \partial x^4}
+    \left(   a_{Q5} - \frac{2}{3}  b_{Q5} - \frac{1}{3}  c_{Q5}  + \frac{10}{27}  \right)   \frac{\partial^2 u}{ \partial x^2} \frac{\partial^3 u}{ \partial x^3}
    \right]   , \\ [2ex]
 K_3 =      
   \frac{1}{16} \left(   a_{Q5} - \frac{1}{3}   c_{Q5}+ \frac{2}{9}    \right)  \left(  \frac{\partial u}{ \partial x} \right)^2  \frac{\partial^3 u}{ \partial x^3}
+   \frac{     45 b_{Q5} + 14}{2160}\frac{\partial u}{ \partial x}    \left(  \frac{\partial^2 u}{ \partial x^2} \right)^2
  , \quad
   K_4  =  \frac{3}{80} \left( \frac{\partial^2 u}{ \partial x^2}    -1  \right)    \left(  \frac{\partial u}{ \partial x} \right)^3
   , \quad
      K_5 =  \frac{1}{720}    \left(  \frac{\partial u}{ \partial x} \right)^5.
\end{eqnarray} 
Unfortunately, fifth-order accuracy cannot be achieved in general since there are five terms that must be eliminated for the three parameters, 
$ a_{Q5}$, $b_{Q5}$, and $c_{Q5}$. 
However, if the flux is a quadratic function of the solution, we are left with the first two terms, $ K_1\frac{\partial f}{ \partial u} $ and $ K_2\frac{\partial^2 f}{ \partial u^2} $, and can uniquely determine the three 
parameters by solving
\begin{eqnarray} 
 a_{Q5} - \frac{2}{15} = 0 ,
 \quad
 2  a_{Q5} - \frac{1}{3}   c_{Q5}  = 0
 , \quad
   a_{Q5} - \frac{2}{3}  b_{Q5} - \frac{1}{3}  c_{Q5}  + \frac{10}{27} = 0,
\end{eqnarray} 
which yield the values given in Equation (\ref{QFSR5_fR_2d_euler_extra_terms}). This is the QFSR5 scheme. 
It is generally a fourth-order scheme but achieves fifth-order accuracy when the flux is quadratic in the solution variable; then QFSR5 is called QFSR5(Z).
As discussed in Section \ref{FSR:reconst_flux_quadratic_special}, it is applicable to the Euler equations.

\begin{table}[t]
\ra{1.5}
\begin{center}
{
\begin{tabular}{rcclcrlrrr}\hline\hline 
\multicolumn{1}{r}{ \multirow{2}{*}{Scheme} }                                               &
\multicolumn{1}{c}{  \multirow{2}{*}{$\kappa$} }                                              &
\multicolumn{1}{c}{  \multirow{2}{*}{$\kappa_3$} }                                          &
\multicolumn{1}{l}{  \multirow{2}{*}{Fluxes} }                                               &
\multicolumn{1}{r}{  \multirow{2}{*}{Stencil}  }             &
\multicolumn{1}{r}{  \multirow{2}{*}{LSQ}  }             &
\multicolumn{2}{c}{Accuracy}           &
\\
 \cmidrule(lr){7-8}
                                              &
                                    &
                                &
                                        &
        &
        &
\multicolumn{1}{r}{$D \ne 0 $}           &
\multicolumn{1}{r}{$D=0$}                                       
  \\ \hline 
SR2   &     Arbitrary           &   $\kappa-1$      &     ${\bf f}_L={\bf f} ( {\bf w}_L)$, ${\bf f}_R={\bf f} ( {\bf w}_R)$    & 5    &   $ \nabla {\bf w}  $     &     {   $ O(h^2) $  }  &    {\color{black}  $ O(h^2) $  }    
  \\ \hline 
FSR3     &       Arbitrary           &   $0$              &    Eqs.(\ref{FSR3_fL_2d_euler},\ref{FSR3_fR_2d_euler}):  $\theta =  \frac{1}{3}$   & 5    &   $ \nabla {\bf w}, \nabla {\cal F} $      &     {\color{blue}  $ O(h^3) $  }  &    {\color{black}  $ O(h^4) $  }   \\
FSR4     &       Arbitrary           &   $\kappa-1$  &    Eqs.(\ref{FSR3_fL_2d_euler},\ref{FSR3_fR_2d_euler}):  $\theta = \frac{1}{3}$    & 7     &  \makecell[r]{  $ \nabla {\bf w},  \nabla^2 {\bf w} $ \\ $ \nabla {\cal F} $}    &   {\color{black}  $ O(h^4) $  }   &     {\color{black}  $ O(h^4) $  }   \\
FSR5     &       Arbitrary           &   $\kappa-1$  &     Eqs.(\ref{FSR5_fL_2d_euler},\ref{FSR5_fR_2d_euler}):   \makecell[l]{  $\theta = \frac{1}{3}$ , \\ $\theta_3= -\frac{8}{15}$  } & 7    & \makecell[r]{  $ \nabla {\bf w},  \nabla^2 {\bf w} $ \\ $ \nabla {\cal F}, \nabla^2 {\cal F} $}   &  {\color{red}  $ O(h^5) $  } &     {\color{red}  $ O(h^6) $  }   \\ \hline
CFSR3   &       Arbitrary           &   $0$     &      Eqs.(\ref{FSR3_fL_2d_euler},\ref{FSR3_fR_2d_euler},\ref{chain_rule_grad}):  $\theta =  \frac{1}{3}$     & 5   &   $ \nabla {\bf w} $    &       {\color{blue}  $ O(h^3) $  }    &    {\color{black}  $ O(h^4) $  }    \\
CFSR4   &       Arbitrary           &   $\kappa-1$ &      Eqs.(\ref{FSR3_fL_2d_euler},\ref{FSR3_fR_2d_euler},\ref{chain_rule_grad}):  $\theta = \frac{1}{3}$    & 7     &   $ \nabla {\bf w},   \nabla^2 {\bf w}$    &      {\color{black}  $ O(h^4) $  }   &     {\color{black}  $ O(h^4) $  }   \\
CFSR5   &       Arbitrary           &   $\kappa-1$ &    Eqs.(\ref{FSR5_fL_2d_euler}-\ref{chain_rule_grad},\ref{chain_rule_hessian}):   \makecell[l]{  $\theta = \frac{1}{3}$ , \\ $\theta_3= -\frac{8}{15}$  } & 7      &   $ \nabla {\bf w},   \nabla^2 {\bf w}$    &    {\color{black}  $ O(h^4) $  }    &     {\color{black}  $ O(h^4) $  }   \\ \hline
QFSR3  &      $ \frac{1}{3}$        &  $0$  &     Eqs.(\ref{QFSR3_fL_2d_euler},\ref{QFSR3_fR_2d_euler}):  $\theta_2 =  \frac{2}{3}$  & 5      &   $ \nabla {\bf w}$    &     {\color{blue}  $ O(h^3) $  }  &    {\color{black}  $ O(h^4) $  }  \\
QFSR4   &     $ \frac{1}{3}$           &  $- \frac{2}{3}$  &   Eqs.(\ref{QFSR3_fL_2d_euler},\ref{QFSR3_fR_2d_euler}):  $\theta_2 =  \frac{2}{3}$    & 7   &   $ \nabla {\bf w},   \nabla^2 {\bf w}$    &    {\color{black}  $ O(h^4) $  }   &    {\color{black}  $ O(h^4) $  }   \\
QFSR5   &       $ \frac{1}{3}$           &    $- \frac{2}{3}$  &   Eqs.(\ref{QFSR5_fL_2d_euler},\ref{QFSR5_fR_2d_euler}):  $\theta_2 =  \frac{2}{3}$   & 7       &   $ \nabla {\bf w},   \nabla^2 {\bf w}$    &    {\color{black}  $ O(h^4) $  }    &      {\color{black}  $ O(h^4) $  }    \\
QFSR5(Z)   &      $ \frac{1}{3}$           &   $- \frac{2}{3}$   &   Eqs.(\ref{QFSR5_fL_2d_euler},\ref{QFSR5_fR_2d_euler}):  $\theta_2 = \frac{2}{3}$      & 7    &   $ \nabla {\bf z},   \nabla^2 {\bf z}$    &  {\color{red}  $ O(h^5) $  } &   {\color{red}  $ O(h^6) $  } 
  \\  \hline  \hline
\end{tabular}
}
\caption{A summary of FSR schemes and their orders of accuracy on grids with uniform spacing in each coordinate direction.  Stencil indicates the size of the stencil in one dimension (e.g., 5 means 2 neighbors
on each side). LSQ indicates the derivatives that need to be computed by an LSQ method: 2/3 gradient components and 3/6 independent Hessian components in 2D/3D. QFSR5(Z) is QFSR5 applied to a conservation law with fluxes quadratic in a certain set of variables ${\bf z}$.
} 
\label{Tab.FSR}
\end{center}
\end{table}

\subsection{Discussion}
\label{accuracy_discussion}

Table \ref{Tab.FSR} summarizes the FSR schemes. It is clear that the direct flux reconstruction
schemes are much more expensive than others as they require flux gradients even just for third-order accuracy.
In three dimensions, it will require computation and storage for three flux vectors and their gradients (9 additional
vectors); for fifth-order accuracy, it will require also the computation and storage for their second derivatives (18 additional
vectors). However, the FSR5 scheme is the only fifth-order scheme among those presented here that can achieve
fifth-order accuracy for arbitrary target equations. 

In terms of efficiency, the most economical third-order scheme would be the CFSR3 scheme, where only the flux gradient 
computed by the chain rule is required. As will be shown numerically, all the third-order schemes, FSR3, CFSR3, and QFSR3, 
yield almost the same level of errors. All these schemes have the same five-point stencil, which is the same as the second-order schemes.
As we will see, accuracy improvements in terms of the error level is much more dramatic when a scheme involves a wider stencil. 

For fourth-order accuracy, it is less simple to identify the most economical one because the error level varies with schemes. For purely the cost of 
computation, the most economical one would be CFSR4 as it still
requires only the normal flux gradient with the chain rule. QFSR4 requires the second derivative of the normal flux with respect to the
solution variables in addition to the flux Jacobian. CFSR5 is also fourth-order accurate and it requires the second derivative of the 
normal flux as well. However, as we will see, QFSR4 yields the most accurate solution among these fourth-order schemes on regular grids. 
Also, CFSR5 is more accurate than CFSR4, which is more accurate than FSR4. All the fourth-order schemes have the seven-point stencil and are 
generally much more accurate than the schemes of five-point stencil.

For fifth-order accuracy, FSR5 is the only choice for general conservation laws. But for a conservation law whose fluxes are quadratic
in a certain set of variables, QFSR5(Z) achieves fifth-order accuracy without direct flux reconstruction. An economical fifth-order scheme for a general conservation law, which does not require direct flux reconstruction nor third-derivatives of solution variables and fluxes, remains to be discovered. 
These fifth-order schemes have the same seven-point stencil as the fourth-order schemes. In some cases, the error level is comparable to that of the fourth-order schemes on coarse grids. 

Finally, SR2 indicates a family of second-order schemes, where the fluxes are evaluated with reconstructed solutions (which makes the scheme second-order accurate \cite{Nishikawa_FakeAccuracy:2020}), including Fromm's scheme ($\kappa=0, \kappa_3= 0$) and the scheme equivalent to UMUSCL-YH ($\kappa=1/3, \kappa_3= -2/3$). These schemes will be considered for comparison in the numerical experiments. 
 
\section{Results}
\label{results}

In this section, we present accuracy verification results for the FSR schemes. The objective here is to demonstrate the design orders of 
accuracy on regular grids and investigate, to some extent, their accuracy on non-regular grids.
For all problems, we will focus on accuracy in the interior domain; exact solutions will be specified at boundary nodes, their neighbors, and neighbors of the neighbors. High-order accuracy near boundaries is important but numerical results reported in the literature indicate that significant improvements can be achieved without high-order boundary treatment for practical problems \cite{burg_umuscl:AIAA2005-4999,yang_harris:AIAAJ2016,Barakos:IJNMF2018}. 
As mentioned earlier, further numerical tests (e.g., boundary conditions, comparisons with other schemes, shock capturing, practical problems) 
are beyond the scope of this paper and will be reported in a subsequent paper. Accuracy of integrated quantities is also a subject for future studies,
which requires high-order accuracy near boundaries and also high-order geometric approximations. 
In all cases, the effective mesh spacing is computed as the $L_1$ 
norm of the square roots  of control volumes over a grid. For comparison, we tested Fromm's scheme (SR2: $\kappa=\kappa_3=0$) and the 
UMUSCL-YH (SR2: $\kappa=1/3, \kappa_3= -2/3$). The latter will be referred to as YH for brevity. 
{\color{black} Again, in this paper, we will focus on accuracy verification. More detailed studies on practical problems, relative efficiency, dispersion/dissipation 
properties will be discussed in a subsequent paper for three-dimensional flows with a three-dimensional unstructured-grid solver. }

\subsection{One dimension}
\label{results_1D}

\subsubsection{Scalar conservation laws}
\label{results_1D_scalar_steady}

We first verify the accuracy of the FSR schemes for the steady Burgers equation: Equation (\ref{scalar_cl}) with $f=u^2/2$ and
$s(x) = A \sin( A x ) \cos( A x )$, so that the exact steady solution is given by $ u(x) =\sin(A x)$, where $A=1.23$. The steady residuals are solved
by an implicit solver \cite{nishikawa_liu_jcp2018} until the residual is reduced by seven orders of magnitude from an initial norm 
for a series of uniform grids with 32, 64, 128, 256 nodes. The error is computed in the $L_\infty$ norm of the solution over all nodes.

  \begin{figure}[th!]
    \centering
      \hfill  
                \begin{subfigure}[t]{0.48\textwidth}
        \includegraphics[width=\textwidth]{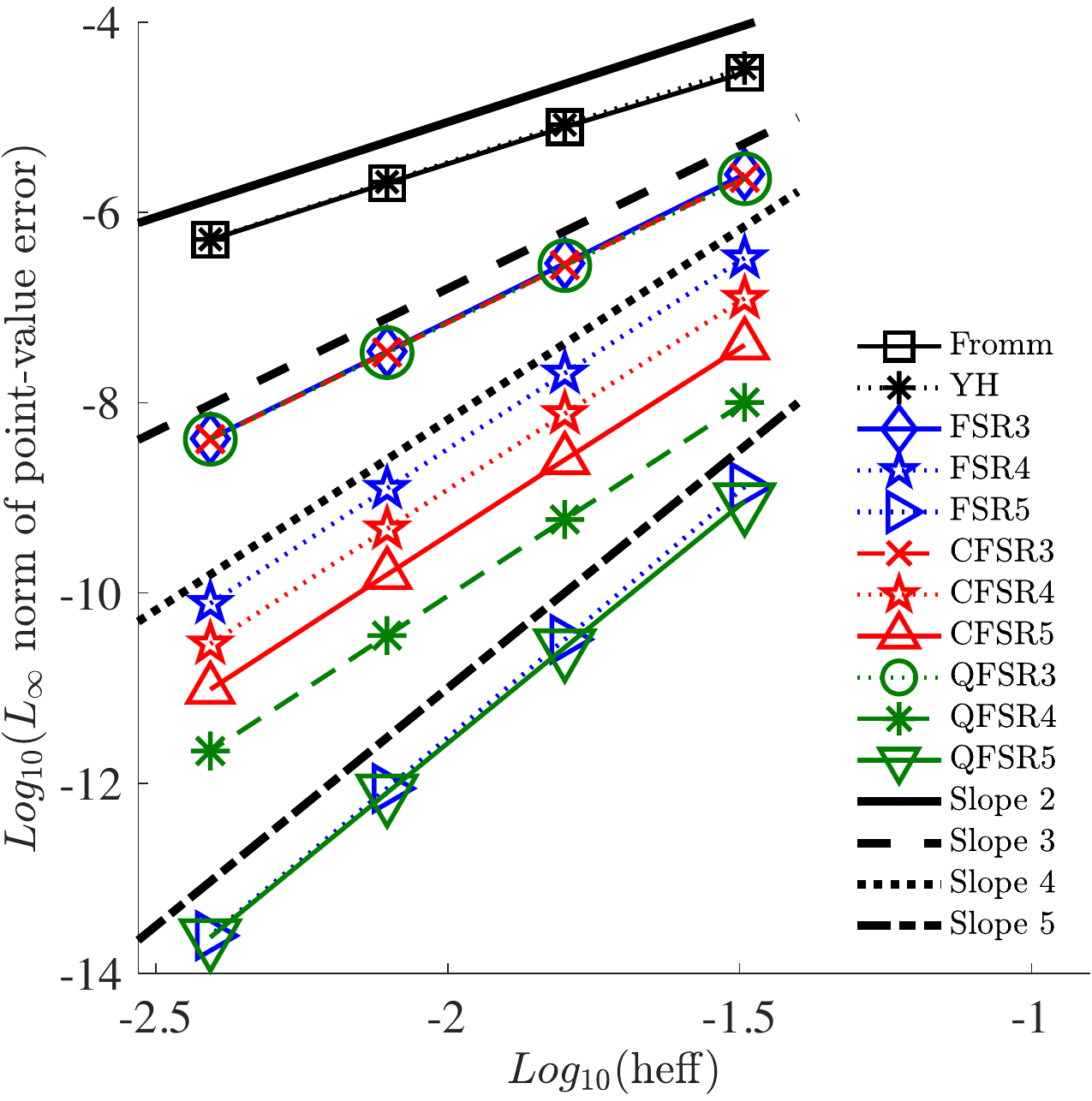}
          \caption{$f = u^2/2$.}
       \label{fig:oned_scalar_steady_de_n2}
      \end{subfigure}
      \hfill
          \begin{subfigure}[t]{0.48\textwidth}
        \includegraphics[width=\textwidth]{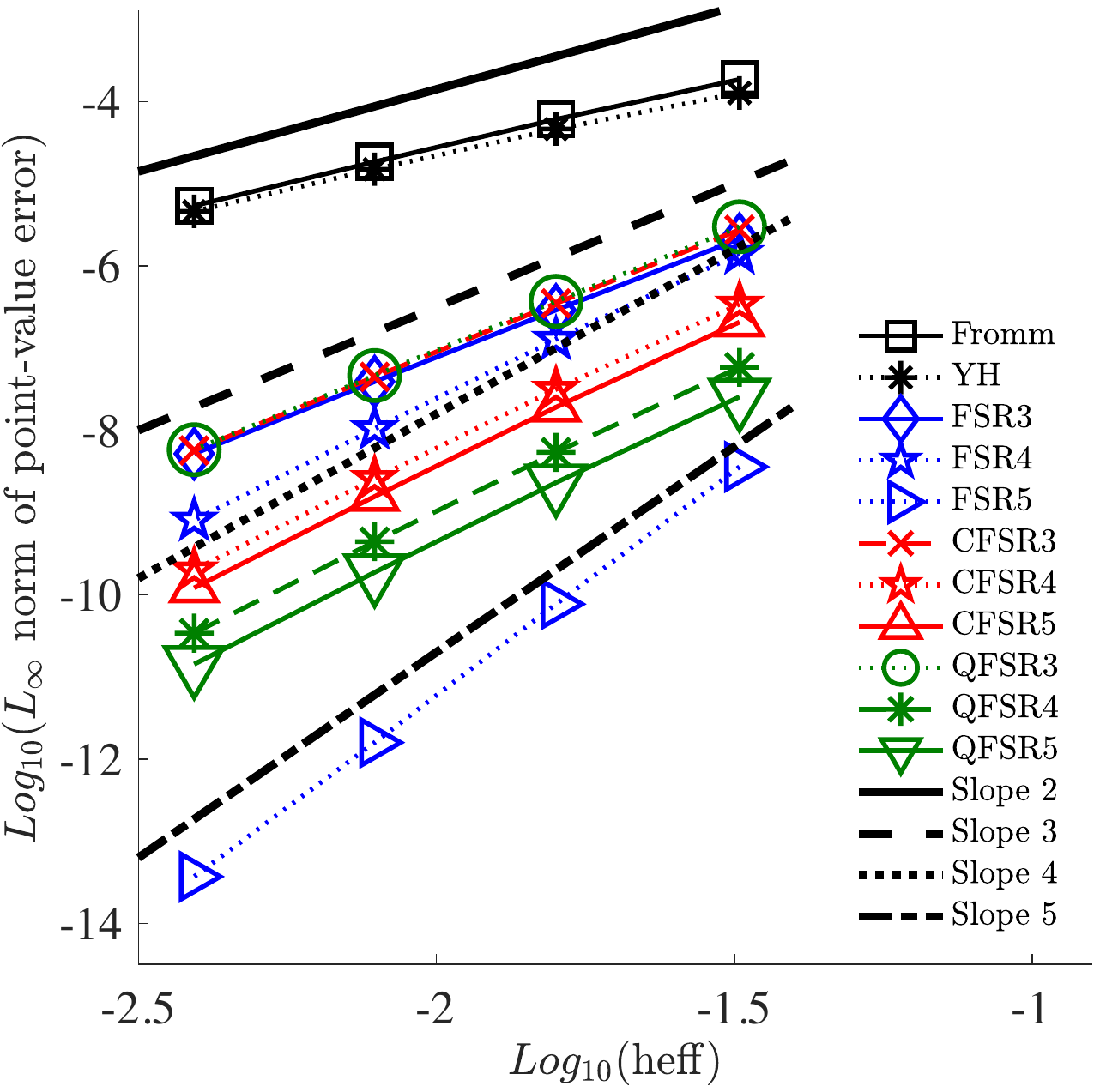}
          \caption{$f = u^3/3$.}
       \label{fig:oned_scalar_steady_de_n3}
      \end{subfigure}
      \hfill    
      \\
            \caption{
\label{fig:oned_scalar_steady}%
Error convergence results for scalar conservation laws in one dimension.
} 
\end{figure}

The results are shown in Figure \ref{fig:oned_scalar_steady_de_n2}. As can be seen, the Fromm's scheme and YH are second-order accurate.
The FSR3, FSR4, and FSR5 schemes are third-, fourth- and fifth-order accurate, respectively, as expected. Also, the CFSR3 scheme is third-order accurate
while the CFSR4 and CFSR5 schemes are fourth-order accurate. CFSR5 has a larger stencil and gives indeed lower errors than CFSR4. Finally, the QFSR schemes are also verified as third-, fourth-, and fifth-order accurate. Here, QFSR5 is fifth-order accurate because the flux is quadratic in the solution variable. 

To demonstrate the impact of non-quadratic fluxes on the QFSR5 scheme, we solve another steady problem for a scalar conservation 
law (\ref{scalar_cl}) with $f=u^3/3$ with the source term defined such that  $u(x) = \sin(1.23 x)$ is the exact solution. Results are shown in
Figure \ref{fig:oned_scalar_steady_de_n3}. As can be clearly seen, QFSR5 reduces to fourth-order accurate as predicted by the analysis although
it is slightly more accurate than QFSR4. Note also that CFSR5 gives, although very slight, lower errors than CFSR4. 

For all the results, the third-order schemes, FSR3, CFSR3, and QFSR3 produce errors at almost the same level. It indicates that the 
most economical third-order scheme might be CFSR3. For fourth-order accurate schemes, the level of errors varies and FSR4 seems to be the least accurate. 
Although CFSR4 and CSFR5 are slightly more efficient in terms of the cost per iteration than QFSR4 and QFSR5, but the QFSR schemes 
may turn out to be more efficient especially in the case of quadratic fluxes. 
{\color{black} Detailed studies for relative efficiency are beyond the scope of
this paper and will be investigated for practical problems with a three-dimensional turbulent-flow solver.}

\subsubsection{Steady Euler equations}
\label{results_1D_euler_steady}

We next consider a steady problem for the Euler equations in one dimension with the exact solution set for the primitive variables as
\begin{eqnarray}
{\bf w}^{exact} (x)
=
\left[
\begin{array}{c}
\rho^{exact} (x) \\
u^{exact} (x)  \\
p^{exact} (x)
\end{array}
\right]
=
\left[
\begin{array}{c}
    1.0 +   0.29 \rho \sin( 2.3 \pi ) +  x^5 \\
 0.2 +   0.30 \rho \sin( 2.0 \pi ) +  x^5 \\
 1.7 +   0.27 \rho \sin( 2.5 \pi ) +  x^5 
\end{array}
\right], 
\label{euler_1d_exact_sol_steady}
\end{eqnarray}
and the forcing term defined by ${\bf s}(x) = \partial_x {\bf f}({\bf w}^{exact}(x))$. Accuracy verification is performed with uniform 
grids of 22, 44, 88, 176 nodes. 
The system of steady residual equations is solved by a pseudo-time integration with the three-stage SSP Runge-Kutta scheme \cite{SSP:SIAMReview2001} 
with a local time step at CFL$=0.99$. The solver is taken to be converged when the residual is reduced by seven orders of 
magnitude in the $L_1$ norm, starting from the initial norm computed with the initial constant solution: ${\bf w}_i = (1, 0.2 , 1.7 )$. Note that the exact solution has been defined to avoid fake high-order accuracy of the SR2 schemes as discussed in Refs.\cite{Nishikawa_3rdMUSCL:2020IJNMF,Nishikawa_3rdQUICK:2020,Nishikawa_FakeAccuracy:2020}.

Results are shown in Figure \ref{fig:oned_euler_steady_de}; these are very similar to those obtained for the Burgers equation. 
Here, we tested both QFSR5 and QFSR5(Z), where the former performs the solution reconstruction with the primitive variables and the latter with the
parameter vector variables. Clearly, only QFSR5(Z) achieves fifth-order accuracy as expected; QFSR5 is fourth-order accurate because the flux is not quadratic in the primitive variables.

  \begin{figure}[th!]
    \centering
      \hfill  
                \begin{subfigure}[t]{0.48\textwidth}
        \includegraphics[width=\textwidth]{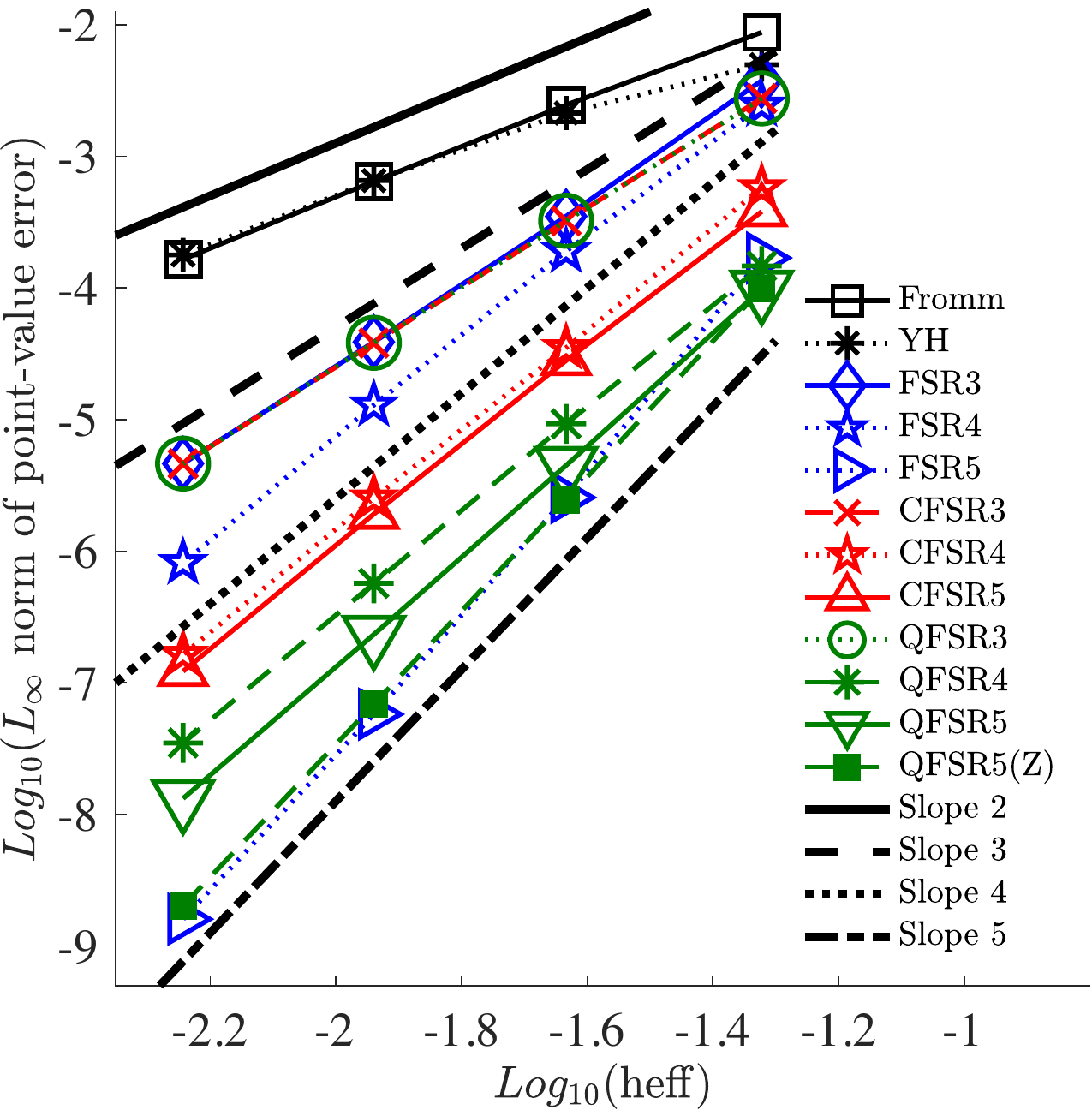}
          \caption{Steady problem.}
       \label{fig:oned_euler_steady_de}
      \end{subfigure}
      \hfill
          \begin{subfigure}[t]{0.48\textwidth}
        \includegraphics[width=\textwidth]{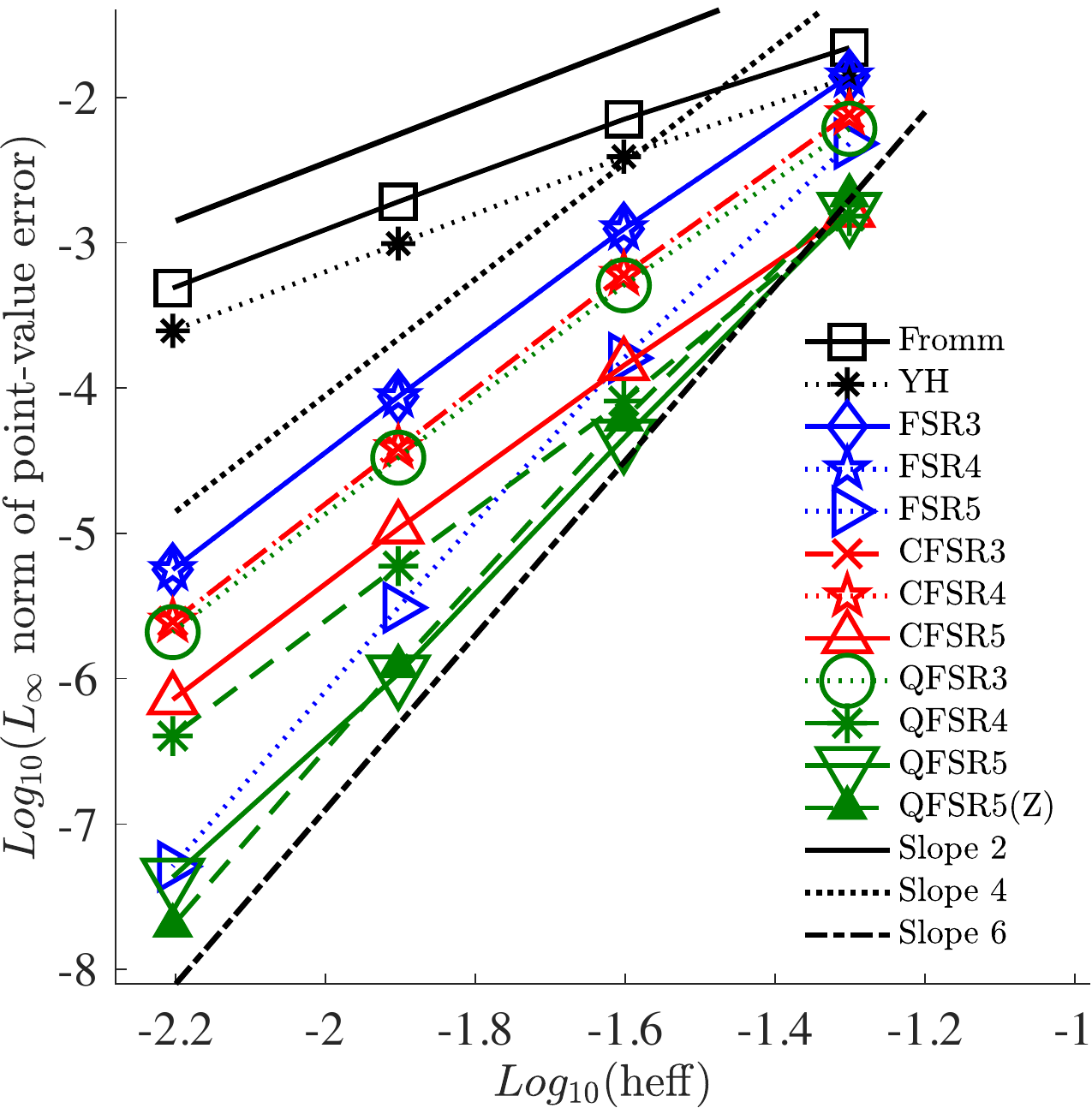}
          \caption{Unsteady problem (zero dissipation schemes).}
       \label{fig:oned_euler_unsteady_de}
      \end{subfigure}
      \hfill    
      \\
            \caption{
\label{fig:oned_euler_steady}%
Error convergence results for the Euler equations in one dimension.
} 
\end{figure}

\subsubsection{Unsteady Euler equations}
\label{results_1D_steady}

We verify the FSR schemes for an unsteady problem with a simple acoustic-wave solution taken from Ref.\cite{idolikeCFD_VOL1_v2p6_pdf}: 
\begin{eqnarray}
\left[
\begin{array}{c}
\rho^{exact} (x,t) \\
u^{exact} (x,t)  \\
p^{exact} (x,t)
\end{array}
\right]
=
\left[
\begin{array}{c}
\displaystyle  \left\{    1 + \frac{\gamma-1}{2} V(x,t)     \right\}^{  \frac{2}{\gamma-1} } \\ [2.5ex]
\displaystyle M_\infty +   V(x,t)  \\ [2.5ex]
\displaystyle  \left\{    1 + \frac{\gamma-1}{2} V(x,t)     \right\}^{  \frac{2 \gamma}{\gamma-1} } \frac{1}{\gamma} 
\end{array}
\right], 
\quad
 V(x,t)  =   \frac{\sin \left[  2 \pi \left\{ x  - ( u + a ) t \right\} \right]}{ \pi t_s (\gamma+1) },
\label{euler_1d_exact_sol_unsteady}
\end{eqnarray}
where $M_\infty = 1.7$, $t_s = 0.2$, $ u = u^{exact} $, and $a= \sqrt{ \gamma p^{exact}  / \rho^{exact}  }$. The exact solution is defined implicitly; it can be easily solved numerically by the fixed-point iteration, where the above exact solution formula is applied repeatedly from the initial values computed with $t=0$. This exact solution is valid until a shock is formed approximately at $t = t_s$. Accuracy verification is performed for the solution at $t=0.07$ at a fixed time step $\Delta t = 2.0$E-$06$ (the total $3,5000$ time steps) over a series of grids with 21, 41, 81, 161 nodes. For this problem, we tested schemes with zero dissipation to verify the orders of accuracy indicated in the right-most column of Table \ref{Tab.FSR}. For time integration, we employ the three-stage SSP Runge-Kutta scheme \cite{SSP:SIAMReview2001} with CFL$=0.99$.

Results are shown in Figure \ref{fig:oned_euler_unsteady_de}. As expected, the second-order schemes remain second-order accuarte. The third-order schemes, FSR3, CFSR3, and QFSR3, are all fourth-order accurate; it verifies that the leading third-order term comes from the dissipation. On the other hand, the fourth-order schemes remain fourth-order accurate. This is expected because the leading fourth-order error comes from the flux average term, not from the dissipation. 
The fifth-order schemes, FSR6 and QFSR5, now achieve sixth-order accuracy as expected. It is interesting to note that QFSR5 is nearly as accurate as FSR5 and QFSR5(Z).

\subsection{Two dimensions}
\label{results_2D}

Finally, we consider the Euler equations in two dimensions. Here, we do not consider FSR3, FSR4, and FSR5 because these schemes are significantly more expensive than other economical FSR schemes and are not worth implementing in a code.

\subsubsection{Steady problem}
\label{results_2D_steady}

We consider the Euler equations in two dimensions with the following exact solution: 
\begin{eqnarray}
\rho                &=&  1.00  + 0.2  \sin[   \pi (  2.3  x + 2.3 y ) ], \\
u                    &=&  0.15 + 0.2  \sin[    \pi ( 2  x + 2 y )  ], \\ 
v                    &=&  0.02  + 0.2 \sin[    \pi ( 2   x + 2 y ) ], \\ 
p                   &=&  1.00  + 0.2 \sin[     \pi ( 2.5  x + 2.5 y ) ].
\end{eqnarray}
 The forcing terms are computed numerically and the residual equations are
solved by an implicit solver (until the residual is reduced by eight orders of magnitude from an initial norm) as described in Ref.\cite{nishikawa_centroid:JCP2020}. 
For this problem, we consider four different types of grids: regular quadrilaterals, right triangles, equilateral triangles, and irregular triangles. 

Figure \ref{fig:twod_euler_steady_quad_reg} shows results for regular quadrilateral grids: $n$$\times$$n$ nodes, where $n=32, 48, 64, 80, 96, 112, 128$. These results are very similar to those of one dimensional problems.  Observe that YH is second-order accurate but more accurate than Fromm's scheme. 
To demonstrate that the accuracy is not affected as long as the mesh spacing is constant in each coordinate direction, we performed the same computations with the $y$-coordinate rescaled as $y \leftarrow 0.1 y$ and the coefficients to $y$ in the sine function of the exact solutions multiplied by 10 (to ensure significant solution variation in the $y$ direction). Hence, the mesh spacing is 10 times larger in $x$-direction than in the $y$-direction. Results are shown in Figure \ref{fig:two_euler_steady_quad_reg_error_ar10}. As can be seen clearly, the results are very similar to those for the isotropic grids. 
{\color{black} Note that the error levels of YH and CFSR5 get higher relative to others compared with the isotropic-grid case. However, there is no apparent reason that these schemes should be more sensitive to the cell aspect ratio. A further investigation is necessary to study the effects of high-aspect-ration grids on these economical schemes, especially for highly thin and curved grids typical in practical turbulent-flow simulations. The same is true for other types of grids considered in the rest of the paper; hence results for high-aspect-ratio grids will be omitted for other grid types. }

  \begin{figure}[th!]
    \centering
      \hfill  
                \begin{subfigure}[t]{0.32\textwidth}
        \includegraphics[width=\textwidth]{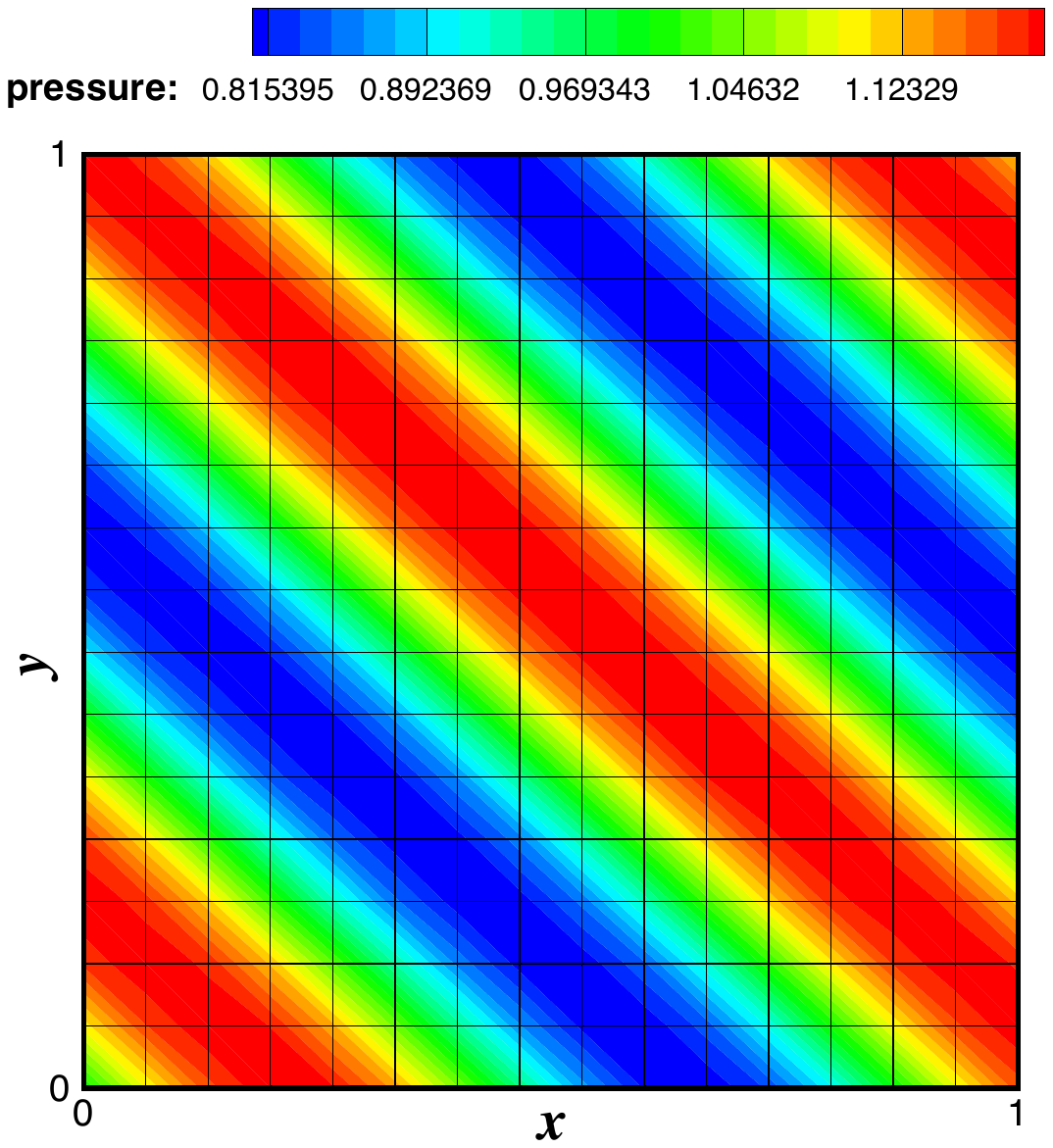}
          \caption{Grid.}
       \label{fig:two_euler_steady_quad_reg_grid}
      \end{subfigure}
      \hfill
          \begin{subfigure}[t]{0.32\textwidth}
        \includegraphics[width=\textwidth]{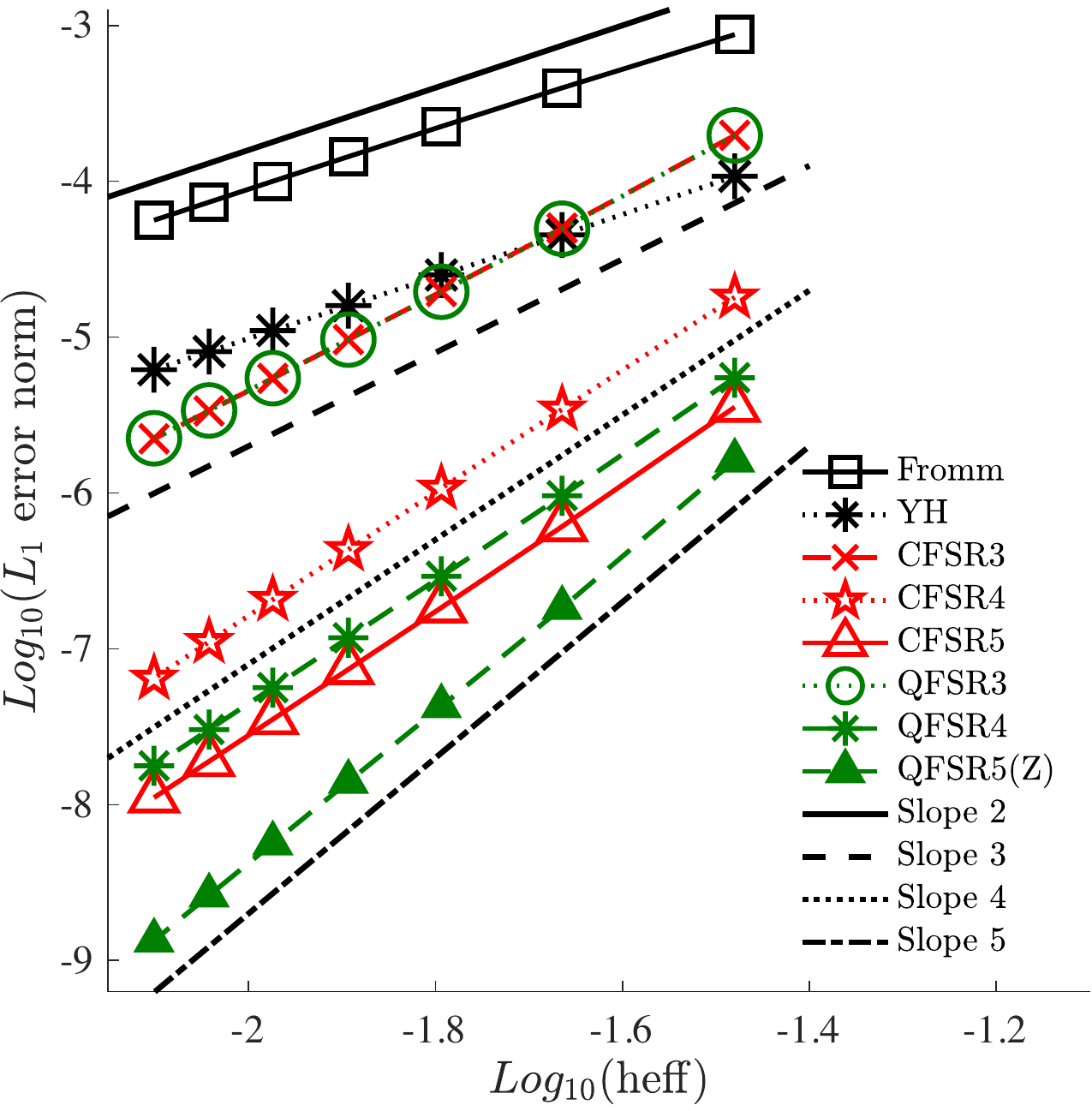}
          \caption{Error convergence: Isotropic grids}
       \label{fig:two_euler_steady_quad_reg_error}
      \end{subfigure}
      \hfill
          \begin{subfigure}[t]{0.32\textwidth}
        \includegraphics[width=\textwidth]{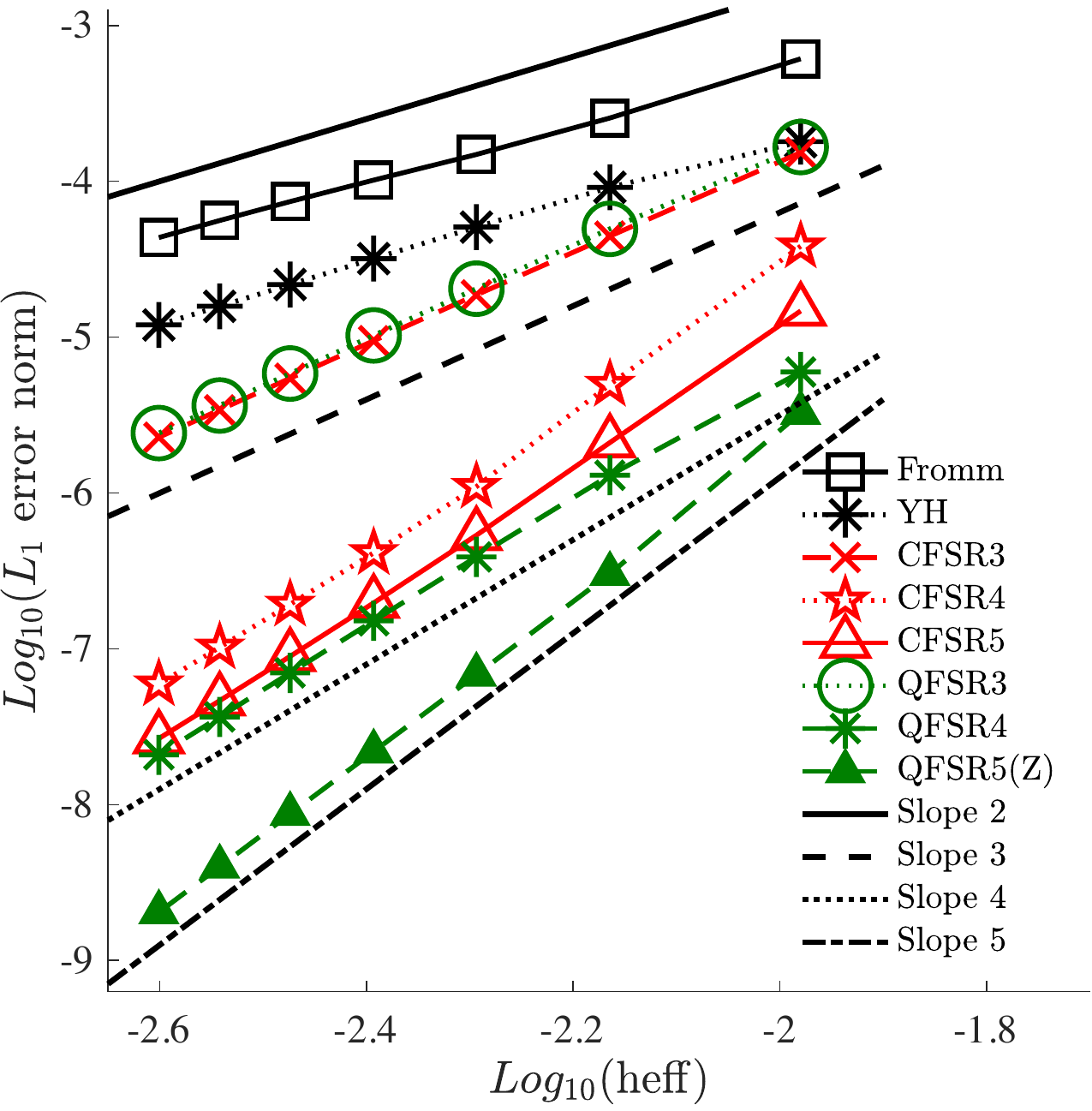}
          \caption{Error convergence: High aspect ratio grids $h_x/h_y = 10$.}
       \label{fig:two_euler_steady_quad_reg_error_ar10}
      \end{subfigure}
            \caption{
\label{fig:twod_euler_steady_quad_reg}%
Regular quadrilateral grids: error convergence results for the steady inviscid problem.
} 
\end{figure}
  \begin{figure}[th!]
    \centering
      \hfill  
                \begin{subfigure}[t]{0.48\textwidth}
        \includegraphics[width=\textwidth]{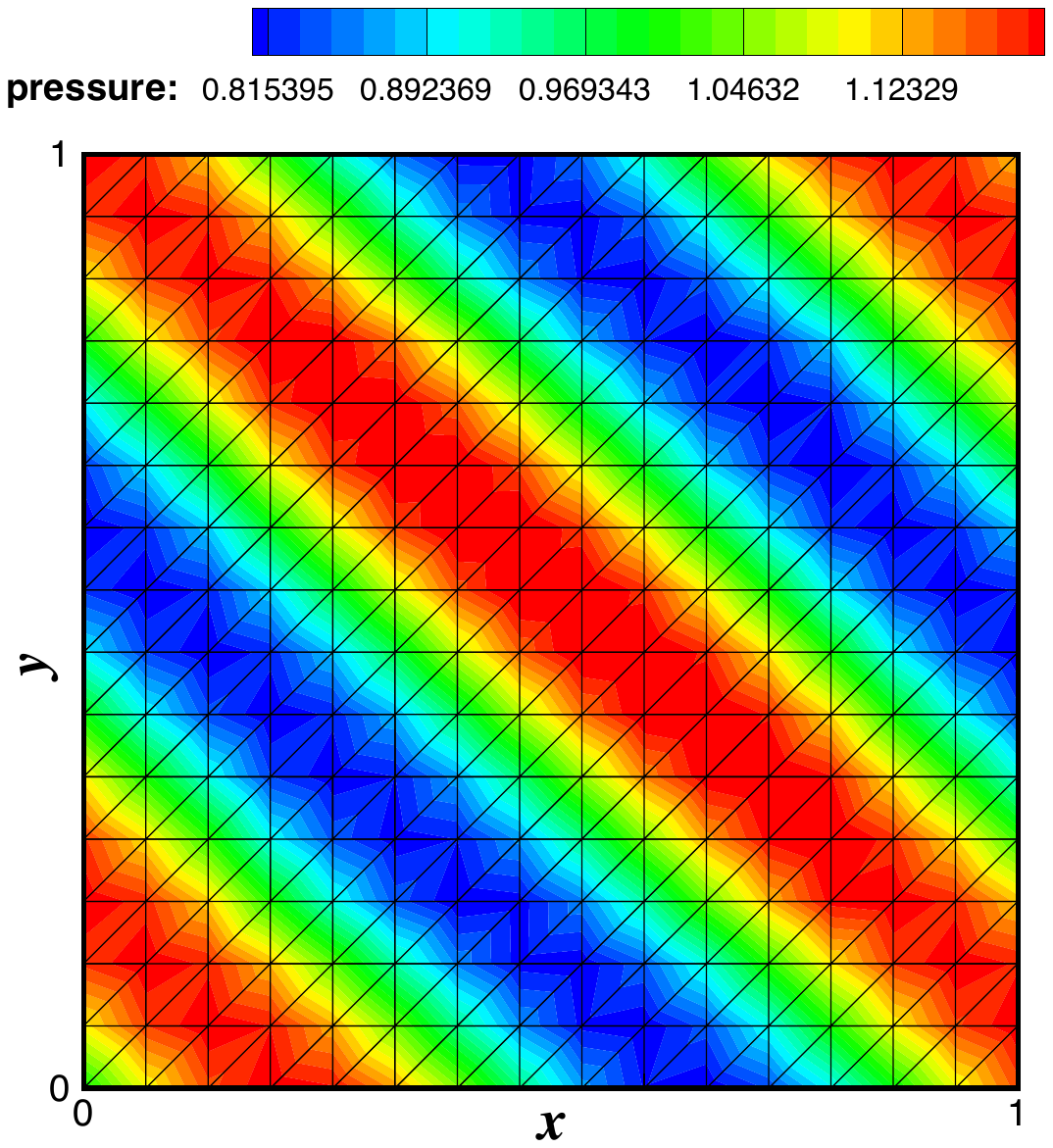}
          \caption{Grid.}
       \label{fig:two_euler_steady_tria_reg_grid}
      \end{subfigure}
      \hfill
          \begin{subfigure}[t]{0.48\textwidth}
        \includegraphics[width=\textwidth]{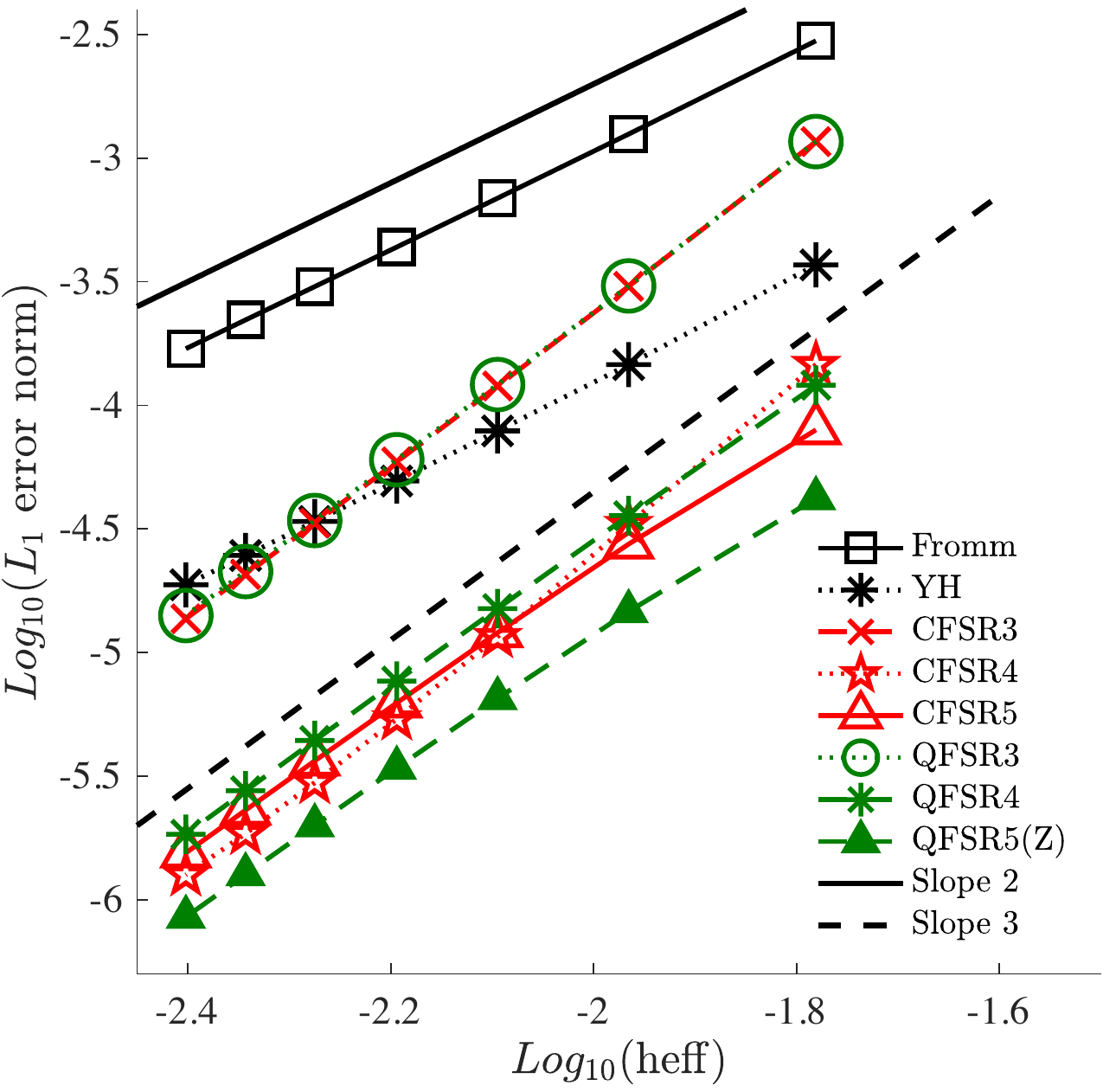}
          \caption{Error convergence.}
       \label{fig:two_euler_steady_tria_reg_error}
      \end{subfigure}
            \caption{
\label{fig:twod_euler_steady_tria_reg}%
Right triangular grids: error convergence results for the steady inviscid problem.
} 
\end{figure}
  \begin{figure}[th!]
    \centering
      \hfill  
                \begin{subfigure}[t]{0.48\textwidth}
        \includegraphics[width=\textwidth]{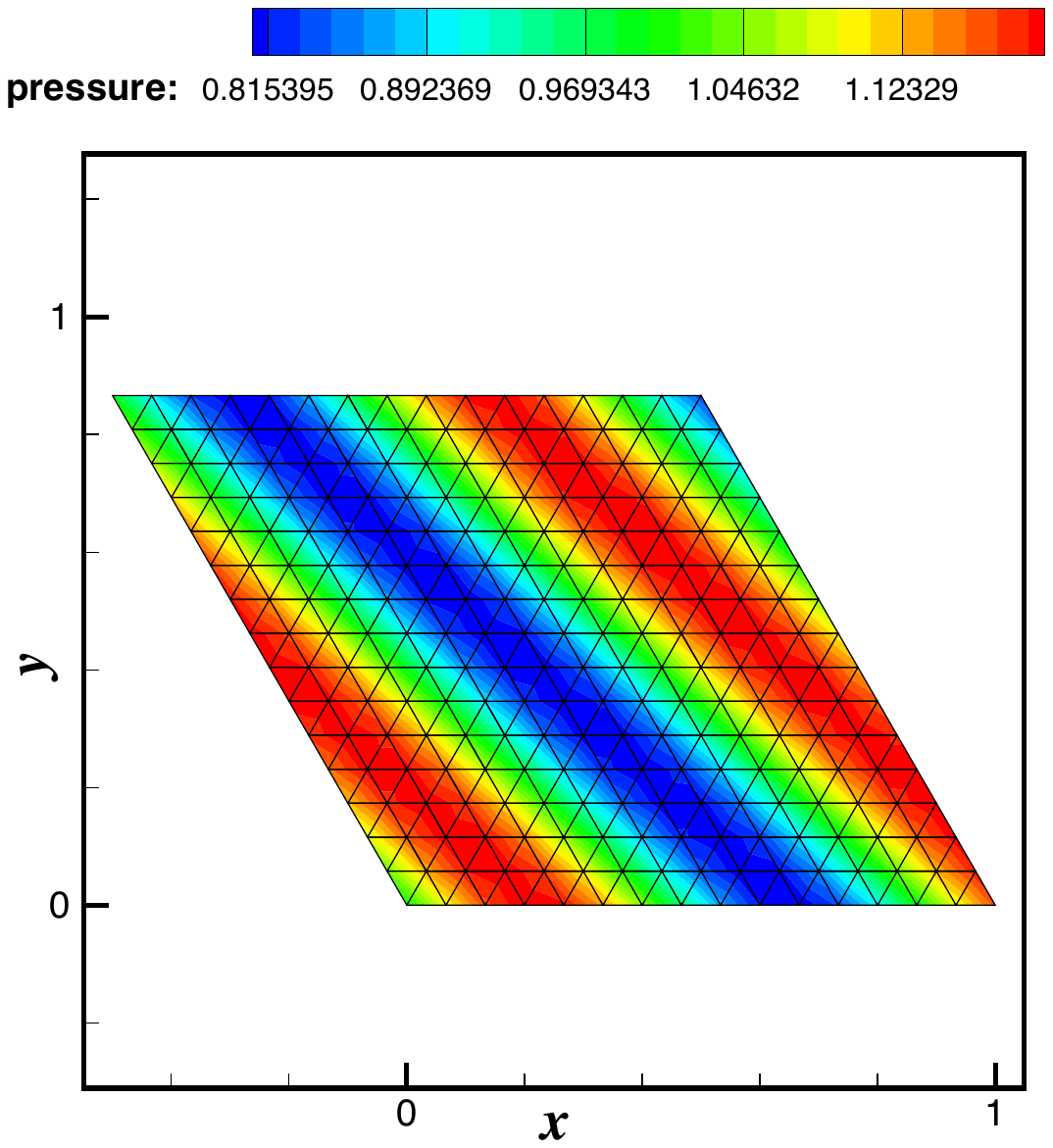}
          \caption{Grid.}
       \label{fig:two_euler_steady_tria_equi_grid}
      \end{subfigure}
      \hfill
          \begin{subfigure}[t]{0.48\textwidth}
        \includegraphics[width=\textwidth]{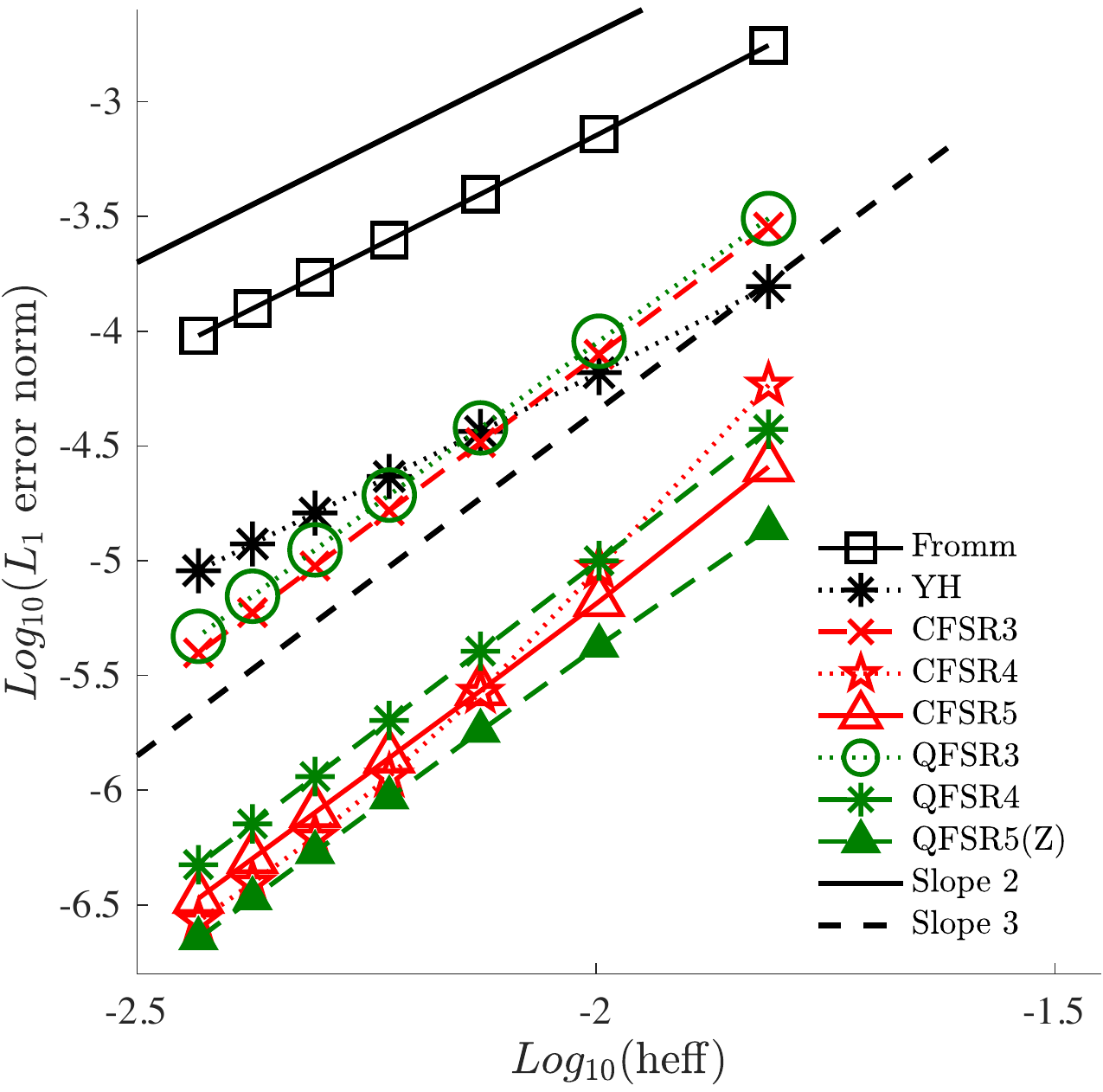}
          \caption{Error convergence.}
       \label{fig:two_euler_steady_tria_equi_error}
      \end{subfigure}
            \caption{
\label{fig:twod_euler_steady_tria_equi}%
Equilateral triangular grids: error convergence results for the steady inviscid problem.
} 
\end{figure}
  \begin{figure}[th!]
    \centering
      \hfill  
                \begin{subfigure}[t]{0.48\textwidth}
        \includegraphics[width=\textwidth]{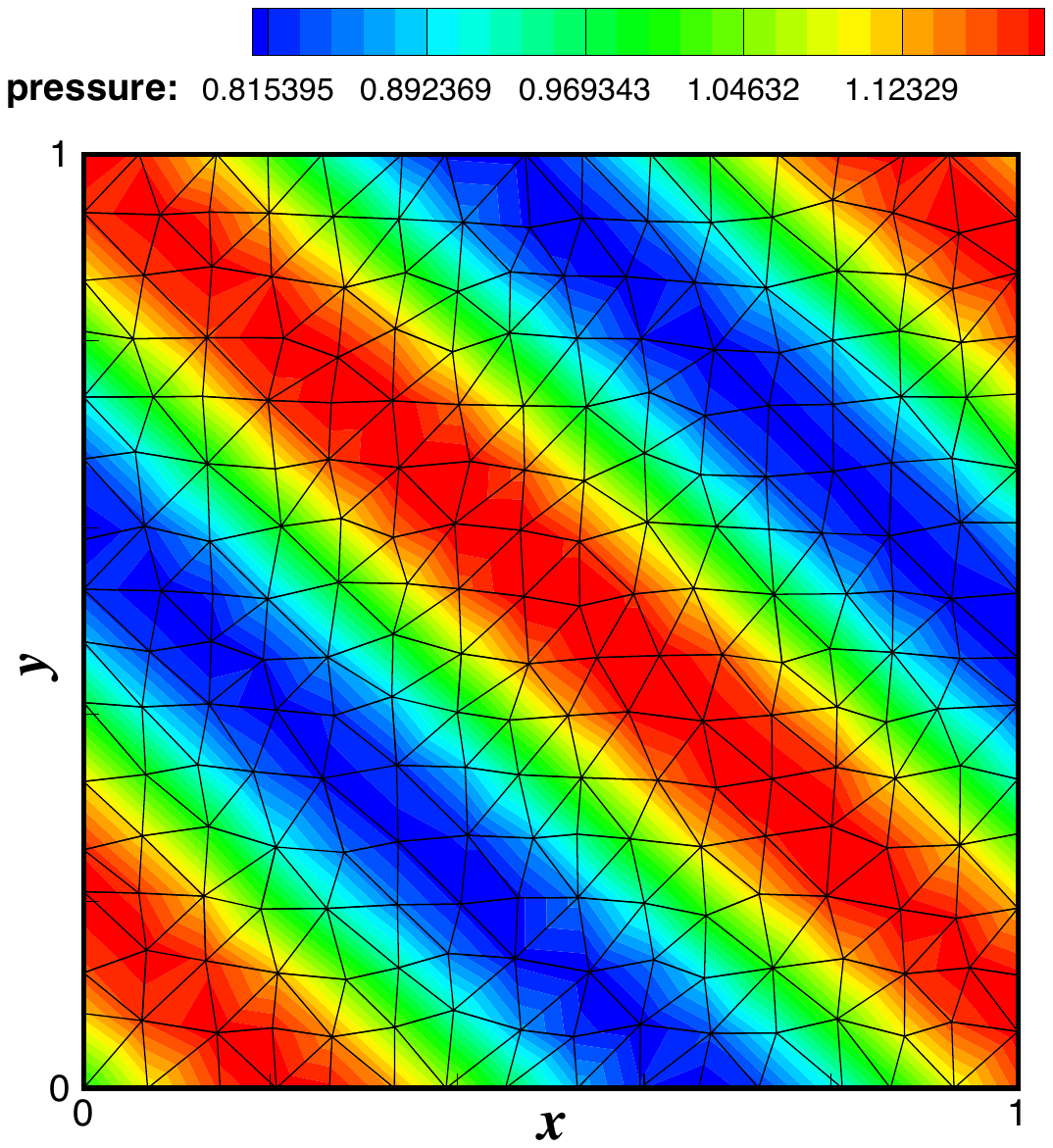}
          \caption{Grid.}
       \label{fig:two_euler_steady_tria_irrg_grid}
      \end{subfigure}
      \hfill
          \begin{subfigure}[t]{0.48\textwidth}
        \includegraphics[width=\textwidth]{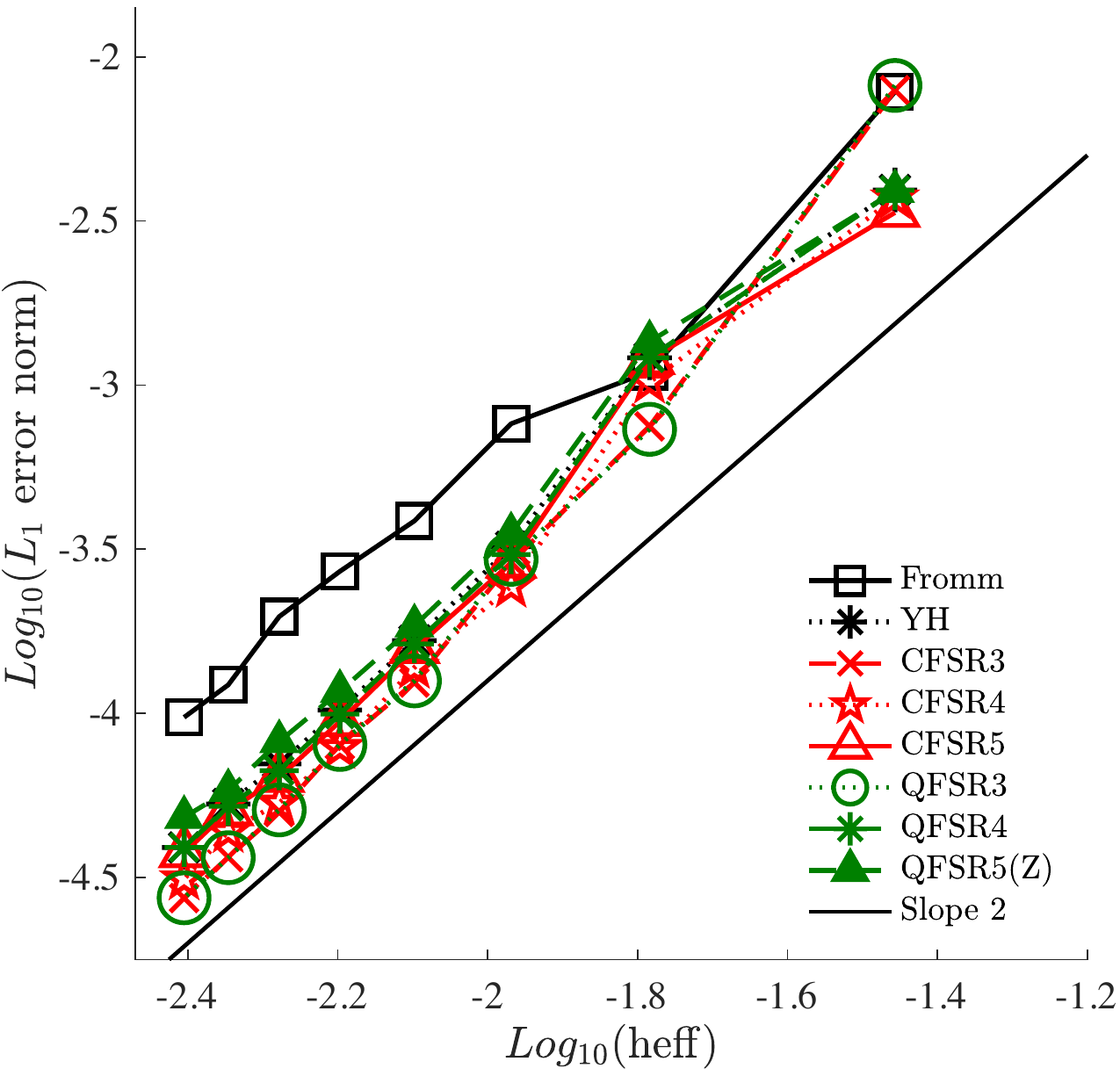}
          \caption{Error convergence.}
       \label{fig:two_euler_steady_tria_irrg_error}
      \end{subfigure}
            \caption{
\label{fig:twod_euler_steady_tria_irrg}%
Iregular triangular grids: error convergence results for the steady inviscid problem.
} 
\end{figure}

Figure \ref{fig:twod_euler_steady_tria_reg} shows results for right triangular grids, which are generated from $n$$\times$$n$ quadrilateral grids, where $n=32, 48, 64, 80, 96, 112, 128$. Again, YH is second-order accurate but more accurate than Fromm's scheme. In this case, YH is more accurate than third-order schemes on coarse grids. For the FSR schemes, we observe that third-order schemes remain third-order accurate but fourth- and fifth-order schemes have lost their design accuracy. These high-order schemes are only third-order accurate on right triangular grids. Apparently, the accuracy deterioration is due to the fact that these schemes do not reduce to their one-dimensional versions along each grid line. Such scheme can be constructed by performing one-dimensional reconstruction along each grid line \cite{AbalakinBakhvalovKozubskaya:IJNMF2015}; but it will require to store a collection of nodes at each edge along its direction. 
Note, however, that these schemes are still significantly more accurate than the second- and third-order schemes.

Figure \ref{fig:twod_euler_steady_tria_equi}} shows results for equilateral triangular grids. 
The grids were generated from the right triangular grids by shifting the nodes to the left as shown in 
Figure \ref{fig:two_euler_steady_tria_equi_grid}. As can be seen, results are very similar to the previous case.
In this case, it was hoped that fourth- and fifth-order schemes would keep the design orders of accuracy, but 
they reduce to third-order accurate because again they do not reduce to purely one-dimensional schemes along each grid line.
It is expected to achieve design orders of accuracy if the LSQ gradient reduces to the centra difference formula at a node when projected
along each grid line. Such a gradient formula has not been discovered yet; it is nor clear if it is even possible.

Figure \ref{fig:twod_euler_steady_tria_irrg} shows results for irregular triangular grids, generated from the right triangular grids with random diagonal
swapping and nodal perturbation.
As can be seen in Figure \ref{fig:two_euler_steady_tria_irrg_grid}, the number of neighbors of each node is random and the grid is irregularly spaced.
As expected, all the schemes reduce to second-order methods as shown in Figure \ref{fig:two_euler_steady_tria_irrg_error}; these economically high-order schemes are high-order accurate only when a grid has regularity, for example, as shown in the previous cases. 
{\color{black} The results indicate that accuracy will be deteriorated on an irregular-grid region of a hybrid composed of a smooth grid and an irregular
grid. Note that the economical high-order schemes are not 
designed for use in fully irregular unstructured grids, these results are presented here merely to demonstrate the deterioration to
 second-order accuracy.} 
For irregular triangular (and tetrahedral) grids, only EB3 can achieve higher than 
second-order accuracy (i.e., third-order accurate) among the economically high-order schemes. See Ref.\cite{nishikawa_liu_source_quadrature:jcp2017}.

\subsubsection{Unsteady problem}
\label{results_2D_steady}

\begin{figure}[htbp!]
    \centering
      \begin{subfigure}[t]{0.48\textwidth}
        \includegraphics[width=\textwidth,trim=0 0 0 0,clip]{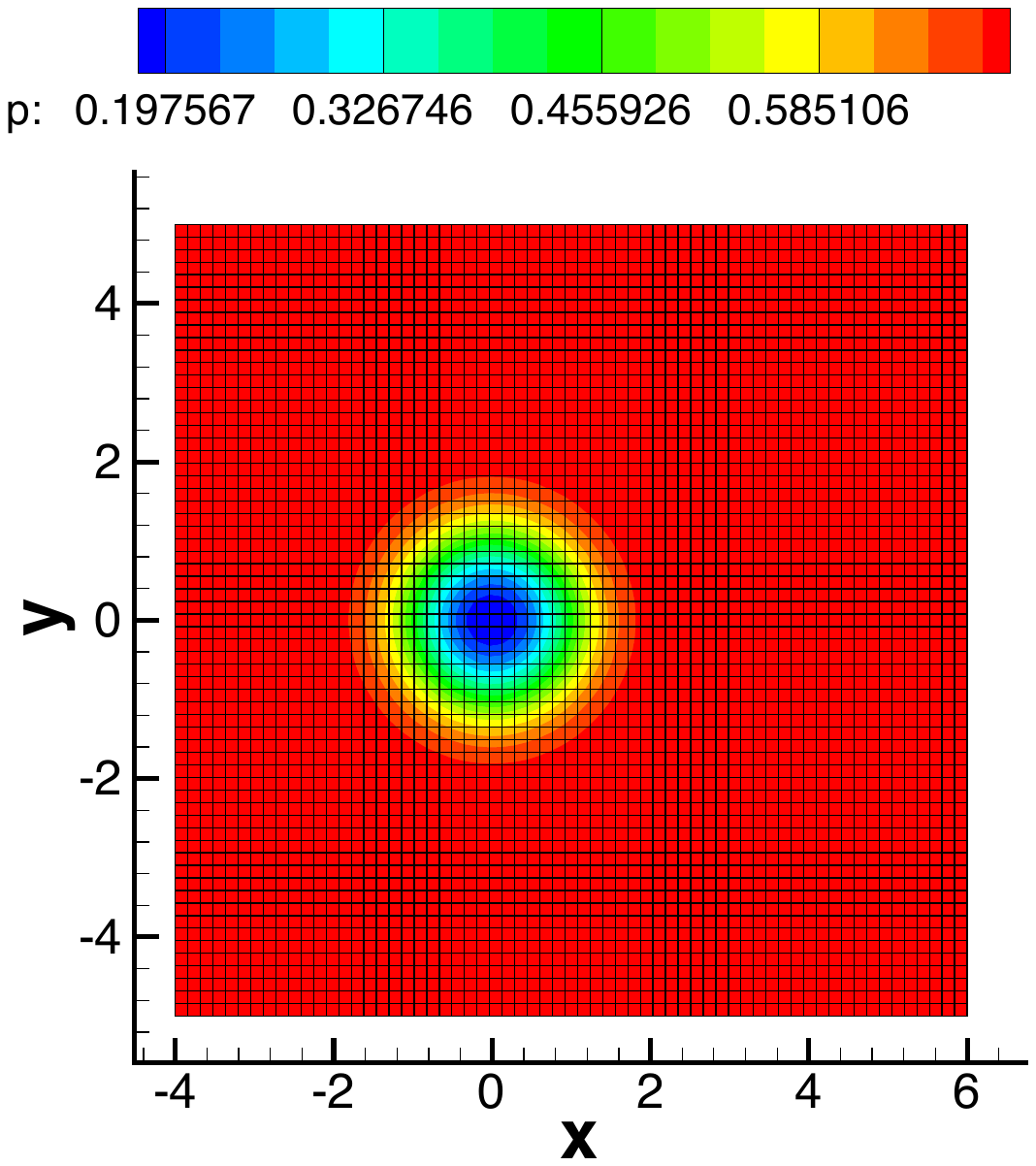}
          \caption{Initial solution on the coarsest grid.}
          \label{fig:inv_uns_err_mms_gridsol}
      \end{subfigure}
      \hfill
      \begin{subfigure}[t]{0.48\textwidth}
        \includegraphics[width=\textwidth,trim=0 0 0 0,clip]{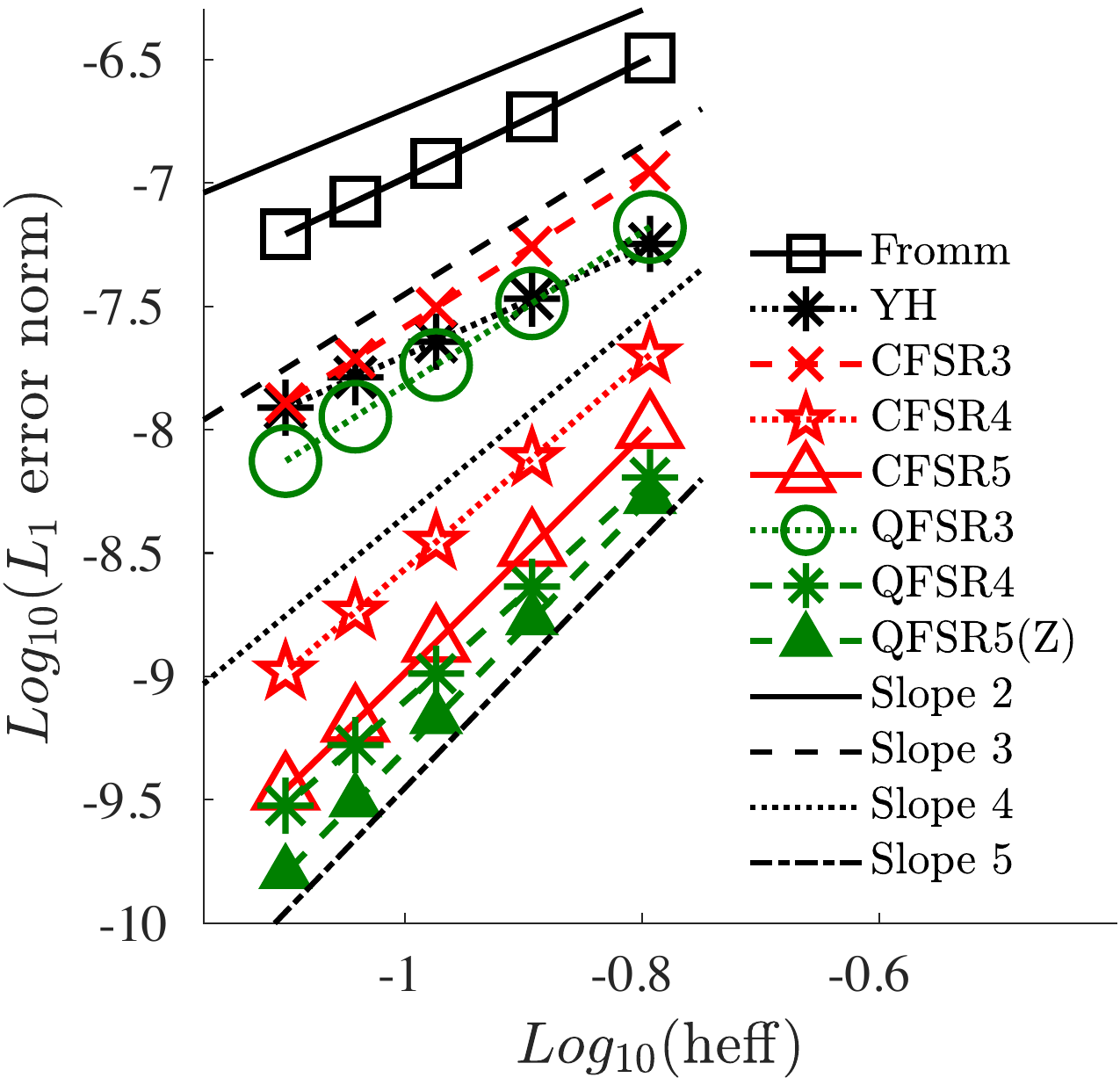}
          \caption{Error convergence.}
          \label{fig:inv_uns_err_mms_error_K5}
      \end{subfigure}
      \caption{Error convergence study for the unsteady Euler equations in two dimensions.}
\label{fig:inv_uns_err_mms} 
\end{figure}

Finally, to demonstrate the design orders of accuracy for an unsteady problem, we consider an inviscid vortex transport problem 
\cite{yang_harris:AIAAJ2016,Barakos:IJNMF2018,ZhongSheng:CF2020,yang_harris:CCP2018,Burg_etal:AIAA2003-3983}, where the exact solution is given by 
\begin{eqnarray}
u =  u_\infty  -  \frac{ K \, \overline{y} }{2 \pi }   \exp \left(  \frac{1-\overline{r} ^2 }{2} \right) ,  \quad
v  = v_\infty   + \frac{ K \, \overline{x} }{2 \pi }   \exp \left(  \frac{1-\overline{r} ^2 }{2} \right),
\end{eqnarray}
and
\begin{eqnarray}
T =  1  -  \frac{ K^2 (\gamma-1)  }{8 \pi^2 }   \exp \left(  1-\overline{r} ^2   \right) , \quad
\rho = T^{  \frac{1}{\gamma-1}  } , \quad
p = \frac{ \rho ^{  \gamma   }  }{\gamma}    ,
\end{eqnarray}
where $\overline{x} = x  - u_\infty t$, $\overline{y} = y  - v_\infty t$, $\overline{r} ^2 =\overline{x}^2 + \overline{y}^2 $, $(u_\infty,v_\infty)=(0.5,0.0)$, and $K=5$ to avoid unexpected linearization as discussed in Ref.\cite{Nishikawa_FakeAccuracy:2020}.  
The initial solution at $t=0$ is shown in Figure \ref{fig:inv_uns_err_mms_gridsol}. For the purpose of accuracy verification, it suffices to perform the calculation for a short time. Therefore, we compute the solution at the final time $t_f=1.0$ with the three-stage SSP Runge-Kutta scheme \cite{SSP:SIAMReview2001} for the total of 1000 time steps with a constant time step $\Delta t = 0.001$, which is small enough for errors to be dominated by the spatial discretization. To verify the spatial order of accuracy, we perform the computation over a series of $n$$\times$$n$ regular quadrilateral grids, 
where $n=64, 80, 96,112, 128$. 
The coarsest grid is shown with pressure contours in Figure \ref{fig:inv_uns_err_mms_gridsol}. We do not consider other types of elements here because effects of element types have already been studied in the previous section and here the focus is on accuracy for an unsteady problem. More specifically, we will demonstrate that the results for an unsteady problem in one dimension will be valid in two dimensions. That is, neither high-order flux quadrature nor a mass matrix formulation is necessary to achieve high-order accuracy on regular quadrilateral grids.

Error convergence results are shown in Figure \ref{fig:inv_uns_err_mms_error_K5}. As expected, Fromm's scheme and YH are second-order accurate but 
YH is significantly more accurate than Fromm's scheme; YH gives errors comparable to those of third-order schemes, CFSR3 and QFSR3, at least for grids tested. For highly refined grids, the true third-order schemes are expected to give much lower errors for a given grid. CFSR4 and CFSR5 are both fourth-order accurate, but CFSR5 is more accurate as we have seen in other tests. QFSR4 and QFSR5(Z) are fourth- and fifth-order accurate as expected. The fourth- and fifth-order schemes are significantly more accurate than Fromm's scheme and YH even on these coarse grids.

{\color{black}
Finally, CPU time required to perform a full time step (i.e., three stages of SSP Runge-Kutta scheme) is measured as the average over 1000 time steps 
 for all the schemes in the finest grid  and compared in Table \ref{Tab.CPUtime_InviscidVortex}.
First, comparing Fromm and YH, we see that the cost of computing the second derivatives and adding a curvature term to the reconstruction 
is very small compared with the total cost per time step: only 7 \% increase, which is comparable with 10 \% increase as reported in 
Ref.\cite{yang_harris:AIAAJ2016}. Second, the CFSR schemes are only slightly more
expensive than Fromm and YH, indicating the cost of the flux reconstruction is also considered as small compared with the total cost per time step.
Finally, the QFSR schemes are more expensive than others but the cost increase is only around 20 \% over other schemes. Note that these results are 
obtained with the explicit time stepping scheme, no special code optimization was performed, and no limiter functions are used. More importantly, a larger 
per-time-step cost does not necessarily indicate a disadvantage if it comes with lower levels of errors. More detailed studies need to be done especially 
for a three-dimensional solver with an implicit time-stepping scheme for practical simulations. 
}
 
%

\begin{table}[htbp!]
\ra{1.0}
\begin{center}
\begin{tabular}{lcccccccc}\hline\hline 
\multicolumn{1}{c}{  }                                               &
\multicolumn{1}{c}{{Fromm} }                                               &
\multicolumn{1}{r}{{YH} }    &
\multicolumn{1}{r}{{CFSR3} }    &
\multicolumn{1}{r}{{CFSR4} }    &
\multicolumn{1}{r}{{CFSR5} }    &
\multicolumn{1}{r}{{QFSR3} }    &
\multicolumn{1}{r}{{QFSR4} }   &
\multicolumn{1}{r}{{QFSR5(Z)} }    
\\ \hline 
\\ [-0.8em]
 Averaged CPU time per time step  & 4.14   & 4.17 &       4.18 & 4.21 &  4.22        & 4.94 & 4.99 & 5.07 
\\ [-0.8em]  \\  \hline  \hline  
\end{tabular}
\caption{ 
{\color{black} Averaged CPU time per time step in seconds measured in the inviscid vortex transport problem on the finest grid.}
}
\label{Tab.CPUtime_InviscidVortex}
\end{center}
\end{table}


\section{Conclusions}
\label{conclusions}

We have provided the clarification of various reconstruction-based unstructured-grid discretization approaches, 
which have been often confused and some of the most economical high-order schemes, which achieve high-order
accuracy with a single flux evaluation per face on regular grids, have gone unnoticed. The key to discovering such schemes is to realize that 
any such scheme must be finite-difference, not finite-volume. Then, we immediately notice that the flux must be reconstructed in order
to achieve high-order accuracy for nonlinear conservation laws. This consideration has led us to the development of a  
general flux-solution-reconstruction schemes (FSR), where reconstructed fluxes and solutions are used to evaluate the consistent part (the averaged flux part)
and the dissipative part of a numerical flux, respectively. 
The FSR schemes can achieve third-, fourth-, and fifth-order accuracy (and fourth- and sixth-order accuracy with zero dissipation) on regular grids.
However, the direct flux reconstruction can be very expensive for unstructured grids in multi-dimensions. 
To develop economical FSR schemes, we have introduced a chain-ruled-based flux reconstruction (CFSR) and a quadratic 
flux reconstruction in terms of solution (QFSR). These economical versions can preserve third- and fourth-order accuracy, but 
fifth-order accuracy is difficult to preserve except cases where the flux is a quadratic in the solution variables. 
For the Euler equations, we have demonstrated that fifth-order accuracy can be achieved by using the parameter vector variables. 
Focusing on the formal accuracy verification, in this paper, we have verified the accuracy of these FSR schemes for one- and two-dimensional
problems. The economical flux reconstruction techniques proposed in this paper may find their wider application, e.g., in generalized finite-difference schemes \cite{HighOrder_GFD:JSC2020} and flux reconstruction schemes \cite{Huynh_FR_adv:AIAA2007,Huynh_FR:AIAA2009}.

Future developments will focus on three areas: boundary effects, shock capturing, and application to practical three-dimensional turbulent-flow 
problems. First, unlike finite-difference schemes, unstructured-grid schemes do not have access to nodes/cells along any possible grid line and 
therefore it is generally difficult to construct high-order discretizations near a boundary. While similar schemes (even though they are not even 
high-order accurate) have been shown to improve resolution of 
a flow solution without high-order boundary treatments \cite{burg_umuscl:AIAA2005-4999,yang_harris:AIAAJ2016,Barakos:IJNMF2018}, 
it would be worth investigating wether high-order accuracy can be achieved with the local grid data available at a boundary.
Second, shock capturing mechanisms need to be incorporated in the FSR schemes for them to be useful for flows with discontinuities. If a limiter is 
applied to the solution reconstruction, it should be applied to the flux reconstruction as well, possibly except for the QFSR schemes, in 
which the flux reconstruction will be automatically limited when the solution reconstructed is limited. Details of shock-capturing 
FSR schemes will be investigated in the future work. Finally, we plan to implement the economical FSR schemes in practical unstructured-grid 
solvers and investigate them for realistic turbulent-flow problems over complex geometries.

\addcontentsline{toc}{section}{Acknowledgments}
\section*{Acknowledgments}

The author gratefully acknowledges support from Software CRADLE, part of Hexagon, the U.S. Army Research Office 
under the contract/grant number W911NF-19-1-0429 with Dr. Matthew Munson as the program manager, and 
 the Hypersonic Technology Project, through the Hypersonic Airbreathing Propulsion Branch of the NASA Langley
 Research Center, under Contract No. 80LARC17C0004. 
 
\addcontentsline{toc}{section}{Data Availability Statement}
\section*{Data Availability Statement}

Data sharing not applicable to this article as no datasets were generated or analysed during the current study.

\addcontentsline{toc}{section}{References}
\bibliography{./bibtex_nishikawa_database}

\begin{thebibliography}{10}
\newcommand{\enquote}[1]{``#1''}

\bibitem{Nishikawa_3rdMUSCL:2020IJNMF}
Nishikawa, H., \enquote{A Truncation Error Analysis of Third-Order {MUSCL}
  Scheme for Nonlinear Conservation Laws,} {\em Int. J. Numer. Meth. Fluids\/},
  Vol.~93, April 2021, pp.~1031--1052.

\bibitem{Nishikawa_3rdQUICK:2020}
Nishikawa, H., \enquote{The {QUICK} Scheme is a Third-Order Finite-Volume
  Scheme with Point-Valued Numerical Solutions,} Vol.~93, April 2021, pp.
  2311--2388.

\bibitem{Nishikawa_FakeAccuracy:2020}
Nishikawa, H., \enquote{On False Accuracy Verification of {UMUSCL} Scheme,} \it
  Commun. Compt. Phys., 2021, in press.

\bibitem{burg_umuscl:AIAA2005-4999}
Burg, C. O.~E., \enquote{Higher Order Variable Extrapolation for Unstructured
  Finite Volume {RANS} Flow Solvers,} {AIAA} Paper 2005-4999, 2005.

\bibitem{yang_harris:AIAAJ2016}
Yang, H.~Q. and Harris, R.~E., \enquote{Development of Vertex-Centered
  High-Order Schemes and Implementation in FUN3D,} {\em {AIAA} J.\/}, Vol.~54,
  2016, pp.~3742--3760.

\bibitem{Barakos:IJNMF2018}
Jimenez-{G}arcia, A. and Barakos, G.~N., \enquote{Assessment of a High-Order
  {MUSCL} Method for Rotor Flows,} {\em Int. J. Numer. Meth. Fluids\/},
  Vol.~87, 2018, pp.~292--327.

\bibitem{ZhongSheng:CF2020}
Zhong, D. and Sheng, C., \enquote{A New Method Towards High-Order Weno Schemes
  on Structured and Unstructured Grids,} {\em Comput. Fluids\/}, Vol.~200,
  2020, pp.~104453.

\bibitem{WhiteNishikawaBaurle_scitech2020}
White, J., Nishikawa, H., and Baurle, R., \enquote{A 3-{D} Nodal-Averaged
  Gradient Approach for Unstructured-grid Cell-centered Finite-volume Methods
  for Application to Turbulent Hypersonic Flow,} {\em SciTech 2020 Forum\/},
  {AIAA} Paper 2020-0652, Orlando, FL, 2020.

\bibitem{scFLOW:Aviation2020}
Higo, Y., Nakashima, Y., Fujiyama, K., Irie, T., and Nishikawa, H.,
  \enquote{{RANS} Solutions on Three-Dimensional Benchmark Configurations with
  {scFLOW}, a Polyhedral Finite-Volume Solver,} {\em {AIAA} Aviation 2020
  Forum\/}, {AIAA} Paper 2020-3029, 2020.

\bibitem{fun3d_website}
\enquote{{FUN3D} online manual,} \url{http://fun3d.larc.nasa.gov}.

\bibitem{nishikawa_liu_source_quadrature:jcp2017}
Nishikawa, H. and Liu, Y., \enquote{Accuracy-Preserving Source Term Quadrature
  for Third-Order Edge-Based Discretization,} {\em J. Comput. Phys.\/},
  Vol.~344, 2017, pp.~595--622.

\bibitem{katz_sankaran:JSC_DOI}
Katz, A. and Sankaran, V., \enquote{An Efficient Correction Method to Obtain a
  Formally Third-Order Accurate Flow Solver for Node-Centered Unstructured
  Grids,} {\em J. Sci. Comput.\/}, Vol.~51, 2012, pp.~375--393.

\bibitem{Boris_Jim_NIA2007-08}
Diskin, B. and Thomas, J.~L., \enquote{Accuracy Analysis for Mixed-Element
  Finite-Volume Discretization Schemes,} {\em NIA Report No. 2007-08\/}, 2007.

\bibitem{ZhangZhangShu2011}
Zhang, R., Zhang, M., and Shu, C.-W., \enquote{On the Order of Accuracy and
  Numerical Performance of Two Classes of Finite Volume {WENO} Schemes,} {\em
  Commun. Comput. Phys.\/}, Vol.~9, No.~3, 2011, pp.~807--827.

\bibitem{NLV6_INRIA_report:2008}
Koobus, B., Wornom, S., Camarri, S., Salvetti, M.-V., and Dervieux, A.,
  \enquote{NonLinear {V6} Schemes for Compressible Flow,} INRIA-00224120v2,
  2008.

\bibitem{yang_harris:CCP2018}
Yang, H.~Q. and Harris, R.~E., \enquote{High-Order Vertex-Centered {U-MUSCL}
  Schemes for Turbulent Flows,} {\em Commun. Comput. Phys.\/}, Vol.~24, No.~2,
  2018, pp.~356--382.

\bibitem{DementRuffin:aiaa2018-1305}
Dement, D.~C. and Ruffin, S.~M., \enquote{Higher Order Cell Centered Finite
  Volume Schemes for Unstructured Cartesian Grids,} {\em 56th {AIAA} Aerospace
  Sciences Meeting\/}, {AIAA} Paper 2018-1305, Kissimmee, Florida, 2018.

\bibitem{jalali_gooch:CF2017}
Jalali, A. and Ollivier-Gooch, C., \enquote{Higher-Order Unstructured Finite
  Volume {RANS} Solution of Turbulent Compressible Flows,} {\em Comput.
  Fluids\/}, Vol.~143, 2017, pp.~32--47.

\bibitem{Tsoutsanis-HOFV:JCP2018}
Tsoutsanis, P., \enquote{Extended Bounds Limiter for High-Order Finite-Volume
  Schemes on Unstructured Meshes,} {\em J. Comput. Phys.\/}, Vol.~362, 2018,
  pp.~69--94.

\bibitem{FVWENO:JSC2014}
Buchm\"{u}ller, P. and Helzel, C., \enquote{Improved Accuracy of High-Order
  {WENO} Finite Volume Methods on {Cartesian} Grids,} {\em J. Sci. Comput.\/},
  Vol.~61, 2014, pp.~343--368.

\bibitem{FVWENO:AMC2016}
Buchm\"{u}ller, P., Drehe, J., and Helzel, C., \enquote{Finite Volume {WENO}
  Methods for Hyperbolic Conservation Laws on {Cartesian} Grids with Adaptive
  Mesh Refinement,} {\em Applied Mathematics and Computation\/}, Vol.~272,
  2016, pp.~460--478.

\bibitem{TamakiImamura:CF2017}
Tamaki, Y. and Imamura, T., \enquote{Efficient Dimension-by-Dimension Higher
  Order Finite-Volume Methods for a {Cartesian} grid with Cell-Based
  Refinement,} {\em Comput. Fluids\/}, Vol.~144, 2017, pp.~74--85.

\bibitem{Tamaki_PhDThesis2018}
Tamaki, Y., {\em Turbulent Flow Simulations around Aircraft using Hierarchical
  Cartesian Grids and the Immersed Boundary Method\/}, Ph.D. thesis, University
  of Tokyo, March 2018.

\bibitem{Denaro:IJNMF1996}
Denaro, F.~M., \enquote{Towards a New Model-Free Simulation of
  High-{R}eynolds-Flows: Local Average Direct Numerical,} {\em Int. J. Numer.
  Meth. Fluids\/}, Vol.~23, 1996, pp.~125--142.

\bibitem{Denaro:IJNMF2002}
Denaro, F.~M. and Sarghini, F., \enquote{2-{D} Transmitral Flows Simulations by
  Means of the Immersed Boundary Method on Unstructured Grids,} {\em Int. J.
  Numer. Meth. Fluids\/}, Vol.~38, 2002, pp.~1133--1157.

\bibitem{FeliceDenaroMoela:NHT1993}
de~Felice, D. and Moela, F. M. D.~C., \enquote{Multidimensional Single-Step
  Vector Upwind Schemes for Highly Convective Transport Problems,} {\em
  Numerical Heat Transfer, {P}art {B}: Fundamentals: An International Journal
  of Computation and Methodology\/}, Vol.~23, No.~4, 1993, pp.~425--460.

\bibitem{nishikawa:AIAA2010}
Nishikawa, H., \enquote{Beyond Interface Gradient: A General Principle for
  Constructing Diffusion Schemes,} {\em Proc. of 40th {AIAA} Fluid Dynamics
  Conference and Exhibit\/}, {AIAA} Paper 2010-5093, Chicago, 2010.

\bibitem{DiskinThomas:ANM2010}
Diskin, B. and Thomas, J.~L., \enquote{Notes on Accuracy of Finite-Volume
  Discretization Schemes on Irregular Grids,} {\em Appl. Numer. Math.\/},
  Vol.~60, 2010, pp.~224--226.

\bibitem{diskin_thomas:AIAA2012-0609}
Diskin, B. and Thomas, J.~L., \enquote{Effects of Mesh Regularity on Accuracy
  of Finite-Volume Schemes,} {\em Proc. of 50th AIAA Aerospace Sciences
  Meeting\/}, {AIAA} Paper 2012-0609, Nashville, Tennessee, 2012.

\bibitem{Hong_Meshless:JCP2007}
Luo, H., Baum, J.~D., and L\"{o}hner, R., \enquote{A Hybrid {C}artesian Grid
  and Gridless Method for Compressible Flows,} {\em J. Comput. Phys.\/},
  Vol.~214, 2006, pp.~618--632.

\bibitem{aiaa2009-596}
Katz, A. and Jameson, A., \enquote{A Comparison of Various Meshless Schemes
  Within a Unified Algorithm,} {\em Proc. of 47th AIAA Aerospace Sciences
  Meeting\/}, {AIAA} Paper 2009-596, Orlando, FL, 2009.

\bibitem{Chiu_etal_meshless_SIAM_2014}
Chiu, E.~K., Wang, Q., Hu, R., and Jameson, A., \enquote{A Conservative
  Mesh-Free Scheme and Generalized Framework for Conservation Laws,} {\em SIAM
  J. Sci. Comput.\/}, Vol.~36, 2014, pp.~A2896--A2916.

\bibitem{HighOrder_GFD:JSC2020}
Li, X.-L. and Ren, Y.-X., \enquote{High Order Compact Generalized Finite
  Difference Methods for Solving Inviscid Compressible Flows,} {\em J. Sci.
  Comput.\/}, Vol.~82, 2020, pp.~18.

\bibitem{Nishikawa_aiaa2020-1786}
Nishikawa, H., \enquote{A Face-Area-Weighted Centroid Formula for Reducing Grid
  Skewness and Improving Convergence of Edge-Based Solver on Highly-Skewed
  Simplex Grids,} {\em AIAA Scitech 2020 Forum\/}, {AIAA} Paper 2020-1786,
  Orlando, FL, 2020.

\bibitem{nishikawa_centroid:JCP2020}
Nishikawa, H., \enquote{A Face-Area-Weighted Centroid Formula for Finite-Volume
  Method That Improves Skewness and Convergence on Triangular Grids,} {\em J.
  Comput. Phys.\/}, Vol.~401, 2020, pp.~109001.

\bibitem{liu_nishikawa_aiaa2017-0738}
Liu, Y. and Nishikawa, H., \enquote{Third-Order Inviscid and Second-Order
  Hyperbolic {N}avier-{S}tokes Solvers for Three-Dimensional Unsteady Inviscid
  and Viscous Flows,} {\em 55th {AIAA} Aerospace Sciences Meeting\/}, {AIAA}
  Paper 2017-0738, Grapevine, Texas, 2017.

\bibitem{nishikawa_liu_aiaa2018-4166}
Nishikawa, H. and Liu, Y., \enquote{Third-Order Edge-Based Scheme for Unsteady
  Problems,} {\em AIAA 2018 Fluid Dynamics Conference\/}, {AIAA} Paper
  2018-4166, Atlanta, Georgia, 2018.

\bibitem{liu_nishikawa_aiaa2016-3969}
Liu, Y. and Nishikawa, H., \enquote{Third-Order Inviscid and Second-Order
  Hyperbolic {N}avier-{S}tokes Solvers for Three-Dimensional Inviscid and
  Viscous Flows,} {\em 46th {AIAA} Fluid Dynamics Conference\/}, {AIAA} Paper
  2016-3969, Washington, D.C., 2016.

\bibitem{nishikawa_boundary_formula:JCP2014}
Nishikawa, H., \enquote{Accuracy-Preserving Boundary Flux Quadrature for
  Finite-Volume Discretization on Unstructured Grids,} 2014, in review.

\bibitem{NishikawaPadway:Aviation2020}
Nishikawa, H. and Padway, E., \enquote{An Adaptive Space-Time Edge-Based Solver
  for Two-Dimensional Unsteady Inviscid Flows,} {\em {AIAA} Aviation 2020
  Forum\/}, {AIAA} Paper 2020-3024, 2020.

\bibitem{LiuNishikawa_2017-3443}
Liu, Y. and Nishikawa, H., \enquote{Third-Order Edge-Based Hyperbolic
  {N}avier-{S}tokes Scheme for Three-Dimensional Viscous Flows,} {\em 23rd
  {AIAA} Computational Fluid Dynamics Conference\/}, {AIAA} Paper 2017-3443,
  Denver, Colorado, 2017.

\bibitem{Merryman:JSC2003}
Merriman, B., \enquote{Understanding the {S}hu-{O}sher Conservative Finite
  Difference Form,} {\em J. Sci. Comput.\/}, Vol.~19, 2003, pp.~309--322.

\bibitem{VLeer_Ultimate_III:JCP1977}
{van Leer}, B., \enquote{Towards the Ultimate Conservative Difference Scheme.
  {III}. Upstream-centered Finite Difference Schemes for Ideal Compressible
  Flow,} {\em J. Comput. Phys.\/}, Vol.~23, 1977, pp.~263--275.

\bibitem{Shu_Osher_Efficient_ENO_II_JCP1989}
Shu, C.-W. and Osher, S.~J., \enquote{Efficient Implementation of Essentially
  Non-Oscillatory Shock-Capturing Schemes, {II},} {\em J. Comput. Phys.\/},
  Vol.~83, 1989, pp.~32--78.

\bibitem{AbalakinBakhvalovKozubskaya:IJNMF2015}
Abalakin, I., Bakhvalov, P., and Kozubskaya, T., \enquote{Edge-Based
  Reconstruction Schemes for Unstructured Tetrahedral Meshes,} {\em Int. J.
  Numer. Meth. Fluids\/}, Vol.~81, 2015, pp.~331--356.

\bibitem{Nishikawa_FANG_AQ:Aviation2020}
Nishikawa, H., \enquote{A Face-Averaged Nodal Gradient Cell-Centered
  Finite-Volume Method for Mixed Grids,} {\em {AIAA} Aviation 2020 Forum\/},
  {AIAA} Paper 2020-3049, 2020.

\bibitem{Roe_JCP_1981}
Roe, P.~L., \enquote{Approximate {R}iemann Solvers, Parameter Vectors, and
  Difference Schemes,} {\em J. Comput. Phys.\/}, Vol.~43, 1981, pp.~357--372.

\bibitem{Nishikawa_RobustFluxes:jcp2020}
Nishikawa, H., \enquote{Robust Numerical Fluxes for Unrealizable States,} {\em
  J. Comput. Phys.\/}, Vol.~408, 2020, pp.~109244.

\bibitem{nishikawa_LP_UMUSCL:JCP2020}
Nishikawa, H., \enquote{On the Loss and Recovery of Second-Order Accuracy with
  {U-MUSCL},} {\em J. Comput. Phys.\/}, Vol.~417, 2020, pp.~109600.

\bibitem{nishikawa_stencil:JCP2019}
Nishikawa, H., \enquote{Efficient Gradient Stencils for Robust Implicit
  Finite-Volume Solver Convergence on Distorted Grids,} {\em J. Comput.
  Phys.\/}, Vol.~386, 2019, pp.~486--501.

\bibitem{idolikeCFD_VOL1_v2p6_pdf}
Masatsuka, K., \enquote{I do like {CFD}, {VOL}.1, {S}econd {E}dition, Version
  2.6,} \url{http://www.cfdbooks.com}, 2018.

\bibitem{Multigrid_book_2001}
Trottenberg, U., Oosterlee, C.~W., and Sch\"uller, A., {\em Multigrid\/},
  Academic Press, 2000.

\bibitem{nishikawa_hyperbolic_poisson:jcp2020}
Nishikawa, H., \enquote{A hyperbolic {P}oisson Solver for Tetrahedral Grids,}
  {\em J. Comput. Phys.\/}, Vol.~409, May 2020, pp.~109358.

\bibitem{nishikawa_liu_jcp2018}
Nishikawa, H. and Liu, Y., \enquote{Hyperbolic Advection-Diffusion Schemes for
  High-{R}eynolds-Number Boundary-Layer Problems,} {\em J. Comput. Phys.\/},
  Vol.~352, 2018, pp.~23--51.

\bibitem{SSP:SIAMReview2001}
Gottlieb, S., Shu, C.-W., and Tadmor, E., \enquote{Strong Stability-Preserving
  High-Order Time Discretization Methods,} {\em SIAM Rev.\/}, Vol.~43, No.~1,
  2001, pp.~89--112.

\bibitem{Burg_etal:AIAA2003-3983}
Burg, O.~E., Sheng, C., Newman, J.~C., Brewer, W., Blades, E., and Marcum,
  D.~L., \enquote{Verification and Validation of Forces Generalized by an
  Unstructured Flow Solver,} {\em Proc. of 16th {AIAA} Computational Fluid
  Dynamics Conference\/}, {AIAA} Paper 2003-3983, Orlando, Florida, 2003.

\bibitem{Huynh_FR_adv:AIAA2007}
{O}llivier {G}ooch, C., Nejat, A., and Michalak, K., \enquote{A Flux
  Reconstruction Approach to High-Order Schemes Including Discontinuous
  {G}alerkin Methods,} {\em Proc. of 18th {AIAA} Computational Fluid Dynamics
  Conference\/}, {AIAA} Paper 2007-4079, Miami, 2007.

\bibitem{Huynh_FR:AIAA2009}
Huynh, H.~T., \enquote{A Reconstruction Approach to High-Order Schemes
  Including Discontinuous {G}alerkin for Diffusion,} {\em Proc. of 47th AIAA
  Aerospace Sciences Meeting\/}, {AIAA} Paper 2009-403, Orlando, Florida, 2009.

\end{thebibliography}
\bibliographystyle{aiaa}

 
\end{document}